\documentclass[a4paper,10pt]{article}
\usepackage[utf8]{inputenc}
\usepackage[left=4cm, right=3cm, bottom=3.1cm, top=3.1cm]{geometry}

\usepackage{amsmath, amssymb, amsthm}
\usepackage{mathtools}
\usepackage{mathrsfs}
\usepackage{scalerel}

\usepackage{enumitem}
\usepackage[all]{xy}
\usepackage{url}
\usepackage{xcolor}
\usepackage{mdframed}
\usepackage{listings}
\usepackage{bbm}
\usepackage{hyperref}
\usepackage{tikz-cd}
\usetikzlibrary{decorations.pathmorphing}
\tikzcdset{arrow style = tikz,diagrams={>={Classical TikZ Rightarrow}}}

\usepackage{stmaryrd}

\usepackage{marginnote}
\usepackage{comment}

\usepackage{letltxmacro}

\numberwithin{equation}{subsection}

\theoremstyle{definition}

\newtheorem{defn}{Definition}[section]
\newtheorem{ex}[defn]{Example}
\newtheorem{non-ex}[defn]{Non-Example}
\newtheorem{new-unob-ex}[defn]{New Unobstructed Example}
\newtheorem{sitn}[defn]{Situation}
\newtheorem{var-def}[defn]{Variant}
\newtheorem{ass}[defn]{Assumption}

\newtheorem{constr}[defn]{Construction}
\newtheorem{rem}[defn]{Remark}

\newtheorem{dat}[defn]{Datum}
\newtheorem{strat}[defn]{Strategy}

\theoremstyle{plain}

\newtheorem{prop}[defn]{Proposition}
\newtheorem{cor}[defn]{Corollary}
\newtheorem{lem}[defn]{Lemma}
\newtheorem{univ-prop}[defn]{Universal Property}
\newtheorem{thm}[defn]{Theorem}
\newtheorem{lem-def}[defn]{Lemma and Definition}

\newtheorem{quest}[defn]{Question}

\newtheorem{obs}[defn]{Observation}

\theoremstyle{remark}

\newenvironment{enum:arabic}{\enumerate[label=\textnormal{(\arabic*)}]}{\endenumerate}
\newenvironment{enum:alph}{\enumerate[label=\textnormal{(\alph*)}]}{\endenumerate}
\newenvironment{enum:roman}{\enumerate[label=\textnormal{(\roman*)}]}{\endenumerate}

\DeclareFontFamily{U}  {MnSymbolC}{}
\DeclareSymbolFont{MnSyC}         {U}  {MnSymbolC}{m}{n}
\SetSymbolFont{MnSyC}       {bold}{U}  {MnSymbolC}{b}{n}
\DeclareFontShape{U}{MnSymbolC}{m}{n}{
	<-6>  MnSymbolC5
	<6-7>  MnSymbolC6
	<7-8>  MnSymbolC7
	<8-9>  MnSymbolC8
	<9-10> MnSymbolC9
	<10-12> MnSymbolC10
	<12->   MnSymbolC12}{}
\DeclareFontShape{U}{MnSymbolC}{b}{n}{
	<-6>  MnSymbolC-Bold5
	<6-7>  MnSymbolC-Bold6
	<7-8>  MnSymbolC-Bold7
	<8-9>  MnSymbolC-Bold8
	<9-10> MnSymbolC-Bold9
	<10-12> MnSymbolC-Bold10
	<12->   MnSymbolC-Bold12}{}

\DeclareMathSymbol{\invneg}{\mathord}{MnSyC}{181}
\DeclareMathSymbol{\blackdiamond}{\mathord}{MnSyC}{109}
\DeclareMathSymbol{\largestar}{\mathord}{MnSyC}{131}
\DeclareMathSymbol{\splus}{\mathord}{MnSyC}{20}
\DeclareMathSymbol{\stimes}{\mathord}{MnSyC}{21}

\gdef\notesoff{\gdef\note##1{}}
\notesoff

\DeclareFontFamily{U}{mathx}{}
\DeclareFontShape{U}{mathx}{m}{n}{ <-> mathx10 }{}
\DeclareSymbolFont{mathx}{U}{mathx}{m}{n}
\DeclareFontSubstitution{U}{mathx}{m}{n}
\DeclareMathAccent{\widecheck}{0}{mathx}{"71}

\DeclareSymbolFont{itgreek}{OML}{cmm}{m}{it}

\DeclareMathSymbol{\itGamma}  {\mathalpha}{itgreek}{"00}
\DeclareMathSymbol{\itDelta}  {\mathalpha}{itgreek}{"01}
\DeclareMathSymbol{\itTheta}  {\mathalpha}{itgreek}{"02}
\DeclareMathSymbol{\itLambda} {\mathalpha}{itgreek}{"03}
\DeclareMathSymbol{\itXi}     {\mathalpha}{itgreek}{"04}
\DeclareMathSymbol{\itPi}     {\mathalpha}{itgreek}{"05}
\DeclareMathSymbol{\itSigma}  {\mathalpha}{itgreek}{"06}
\DeclareMathSymbol{\itUpsilon}{\mathalpha}{itgreek}{"07}
\DeclareMathSymbol{\itPhi}    {\mathalpha}{itgreek}{"08}
\DeclareMathSymbol{\itPsi}    {\mathalpha}{itgreek}{"09}
\DeclareMathSymbol{\itOmega}  {\mathalpha}{itgreek}{"0A}

\newcommand{\kk}{\Bbbk}
\newcommand{\kq}{\mathbbm{q}}

\newcommand{\IA}{\mathbb{A}}

\newcommand{\CC}{\mathbb{C}}

\newcommand{\FF}{\mathbb{F}}

\newcommand{\HH}{\mathbb{H}}

\newcommand{\NN}{\mathbb{N}}

\newcommand{\PP}{\mathbb{P}}
\newcommand{\QQ}{\mathbb{Q}}
\newcommand{\RR}{\mathbb{R}}
\newcommand{\VV}{\mathbb{V}}

\newcommand{\ZZ}{\mathbb{Z}}

\newcommand{\cA}{\mathcal{A}}

\newcommand{\cC}{\mathcal{C}}

\newcommand{\cE}{\mathcal{E}}
\newcommand{\cF}{\mathcal{F}}
\newcommand{\cG}{\mathcal{G}}
\newcommand{\cH}{\mathcal{H}}
\newcommand{\cI}{\mathcal{I}}

\newcommand{\cK}{\mathcal{K}}
\newcommand{\cL}{\mathcal{L}}
\newcommand{\cM}{\mathcal{M}}
\newcommand{\cN}{\mathcal{N}}
\newcommand{\cO}{\mathcal{O}}
\newcommand{\cP}{\mathcal{P}}
\newcommand{\cQ}{\mathcal{Q}}

\newcommand{\cS}{\mathcal{S}}
\newcommand{\cT}{\mathcal{T}}
\newcommand{\cU}{\mathcal{U}}
\newcommand{\cV}{\mathcal{V}}
\newcommand{\cW}{\mathcal{W}}

\newcommand{\scA}{\mathscr{A}}

\newcommand{\scC}{\mathscr{C}}
\newcommand{\scD}{\mathscr{D}}
\newcommand{\scE}{\mathscr{E}}

\newcommand{\scH}{\mathscr{H}}
\newcommand{\scI}{\mathscr{I}}

\newcommand{\scK}{\mathscr{K}}

\newcommand{\scS}{\mathscr{S}}
\newcommand{\scT}{\mathscr{T}}

\newcommand{\scV}{\mathscr{V}}
\newcommand{\scW}{\mathscr{W}}
\newcommand{\scX}{\mathscr{X}}

\newcommand{\scZ}{\mathscr{Z}}

\newcommand{\fX}{\mathfrak{X}}

\newcommand{\bA}{\mathbf{A}}
\newcommand{\bB}{\mathbf{B}}

\DeclareMathOperator{\Spec}{Spec}
\DeclareMathOperator{\Spf}{Spf}

\DeclareMathOperator{\Proj}{Proj}

\newcommand{\eps}{\varepsilon}
\newcommand{\cp}{\varphi}

\newcommand{\lookup}{\colorbox{yellow}{Look it up!}}

\newcommand{\m}[1]{\mathrm{#1}}

\let\lsem\llbracket
\let\rsem\rrbracket

\newcommand*{\qeddef}{%
	\leavevmode\unskip\penalty9999 \hbox{}\nobreak\hfill
	\quad\hbox{$\blacklozenge$}%
}

\newcommand*{\qedex}{%
	\leavevmode\unskip\penalty9999 \hbox{}\nobreak\hfill
	\quad\hbox{$\lozenge$}%
}

\DeclareRobustCommand{\punct}[0]{{\scaleobj{0.85}{\odot}}}

\DeclareFontFamily{U}{mathx}{\hyphenchar\font45}
\DeclareFontShape{U}{mathx}{m}{n}{
	<-6> mathx5 <6-7> mathx6 <7-8> mathx7
	<8-9> mathx8 <9-10> mathx9
	<10-12> mathx10 <12-> mathx12
}{}
\DeclareSymbolFont{mathx}{U}{mathx}{m}{n}
\DeclareFontSubstitution{U}{mathx}{m}{n}
\DeclareMathAccent{\widebar}{0}{mathx}{"73}
\title{Log Hodge--de Rham degeneration \\ beyond the toroidal case}
\author{Simon Felten and Matthias Zach}

\begin{document}

\maketitle

\begin{abstract}
	Let $\phi\colon X \to S$ be a flat and proper one-parameter degeneration of algebraic varieties with a single degenerate fiber $X_0$. If $\phi\colon (X|X_0) \to (S|0)$ is log smooth and saturated, the log Hodge--de Rham spectral sequence is known to degenerate at the page $E_1$. Adding to that some further well-known related properties, we say that the log de Rham complex $\Omega^{\bullet,\m{log}}_{X/S}$ is \emph{cohomologically flat as a filtered complex}, short \emph{cffc}. Furthermore, if $\phi\colon (X|X_0) \to (S|0)$ has only toroidal log singularities, it is known that the \emph{reflexive log de Rham complex} $\cW^\bullet_{X/S}$ is cffc.
        A major motivation to study this property is that in the Calabi--Yau case, $\cW^\bullet_{X/S}$ being cffc implies that a certain logarithmic deformation functor for $\phi_0\colon X_0^\partial \to S_0^\dagger$ is unobstructed. 
	
	In this article, we develop a computational criterion to determine, under very mild assumptions on $\phi\colon X \to S$, whether the reflexive log de Rham complex $\cW^\bullet_{X/S}$ is cffc. Applying this criterion, we found many explicit examples of degenerations of surfaces where $\cW^\bullet_{X/S}$ is cffc, even though the log singularities are worse than toroidal. We did not find any new examples of degenerations of threefolds such that $\cW^\bullet_{X/S}$ is cffc.
	
\end{abstract}

\tableofcontents

\newpage

\section{Introduction}

In this article, we study the cohomology and deformation theory of degenerating families of algebraic varieties of the following kind:

\begin{sitn}\label{sitn:basic}\note{sitn:basic}
	Let $\kk$ be an algebraically closed field with $\m{char}(\kk) = 0$ and let $\widebar S = \IA^1_t = \Spec \kk[t]$. Let $\widebar X \subset \PP^n \times \widebar S$ be a closed subscheme and assume that the induced morphism $\phi\colon \widebar X \to \widebar S$ is flat. Fix $d \geq 1$ and assume that all fibers of $\phi\colon \widebar X \to \widebar S$ are pure of dimension $d$. Let $0 \in S \subseteq \widebar S$ be an open subset, let $S^\star = S \setminus \{0\}$, let $X^\star = \phi^{-1}(S^\star) \subseteq X$, and assume that $\phi\colon X^\star \to S^\star$ is smooth. Now, $\phi\colon X \to S$ is a degenerating family of algebraic varieties. \qeddef
\end{sitn}
\begin{ex}
%	Let $S = \widebar S \setminus \{-\frac{4}{27}\}$ and
	\[\widebar X = \{y^2z - x^3 - x^2z - tz^3 = 0\} \subset \PP^2 \times \widebar S.\]
	Then, $\phi\colon X \to S$ is an instance of Situation~\ref{sitn:basic}. It is a semistable degeneration of elliptic curves. \qedex
\end{ex}

We obtain a morphism of log schemes $\phi\colon X^\dagger \to S^\dagger$ by using the divisorial log structures defined by $\{0\} \subset S$, respectively $X_0 \subset X$. If $\phi\colon X^\dagger \to S^\dagger$ is log smooth or has only toroidal log singularities, it is known that a certain log de Rham complex $\cA^\bullet_{X/S} = \cW^\bullet_{X/S}$ is cohomologically flat as a filtered complex in the sense of Definition~\ref{defn:coh-flat}. As a consequence, if $\phi\colon X^\dagger \to S^\dagger$ is log Calabi--Yau, a certain logarithmic deformation functor $\m{ELD}_{X_0/S_0}^\scD$ for the central fiber $\phi_0\colon X_0^\dagger \to S_0^\dagger$ is unobstructed.

Our main contribution in this article is the experimental observation that $\cA^\bullet_{X/S}$ is cohomologically flat as a filtered complex in many explicit examples of degenerations of \emph{surfaces} with \emph{non-toroidal} log singularities, see Observation~\ref{obs:two-dim-cffc}. In the Calabi--Yau case, this implies again that $\m{ELD}_{X_0/S_0}^\scD$ is unobstructed, see Theorem~\ref{thm:ELD-unobstructed}. We did not find any example of a degeneration of \emph{threefolds} with non-toroidal log singularities and $\cA^\bullet_{X/S}$ being cohomologically flat as a filtered complex. For motivation to study non-toroidal log singularities, see Section~\ref{sect:non-toroidal-log-singularities-intro}, Section~\ref{sect:intrinsic-mirror}, and Section~\ref{sect:Fano-smoothing}.

\subsection{The smooth case}

Suppose that $X_0$ is smooth in Situation~\ref{sitn:basic}. Each piece $\Omega^i_{X/S}$ of the de Rham complex is finite locally free, in particular flat over $S$. The Hodge sheaves $R^j\phi_*\Omega^i_{X/S}$ are finite locally free $\cO_S$-modules, and for every Cartesian diagram 
\[
\xymatrix{
	Y \ar[d]^\psi \ar[r]^c & X \ar[d]^\phi \\
	T \ar[r]^b & S \\
}
\]
of schemes, the natural map $b^*R^j\phi_*\Omega^i_{X/S} \to R^j\psi_*(c^*\Omega^i_{X/S})$ is an isomorphism. Moreover, the Hodge--de Rham spectral sequence
\[R^j\phi_*\Omega^i_{X/S} \Rightarrow R^{i + j}\phi_*\Omega^\bullet_{X/S}\]
degenerates at the page $E_1$. We say that the de Rham complex $\Omega^\bullet_{X/S}$ is \emph{cohomologically flat as a filtered complex}, and we shall abbreviate this as \emph{cffc}. For more details, see Definition~\ref{defn:coh-flat}.

We say that $X_0$ is \emph{Calabi--Yau} if $\Omega^d_{X_0} \cong \cO_{X_0}$. In this case, the celebrated Bogomolov--Tian--Todorov theorem guarantees that the flat deformation functor $\m{Def}_{X_0}$ is unobstructed. Various proofs of this fact are known, see for example  \cite{Bogomolov1978,Tian1987,Todorov1989,Ran1992,BarannikovKontsevich1998,IaconoManetti2010}.

\subsection{The log smooth case}

We continue in Situation~\ref{sitn:basic} but do not assume that $X_0$ is smooth. As a tool to study the degeneration, we use \emph{logarithmic structures}. A log structure on a scheme $Y$ is a sheaf of commutative monoids $\cM_Y$ in the \'etale topology of $Y$, together with a monoid homomorphism $\alpha\colon \cM_Y \to (\cO_Y,\cdot)$ such that $\alpha^{-1}(\cO_Y^\times) \to \cO_Y^\times$ is an isomorphism. For a comprehensive treatment of log structures, we refer the reader to \cite{LoAG2018}, for a quick introduction to \cite[\S7]{FeltenGLDT}.

We endow $X$ with the divisorial logarithmic structure $\cM_X$ defined by $X_0 \subset X$; we write $X^\dagger = (X,\cM_X)$ for the resulting log scheme. Concretely, for an \'etale morphism $h\colon W \to X$, we have 
\[\Gamma(W,\cM_X) = \{f \in \Gamma(W,\cO_W) \ | \ f|_{W \setminus h^{-1}(X_0)} \text{ is invertible}\},\]
together with the obvious map $\alpha\colon \cM_X \to \cO_X$. Similarly, we endow $S$ with the divisorial logarithmic structure $\cM_S$ defined by $\{0\} \subset S$, and we obtain a morphism of log schemes $\phi\colon X^\dagger \to S^\dagger$.

In the best situations, $\phi\colon X^\dagger \to S^\dagger$ has the infinitesimal lifting property in the category of log schemes---it is \emph{log smooth} in the sense of \cite[Defn.~IV.3.1.1]{LoAG2018}. For technical reasons, we also wish $\phi\colon X^\dagger \to S^\dagger$ to be \emph{saturated} in the sense of \cite[Defn.~III.2.5.1]{LoAG2018}.\footnote{When we say that $\phi\colon X^\dagger \to S^\dagger$ is log smooth and saturated at a point $p \in X$, we implicitly understand that the log structures of $X^\dagger$ and $S^\dagger$ are fine and saturated at $p$ respectively $\phi(p)$.} If $\phi\colon X \to S$ is semistable, both are satisfied.

The family $\phi\colon X^\dagger \to S^\dagger$ has a \emph{log de Rham complex} $\Omega^{\bullet,\m{log}}_{X/S}$, see \cite[Prop.~V.2.1.1]{LoAG2018}. We always have $\Omega^{0,\m{log}}_{X/S} = \cO_X$. If $\phi\colon X^\dagger \to S^\dagger$ is log smooth and saturated, $\Omega^{1,\m{log}}_{X/S}$ is a vector bundle of rank $d$, and $\Omega^{\bullet,\m{log}}_{X/S}$ is concentrated in degrees $[0,d]$.

\begin{thm}\label{thm:log-smooth-cffc}\note{thm:log-smooth-cffc}
	In Situation~\ref{sitn:basic}, assume that $\phi\colon X^\dagger \to S^\dagger$ is log smooth and saturated. Then, the log de Rham complex $\Omega^{\bullet,\m{log}}_{X/S}$ is cffc.
\end{thm}
\begin{proof}
	See \cite[Thm.~4.12]{kkatoFI}.
\end{proof}

\subsection{The log smooth Bogomolov--Tian--Todorov theorem}

We continue in Situation~\ref{sitn:basic} and assume that $\phi\colon X^\dagger \to S^\dagger$ is log smooth and saturated. In this section, we discuss the log smooth deformation functor of the central fiber. Its unobstructedness in the Calabi--Yau case can be framed as a consequence of Theorem~\ref{thm:log-smooth-cffc}.

Let $S_0^\dagger = \Spec(\NN \to \kk)$ be the \emph{standard log point}, the spectrum of the pre-log ring $\beta\colon \NN \to \kk$ given by $\beta(0) = 1$ and $\beta(n) = 0$ for $n \geq 1$. Explicitly, it is the scheme $S_0 = \Spec \kk$ together with the log structure $\cM_{S_0}$ given by $\Gamma(S_0,\cM_{S_0}) = \NN \oplus \kk^\times$ and $\alpha(n,\lambda) = 0^n \cdot \lambda$, where $0^0 = 1$ and $0^n = 0$ for $n \geq 1$.

We have a strict closed immersion $b_0\colon S_0^\dagger \to S^\dagger$; it gives rise to the log-geometric central fiber $\phi_0\colon X_0^\dagger \to S_0^\dagger$ of $\phi\colon X^\dagger \to S^\dagger$ as the base change along $b_0\colon S_0^\dagger \to S^\dagger$ in the category of log schemes. The underlying scheme of $X_0^\dagger$ is $X_0 = \phi^{-1}(0)$.

Let $\m{Art}_{\kk\lsem t\rsem\sslash\kk}$ be the category whose objects are Artinian local $\kk\lsem t\rsem$-algebras $A$ with residue field $\kk$. For the exact definition, see \cite[Defn.~A.1]{FeltenGLDT}. For each $A \in \m{Art}_{\kk\lsem t\rsem\sslash\kk}$, let $S_A^\dagger = \Spec(\NN \to A)$ be the spectrum of the pre-log ring $\beta\colon \NN \to A$ given by $\beta(n) = t^n$; its underlying scheme is $S_A = \Spec A$. Then, $S_0^\dagger \to S_A^\dagger$ is an infinitesimal thickening of log schemes.

The \emph{log smooth deformation functor}
\[\m{LD}_{X_0/S_0}\colon \enspace \m{Art}_{\kk\lsem t\rsem\sslash\kk} \to \m{Set}\]
classifies, for each $A \in \m{Art}_{\kk\lsem t\rsem\sslash\kk}$, log smooth thickenings $\psi_A\colon Y_A^\dagger \to S_A^\dagger$ of $\phi_0\colon X_0^\dagger \to S_0^\dagger$, up to isomorphism. For more details, see \cite[\S7.6]{FeltenGLDT}; the original reference is \cite{Kato1996}.

We say that $\phi_0\colon X_0^\dagger \to S_0^\dagger$ is \emph{log Calabi--Yau} if $\Omega^{d,\m{log}}_{X_0/S_0} \cong \cO_{X_0}$. In this case, $\m{LD}_{X_0/S_0}$ is unobstructed by \cite[Cor.~1.103]{FeltenGLDT}, see also \cite{FFR2021,Felten2022,FeltenPetracci2022,ChanLeungMa2023}. In the proof, the degeneration of the log Hodge--de Rham spectral sequence and the flatness of the universal log de Rham cohomology (a notion that behaves similarly to infinitesimal and crystalline cohomology) play a crucial role; in Situation~\ref{sitn:basic}, both are consequences of Theorem~\ref{thm:log-smooth-cffc}.

\subsection{Log singularities}

We continue in Situation~\ref{sitn:basic}. For many $\phi\colon X \to S$, the map $\phi\colon X^\dagger \to S^\dagger$ is not log smooth and saturated. The \emph{log singular locus} $\itSigma_X \subset X$ is the smallest closed subset such that $(X \setminus \itSigma_X)^\dagger \to S^\dagger$ is log smooth and saturated.

\begin{ex}\label{ex:quartic-degen}\note{ex:quartic-degen}
	Let
	\[\widebar X = \{xyzw - t(x^4 + y^4 + z^4 + w^4) = 0\} \subset \PP^3 \times \widebar S,\]
	and let $0 \in S \subseteq \widebar S$ be an open neighborhood such that $\phi^{-1}(0)$ is the only singular fiber of $\phi\colon X \to S$. Then, $\phi\colon X \to S$ is semistable except in the $24$ points of
	\[Z = \m{Sing}(X_0) \cap \{x^4 + y^4 + z^4 + w^4 = 0\}.\]
	One may show that $\phi\colon X^\dagger \to S^\dagger$ is not log smooth at any $p \in Z$, in particular, $\itSigma_X = Z$. For more details, see \cite[Ex.~1.79, Ex.~1.80]{FeltenGLDT}. \qedex
\end{ex}

Due to theoretical difficulties in dealing with the log structure in $\itSigma_X$ (for instance, $\cM_X$ does not admit a chart at $p \in \itSigma_X$ in Example~\ref{ex:quartic-degen}), we have adopted the point of view that we choose a closed subset $Z_X \supset \itSigma_X$ and \emph{delete} the log structure in $Z_X$. After setting $U = X \setminus Z_X$, we obtain a \emph{partial log scheme} $X^\partial = (X,U,\cM_U)$, whose log structure is defined only on $U$, and a morphism of partial log schemes $\phi\colon X^\partial \to S^\dagger$. Now, $Z_X$ is the \emph{log indeterminate locus} as there is no log structure determined in $Z_X$; its complement $U$ is the \emph{log determinate locus}. For more details on partial log schemes, see \cite[Defn.~1.75]{FeltenGLDT}. The terminology ``log (in-)determinate locus'' has been introduced in \cite{Felten-gls-gr}.

When certain conditions are met, we say that a morphism of partial log schemes $\phi\colon X^\partial \to S^\dagger$ is \emph{generically log smooth}; most notably, $Z_X \subset X$ has relative codimension $\geq 2$, we have $j_*\cO_U = \cO_X$ for the inclusion $j\colon U \to X$, and $\phi\colon U^\dagger \to S^\dagger$ is log smooth and saturated. Details can be found in \cite[Defn.~8.2]{FeltenGLDT}. One should think of generic log smoothness as a logarithmic analog of a flat family with normal fibers.

In Situation~\ref{sitn:basic}, suppose that we have homogeneous polynomials $f_1,\ldots,f_r$ cutting out $\widebar X \subset \PP^n \times \widebar S$. In \cite{Felten-gls-gr}, the first author developed algorithmic tools to compute (in many cases) some $Z_X \subset X$ such that $\phi\colon (X \setminus Z_X)^\dagger \to S^\dagger$ is log smooth and saturated, and to determine whether $\phi\colon X^\partial \to S^\dagger$ is generically log smooth. Due to practical difficulties in determining the log singular locus $\itSigma_X$ of $\phi\colon X^\dagger \to S^\dagger$, we usually have $Z_X \supsetneq \itSigma_X$ in this setup.

\subsection{Toroidal log singularities}

In Situation~\ref{sitn:basic}, suppose that $\phi\colon X^\dagger \to S^\dagger$ is not log smooth, but that we have some $Z_X \subset X$ such that $\phi\colon X^\partial \to S^\dagger$ is generically log smooth. By the results of \cite{FFR2021,FeltenGLDT}, we still have some handle on the log de Rham complex of $\phi\colon X^\partial \to S^\dagger$ and the deformation theory of $\phi_0\colon X_0^\partial \to S_0^\dagger$ under the condition that the log singularities of $\phi\colon X^\partial \to S^\dagger$ have a specific \'etale local form which arises from toric geometry; we say that $\phi\colon X^\partial \to S^\dagger$ has \emph{toroidal} log singularities. In this section, we explain this notion, and in the next section, we discuss cohomology and deformation theory.
\medskip

Let $\theta\colon \NN \to \itPi$ be an injective and saturated homomorphism of sharp toric monoids. First, it gives rise to a toric morphism of affine toric varieties $\mu\colon A_\itPi \to A_\NN$, where we set $A_\itPi \coloneq \Spec \kk[\itPi]$ and $A_\NN \coloneq \Spec \kk[\NN]$. Now, we endow $A_\NN$ with the divisorial log structure defined by its full toric boundary $\{0\} \subset A_\NN$ and denote the resulting log scheme by $A_\NN^\dagger$. On $A_\itPi$, we use $\mu^{-1}(0) \subset A_\itPi$ rather than the full toric boundary to define a divisorial log structure and obtain a log scheme $A_{\itPi/\NN}^\dagger$. The morphism of log schemes $\mu\colon A_{\itPi/\NN}^\dagger \to A_\NN^\dagger$ is often not log smooth, but we can find an open subset $U_{\itPi/\NN} \subseteq A_{\itPi}$ such that, once we use it as the log determinate locus of the partial log scheme $A_{\itPi/\NN}^\partial = (A_\itPi,U_{\itPi/\NN},\cM_{U_{\itPi/\NN}})$, the morphism $\mu\colon A_{\itPi/\NN}^\partial \to A_\NN^\dagger$ becomes generically log smooth. For more details, see \cite[\S8]{FeltenGLDT}.

\begin{ex}\label{ex:xy-tu}\note{ex:xy-tu}
	Let $\itPi \subseteq \ZZ^3$ be the submonoid generated by $\bar t = (0,1,0)$, $\bar u = (0,0,1)$, $\bar x = (1,0,0)$, and $\bar y = (-1,1,1)$, and let $\theta\colon \NN \to \itPi$ be given by $\theta(1) = \bar t$. Then, $\mu\colon A_\itPi \to A_\NN$ can be identified with $Y \coloneq \Spec \kk[t,x,y,u]/(xy - tu) \to \Spec \kk[t]$. The log singular locus of $\mu\colon A_{\itPi/\NN}^\dagger \to A_\NN^\dagger$ in the total space is $\itSigma_Y = \{t = x = y = u = 0\}$. We have $U_{\itPi/\NN} = A_\itPi \setminus \{x = y = u = 0\}$ so that $Z_Y \supsetneq \itSigma_Y$ for the log indeterminate locus. \qedex
\end{ex}

\begin{defn}\label{defn:toroidal-log-sing}\note{defn:toroidal-log-sing}
	Let $\phi\colon X^\partial \to S^\dagger$ be the morphism of partial log schemes from Situation~\ref{sitn:basic} for some $Z_X \supset \itSigma_X$. Let $\theta\colon \NN \to \itPi$ be an injective and saturated homomorphism of sharp toric monoids. Then, a point $p \in X_0$ is a \emph{toroidal log singularity of type $\theta$} if we can find the following:
	\begin{enum:alph}
		\item a morphism of partial log schemes $h\colon W^\partial \to X^\partial$ over $S^\dagger$ such that $h\colon W \to X$ is \'etale, and such that $h\colon (W \setminus Z_W)^\dagger \to (X \setminus Z_X)^\dagger$ is strict; note that $h\colon W^\partial \to X^\partial$ being a morphism of partial log schemes entails that $W \setminus Z_W \subseteq h^{-1}(X \setminus Z_X)$, but that we may have $W \setminus Z_W \not= h^{-1}(X \setminus Z_X)$;
		\item a point $w \in W$ with $h(w) = p$;
		\item a morphism of partial log schemes $k\colon W^\partial \to A_{\itPi/\NN}^\partial \times_{A_\NN^\dagger} S^\dagger$ over $S^\dagger$ such that the underlying morphism of schemes is \'etale, such that $(W \setminus Z_W)^\dagger \to U_{\itPi/\NN}^\dagger$ is strict, such that $Z_W$ has relative codimension $\geq 2$ with respect to $W \to S$, and such that $k(w) = (o,0)$ for the torus-fixed point $o \in A_\itPi$.
	\end{enum:alph}
	Note that we may have $p \notin \itSigma_X$. We say that $\phi\colon X^\partial \to S^\dagger$ is \emph{log toroidal} if it has only toroidal log singularities. \qeddef
\end{defn}
\begin{rem}
	This is almost the definition of being log toroidal with respect to $S^\dagger \to A_\NN^\dagger$, see \cite[Defn.~8.27]{FeltenGLDT}. The latter definition allows slightly more general local models, but one may use \cite[Lem.~4.13]{FeltenThesis} to show that the two definitions are equivalent in Situation~\ref{sitn:basic}. \qedex
\end{rem}
\begin{rem}
	Toroidal log singularities and log toroidality are defined similarly for $\phi_0\colon X_0^\partial \to S_0^\dagger$ and for $\phi_A\colon X_A^\partial \to S_A^\dagger$ when $A \in \m{Art}_{\kk\lsem t\rsem\sslash\kk}$. \qedex
\end{rem}
\begin{ex}
	In Example~\ref{ex:quartic-degen}, the $24$ points of $Z$ are toroidal log singularities of type $\theta$ for $\theta$ as in Example~\ref{ex:xy-tu}. This is proven in \cite[Ex.~1.80]{FeltenGLDT}. \qedex
\end{ex}

In Situation~\ref{sitn:basic}, the family $\phi\colon X^\partial \to S^\dagger$ frequently has toroidal log singularities of type $\theta$ for $\theta$ as in Example~\ref{ex:xy-tu}. To discuss such and slightly more complicated log singularities more conveniently, we introduce some terminology.

\begin{constr}\label{constr:toroidal-codim-one}\note{constr:toroidal-codim-one}
	Let $\ell \geq 1$ and $m_1,\ldots,m_q \geq 0$. Let
	\[\itPi = \itPi(\ell;m_1,\ldots,m_q) \subseteq \ZZ f \oplus \bigoplus_{i = 0}^q \ZZ e_i^*\]
	be the submonoid generated by $\bar x = f$, $\bar y = -f + \ell e_0^* + \sum_{i = 1}^q m_ie_i^*$, $\bar t = e_0^*$, and $\bar z_i = e_i^*$ for $1 \leq i \leq q$; let $\theta\colon \NN \to \itPi$ be given by $\theta(1) = \bar t$. Then, $\mu\colon A_{\itPi} \to A_\NN$ is isomorphic to 
	\[\Spec \kk[t,x,y,z_1,\ldots,z_q]\left/\middle(xy - t^\ell\prod_{i = 1}^q z_i^{m_i}\right) \to \Spec \kk[t].\]
	For more details, see \cite[Ex.~14.5]{FeltenGLDT}. For $\ell = 1$, $q = 1$, and $m_1 = 1$, we recover the map $\theta\colon \NN \to \itPi$ from Example~\ref{ex:xy-tu}. \qeddef
\end{constr}

\begin{defn}
	In the setup of Definition~\ref{defn:toroidal-log-sing}, a point $p \in X_0$ is a \emph{toroidal log singularity of type $T(\ell;m_1,\ldots,m_q)$} if it is a toroidal log singularity of type \[\theta\colon \enspace \NN \to \itPi(\ell;m_1,\ldots,m_q,0,\ldots,0)\]
	for some $\ell \geq 1$ and $m_1,\ldots,m_q \geq 1$, where we use the notation of Construction~\ref{constr:toroidal-codim-one}.\footnote{``T'' stands for ``toric'' and also for ``type.''} \qeddef
\end{defn}

Now, log singularities as in Example~\ref{ex:xy-tu} are toroidal of type $T(1;1)$, and toroidal log singularities of type $T(2;1)$ are locally modeled on the equation $xy = t^2u$.

\subsection{The log de Rham complex in the toroidal case}

We continue in Situation~\ref{sitn:basic}. Even if $\phi\colon X^\partial \to S^\dagger$ has only toroidal log singularities, the sheaf of log differential forms $\Omega^{1,\m{log}}_{X/S}$ of $\phi\colon X^\dagger \to S^\dagger$ is not always quasicoherent; for instance, this is the case in Example~\ref{ex:xy-tu}. From a different perspective, after deleting the log structure in $Z_X \supset \itSigma_X$, there is no $\Omega^{\bullet,\m{log}}_{X/S}$ defined at all. To circumvent these issues, we set 
\[\cW^\bullet_{X/S} \coloneq j_*\Omega^{\bullet,\m{log}}_{U/S}\]
for any generically log smooth $\phi\colon X^\partial \to S^\dagger$, cf.~\cite[Defn.~8.3]{FeltenGLDT}. Here, $j\colon U \to X$ is the inclusion of the log determinate locus. We say that $\cW^\bullet_{X/S}$ is the \emph{reflexive log de Rham complex}; namely, each $\cW^i_{X/S}$ is a reflexive coherent sheaf, due to the codimension of the complement of $U$ being $\mathrm{codim}(Z_X, X) \geq 2$.

The behavior of the reflexive log de Rham complex under base change is non-trivial. Suppose that we have a Cartesian diagram
\begin{equation}
  \label{eqn:BaseChangeDiagram}
\xymatrix{
	Y^\partial \ar[r]^c \ar[d]^\psi & X^\partial \ar[d]^\phi \\
	T^\dagger \ar[r]^b & S^\dagger \\
}
\end{equation}
in the category of partial log schemes. Then, for every $i \geq 0$, we have a comparison map
\[c^*\cW^i_{X/S} \to \cW^i_{Y/T}\]
of coherent sheaves on $Y$, but there is no a priori reason why this should be an isomorphism. We say that $\phi\colon X^\partial \to S^\dagger$ has the \emph{base change property} if these comparison maps are isomorphisms, see \cite[Defn.~8.47]{FeltenGLDT}.

\begin{thm}\label{thm:base-change-property}\note{thm:base-change-property}
	In Situation~\ref{sitn:basic}, choose $Z_X \supset \itSigma_X$ and assume that $\phi\colon X^\partial \to S^\dagger$ is generically log smooth and has only toroidal log singularities. Then, it has the base change property.
\end{thm}
\begin{proof}
	See \cite[Thm.~8.63]{FeltenGLDT} and the references therein.
\end{proof}

In Theorem~\ref{thm:log-smooth-cffc}, we saw that $\Omega^{\bullet,\m{log}}_{X/S}$ is cffc if $\phi\colon X^\dagger \to S^\dagger$ is log smooth and saturated. The perhaps somewhat surprising fact is that this remains true for $\cW^\bullet_{X/S}$ if $\phi\colon X^\partial \to S^\dagger$ has toroidal log singularities.

\begin{thm}\label{thm:log-toroidal-cffc}\note{thm:log-toroidal-cffc}
	In Situation~\ref{sitn:basic}, choose $Z_X \supset \itSigma_X$ and assume that $\phi\colon X^\partial \to S^\dagger$ is generically log smooth and has only toroidal log singularities. Then, the reflexive log de Rham complex $\cW^\bullet_{X/S}$ is cffc. 
\end{thm}
\begin{proof}
	See \cite[Thm.~8.65]{FeltenGLDT} and the references therein.
\end{proof}

Just as in the log smooth case, we have a logarithmic Bogomolov--Tian--Todorov theorem in the case of toroidal log singularities. Let us assume that $\cW^d_{X/S}$ is a line bundle; the deformation functor 
\[\m{LD}_{X_0/S_0}^\scD\colon \enspace \m{Art}_{\kk\lsem t\rsem\sslash\kk} \to \m{Set}\]
is defined as follows. 
Let $A \in \m{Art}_{\kk\lsem t\rsem\sslash\kk}$, i.e. an Artinian local $\kk$-algebra
with natural maps $\kk\lsem t\rsem \to A \to \kk$.  A deformation $\psi_A\colon
Y_A^\partial \to S_A^\dagger$ of $\phi_0\colon X_0^\partial \to
S_0^\dagger$ is a morphism of partial log schemes, subject to the condition
that $\psi_A\colon Y_A^\partial \to S_A^\dagger$ must be \emph{locally}
isomorphic to the base change $\phi_A\colon X_A^\partial \to S_A^\dagger$
of $\phi\colon X^\partial \to S^\dagger$ along $S_A^\dagger \to
S^\dagger$. The latter map is induced by $\kk[t] \to A$. Now, if
$\cW^d_{X_0/S_0} \cong \cO_{X_0}$ (i.e. in the Calabi-Yau case), then $\m{LD}_{X_0/S_0}^\scD$ is
unobstructed by \cite[Cor.~1.103]{FeltenGLDT}. The symbol $\scD$ refers
to our choice of distinguished deformations to which any deformation must
be locally isomorphic; from a technical perspective, it is a \emph{system
of deformations} in the sense of \cite[Defn.~10.15]{FeltenGLDT}. Again,
the fact that $\cW^\bullet_{X/S}$ is cffc plays a crucial role in the
proof of unobstructedness.

\subsection{Non-toroidal log singularities}\label{sect:non-toroidal-log-singularities-intro}\note{sect:non-toroidal-log-singularities-intro}

For a log singularity, being toroidal is rather special, akin to being log smooth. For instance in intrinsic mirror symmetry (Section~\ref{sect:intrinsic-mirror}), when attempting to smoothen a reducible toric Fano scheme (Section~\ref{sect:Fano-smoothing}), or simply when writing down a family in generators and relations, one often encounters non-toroidal log singularities. Little theory is known for log singularities which are not toroidal; our goal in this article is to understand their cohomology and deformation theory better, by means of studying examples.

\begin{ex}
	Let
	\[\widebar X = \{xyzw - t((x + y + z)w^3 + x^4 + y^4 + z^4) = 0\} \subset \PP^3 \times \widebar S,\]
	and let $0 \in S \subseteq \widebar S$ be an open subset such that $X_0$ is the only singular fiber of $\phi\colon X \to S$. Then, $\{p_w\} = \{t = x = y = z = 0\}$ is a non-toroidal log singularity of $\phi\colon X^\partial \to S^\dagger$. This can be seen by computing $\cA^1_{X_0/S_0} = c_0^*\cW^1_{X/S}$; it is not reflexive at $p_w$, so the base change property does not hold, and $\phi\colon X^\partial \to S^\dagger$ cannot have a toroidal log singularity at $p_w$. For more discussion of this family, see Example~\ref{ex:C21-WS3}. \qedex
\end{ex}

When one examines the proof of Theorem~\ref{thm:log-toroidal-cffc}, one sees that it is a combination of general methods which do not use that the log singularities are toroidal, and explicit constructions and computations in the toric local models. One may therefore speculate that Theorem~\ref{thm:log-toroidal-cffc} can be proven for families $\phi\colon X^\partial \to S^\dagger$ with more general log singularities if the explicit constructions and computations succeed for other local models.

Example~\ref{ex:C31-US1} shows that $\cW^\bullet_{X/S}$ is not always cffc in Situation~\ref{sitn:basic}. Thus, Theorem~\ref{thm:log-toroidal-cffc} cannot be generalized to families with \emph{arbitrary} log singularities.

\begin{quest}
	Which sets of log singularities $\scC$, containing at least one non-toroidal log singularity, have the following property: For any instance of Situation~\ref{sitn:basic} such that $\phi\colon X^\partial \to S^\dagger$ is locally of the form $\scC$, the complex $\cW^\bullet_{X/S}$ is cffc?
\end{quest}

To make the question precise, we define a \emph{class of log singularities} as a set $\scC$ of generically log smooth morphisms of partial log schemes $\mu\colon Y_\mu^\partial \to A_\NN^\dagger$. We say that a morphism of partial log schemes $\phi\colon X^\partial \to S^\dagger$ arising from Situation~\ref{sitn:basic} for a choice of $Z_X \supset \itSigma_X$ is \emph{of class $\scC$} if the analog of the condition in Definition~\ref{defn:toroidal-log-sing} is satisfied with local models in $\scC$. Then, a class of log singularities $\scC$ is \emph{cffc-admissible} if, for any instance of Situation~\ref{sitn:basic} and any $Z_X \supset \itSigma_X$ such that $\phi\colon X^\partial \to S^\dagger$ is of class $\scC$, the complex $\cW^\bullet_{X/S}$ is cffc.

In practice, the situation is usually that we have a class of log singularities $\scC$ specified by a concrete geometric problem or construction, and that we wish to determine if this given class $\scC$ is cffc-admissible. In the next two sections, we discuss two settings where non-toroidal log singularities arise naturally: intrinsic mirror symmetry, and smoothing reducible toric Fano schemes.

\subsection{Intrinsic mirror symmetry}\label{sect:intrinsic-mirror}\note{sect:intrinsic-mirror}

Our first source of non-toroidal log singularities comes from intrinsic mirror symmetry \cite{GrossHackingKeel2015,GS-intrinsic,GS-canonical,ACGS2020,GHKS2022,ACGS2025,GS2026-intrinsic}, in particular the construction in Gross--Hacking--Keel \cite{GrossHackingKeel2015} (see next paragraph). For the present article, this serves as motivation only; we will not prove any result concerning the cohomology or deformation theory specifically of log singularities arising from the intrinsic mirror construction of \cite{GrossHackingKeel2015}.

The construction in \cite{GrossHackingKeel2015} starts with a smooth projective rational surface $Y$ and an anti-canonical cycle of rational curves $D \subset Y$. The \emph{intrinsic mirror} of $(Y,D)$ is a formal flat deformation $\fX \to \Spf \kk\lsem P\rsem$ over the formal spectrum of $\kk\lsem P\rsem$ of some affine surface $\VV_n$ for a certain sharp toric monoid $P$. Here, $n \geq 1$ is the number of irreducible components of $D$, and the scheme $\VV_n$ depends only on $n$.

\begin{ex}
	We discuss the intrinsic mirror of a cubic surface as treated in \cite{GHKS2022}. Let $\widebar D = \widebar D_1 \cup \widebar D_2 \cup \widebar D_3 \subset \PP^2$ be the toric boundary, and let $\rho\colon Y \to \PP^2$ be the blowup in two general points on each $\widebar D_i$. Let $D = D_1 \cup D_2 \cup D_3$ be the strict transform of $\widebar D$; it is an anti-canonical divisor in $Y$. For $1 \leq i \leq 3$, let $L_{ij}$ for $1 \leq j \leq 8$ be the eight lines on $Y$ which intersect $D_i$ and are not contained in $D$.
	
	Note that $\m{Pic}(Y) \cong \ZZ^7$. Let $P = \m{NE}(Y) \subseteq \m{Pic}(Y)$ be the submonoid of effective curve classes; the associated real cone is generated by the classes of the $27$ lines, see \cite[Ch.~4, Prop.~1.12]{Huybrechts2023}. By \cite[Thm.~1.1]{GHKS2022}, the intrinsic mirror to $(Y,D)$ is given by the equation 
	\[\vartheta_1\vartheta_2\vartheta_3 = \sum_{i = 1}^3 \m{z}^{D_i}\vartheta_i^2 + \sum_{i = 1}^3 \left(\sum_{j = 1}^8 \m{z}^{L_{ij}}\right) \m{z}^{D_i} \vartheta_i + \sum_\pi \m{z}^{\pi^*H} + 4\m{z}^D\]
	in $\kk[P][\vartheta_1,\vartheta_2,\vartheta_3]$. Here, $\pi$ ranges over all birational morphisms $\pi\colon Y \to \PP^2$ with the property that $\pi$ induces an isomorphism $D \cong \widebar D$, up to automorphisms of $(\PP^2,\widebar D)$; there are $\lvert W(D_4)\rvert = 192$ of them, cf.~\cite[Ex.~6.13]{GrossHackingKeel2015}. The central fiber is $\VV_3 = \{\vartheta_1\vartheta_2\vartheta_3 = 0\} \subset \IA^3$.
	
	We have $\m{rk}(P) = 7$. To obtain a one-parameter family, we proceed as follows. Let $\phi\colon P \to \NN$ be given by $C \mapsto C \cdot (-K_Y)$, the intersection product. Since $-K_Y$ is ample, we have $C \cdot (-K_Y) \geq 0$, and $\phi$ is well-defined. After applying $\phi$ to the coefficients of the equation and setting $\kk[\NN] = \kk[t]$, we obtain
	\[\vartheta_1\vartheta_2\vartheta_3 = \sum_{i = 1}^3\vartheta_i^2 + \sum_{i = 1}^3 8t^2\vartheta_i + 196t^3.\]
	Namely, the projection formula yields $\phi(\pi^*H) = (\pi^*H) \cdot D = H \cdot (\pi_*D) = 3$. Now, the affine family $\phi\colon X \to \IA^1_t$ defined by this equation has a non-toroidal log singularity at $\{0\} = \{t = \vartheta_1 = \vartheta_2 = \vartheta_3 = 0\}$. One may confirm this by showing that the base change property fails at $0$, so Theorem~\ref{thm:base-change-property} proves that the log singularity is not toroidal. \qedex
\end{ex}

Now, suppose that intrinsic mirror symmetry gives us a modular compactification
\[\pi\colon \enspace \overline\cK \to \overline\cM\]
of some moduli space of K3 surfaces. Gross and Siebert have communicated to the authors that this is, at the time of writing, work in progress together with Hacking and Keel, see also \cite[p.~5]{Gross2025-remarks}.

One would like to have a log de Rham complex $\cA^\bullet_{\overline\cK/\overline\cM}$ which is cffc. In particular, for any one-parameter family arising from a map $b\colon S \to \overline\cM$, we would like to have a log de Rham complex which is cffc, and $\cW^\bullet_{X/S}$ is a natural candidate to start with if the general fiber of $\phi\colon X = \overline\cK \times_{\overline\cM} S \to S$ is smooth. Thus, it is desirable for the class of log singularities $\scC^\m{ims}$ arising from intrinsic mirror symmetry to be cffc-admissible.\footnote{Given that intrinsic mirror symmetry provides formal families rather than algebraic families, one should work with a class of \emph{formal} log singularities and adjust the definition of cffc-admissibility accordingly. Since intrinsic mirror symmetry serves only as motivation in this article, we do not elaborate the details.}

Since $\pi\colon \overline\cK \to \overline\cM$ is (most likely) a Calabi--Yau family, we may expect some form of the logarithmic Bogomolov--Tian--Todorov theorem in this situation. This would show that $\overline\cM$ is log smooth over the trivial log point, in case this is not already evident from the construction of $\pi\colon \overline\cK \to \overline\cM$.

\subsection{Smoothing reducible toric Fano schemes}\label{sect:Fano-smoothing}\note{sect:Fano-smoothing}

This is our second source of non-toroidal log singularities. In \cite[\S1.12.3]{FeltenGLDT}, the first author described a strategy to smoothen reducible toric Fano schemes $\PP_0(\widetilde\Delta)$ of dimension $3$ by means of logarithmic deformation theory, communicated to him by Alessio Corti and Helge Ruddat. Preliminary results toward following this strategy are published in \cite{CHP2024,CGR2025,Grafnitz2025,CortiRuddat2026,Filip2026}.

Let $\Delta \subset \RR^3$ be a reflexive polytope, giving rise to the Gorenstein toric Fano variety $\PP(\Delta)$. From $\Delta$, one constructs a certain unbounded polyhedron $\widetilde\Delta \subset \RR^3 \oplus \RR$ and a toric morphism $\psi\colon \PP(\widetilde\Delta) \to \IA^1_t$. Its general fiber is $\PP(\Delta)$ while the central fiber $\PP_0(\widetilde\Delta)$ has many irreducible components, each of which corresponds to a maximal cell in the central subdivision of $\Delta$. This is called the \emph{Mumford degeneration}.

In the case that $\PP(\Delta)$ is not smooth, we wish to find another flat and proper family $\phi\colon X \to \IA^1_t$ with $\phi^{-1}(0) = \PP_0(\widetilde\Delta)$ and whose general fiber $X_t = \phi^{-1}(t)$ is a smooth Fano variety. We use the following strategy: First, we construct a suitable partial log structure on $V = \PP_0(\widetilde\Delta)$ over the standard log point $S_0^\dagger$, and then we use logarithmic deformation theory to obtain $\phi\colon X \to \IA^1_t$, or at least a formal smoothing $\phi\colon \fX \to \Spf \kk\lsem t\rsem$.

The scheme $V$ carries a sheaf of commutative monoids $\cP_V$ which serves as the ghost sheaf of a not-yet-constructed log structure, together with a global section $\bar\rho \in \Gamma(V,\cP_V)$ which serves as the image of $1 \in \overline\cM_{S_0}$ in $\cP_V$ under a not-yet-constructed morphism of log schemes to $S_0^\dagger = \Spec(\NN \to \kk)$. For later use, we mention that $V$ admits a stratification $\cS_kV = \{v \in V \ | \ \m{rk}(\cP_{V,v}) = k + 1\}$ according to the rank of $\cP_V$. Geometrically, $\cS_0V \subseteq V$ is the open subset where $V$ is normal (and in fact smooth), $\cC_1V = \bigsqcup_{k \geq 1}\cS_kV$ is the (two-dimensional) closed subset where at least two irreducible components intersect, and $\cC_2V = \bigsqcup_{k \geq 2}\cS_kV$ is the (one-dimensional) closed subset where three or more irreducible components intersect. Finally, $\cC_3V = \cS_3V = \{o\}$ consists of a single point which happens to be the locus where all irreducible components intersect.

There is a sheaf $\cL\cS_V$ on $V$ which classifies isomorphism classes of \emph{log smooth structures}, i.e., for an open subset $W \subseteq V$, the set $\Gamma(W,\cL\cS_V)$ consists of (isomorphism classes of) log structures $\cM_W$ on $W$ with an isomorphism $\overline\cM_W \cong \cP_V|_W$ and a log smooth, saturated, and vertical morphism $\chi_0\colon (W,\cM_W) \to S_0^\dagger$ with $\chi_0^\splus(1) = \bar\rho|_W$ under $\chi_0^\splus\colon \chi^{-1}\overline\cM_{S_0} \to \cP_V|_W$.

If we have a global section $\lambda \in \Gamma(V,\cL\cS_V)$, we obtain a log smooth morphism $\chi_0\colon V^\dagger \to S_0^\dagger$ that we can use to obtain a (formal) smoothing via logarithmic deformation theory. Now, the sheaf $\cL\cS_V$ admits an explicit description in terms of $\Delta$, and from this, one can see that $\cL\cS_V$ does \emph{not} have a global section. Instead, we wish to find:
\begin{enum:alph}
	\item a closed subset $Z_0 \subset \cC_1V$ with $\m{dim}(Z_0) = 1$, $\m{dim}(Z_0 \cap \cC_2V) = 0$, and $Z_0 \cap \cC_3V = \varnothing$;
	\item a log smooth structure $\lambda \in \Gamma(U_0,\cL\cS_V)$ for $U_0 = V \setminus Z_0$.
\end{enum:alph}
There is one canonical choice $(Z_0^\m{can},U_0^\m{can},\lambda^\m{can})$ that comes from the Mumford degeneration $\psi\colon \PP(\widetilde\Delta) \to \IA^1_t$. However, if the general fiber $\PP(\Delta)$ of the Mumford degeneration is not smooth, it is difficult (though perhaps not impossible) to use logarithmic deformation theory to obtain a smoothing $\phi\colon X \to \IA^1_t$. Namely, any point $v \in Z_0^\m{can} \cap \cS_1V$ is a toroidal log singularity of type $T(\ell;m)$, and we might have $m \geq 2$. The logarithmic deformation theory of these log singularities is difficult to control, and the obvious deformation (coming from the toric local model) has a singularity of the type $xy = u^m$ in its general fiber.

To solve these problems, we want to find another $(Z_0,U_0,\lambda)$ such that any $v \in Z_0 \cap \cS_1V$ is a toroidal log singularity of type $T(\ell;1)$. Namely, for such log singularities, we have a well-understood deformation theory, and the general fiber is always smooth. However, the price for having only toroidal log singularities of type $T(\ell;1)$ in $Z_0 \cap \cS_1V$ is that (some of) the finitely many log singularities in $Z_0 \cap \cC_2V$ are not toroidal.

In \cite{CGR2025}, a set $\scC_0^\m{adm}$ of (``admissible'') generically log smooth morphisms 
%$\mu_0\colon M_\mu^\partial \to S_0^\dagger$ 
was defined by means of a combinatorial construction. In work in progress, Corti, Ruddat, and Gr\"afnitz show that $V = \PP_0(\widetilde\Delta)$ carries a partial log smooth structure $(Z_0^\m{adm},U_0^\m{adm},\lambda^\m{adm})$ such that $Z_0^\m{adm} \cap \cS_1V$ consists of toroidal log singularities of type $T(\ell;1)$, and such that the finitely many log singularities in $Z_0^\m{adm} \cap \cC_2V$ are of class $\scC_0^\m{adm}$. The construction of $(Z_0^\m{adm},U_0^\m{adm},\lambda^\m{adm})$ depends both on conditions on $\Delta$ and the choice of additional data.

If we want to use $(Z_0^\m{adm},U_0^\m{adm},\lambda^\m{adm})$ to smooth $V$, we should know how each $\mu_0\colon M_\mu^\partial \to S_0^\dagger$ in $\scC_0^\m{adm}$ can be smoothed. Conjecturally, there is a preferred morphism of partial log schemes $\mu\colon Y_\mu^\partial \to A_\NN^\dagger$ with central fiber $\mu_0\colon M_\mu^\partial \to S_0^\dagger$ and smooth general fiber, for each $\mu_0\colon M_\mu^\partial \to S_0^\dagger$ in $\scC^\m{adm}_0$. They form a conjectural class of log singularities $\scC^\m{adm}$, and the smoothing $\phi\colon X \to \IA^1_t$ that we want to construct should have only toroidal log singularities of type $T(\ell;1)$ and log singularities of class $\scC^\m{adm}$.

\begin{ex}[A32-ADM]\label{ex:adm-model}\note{ex:adm-model}
	Consider the affine family
	\[\mu\colon Y = \Spec \kk[t,x,y,z,w,u]/(xy - z^2 - t^2,\: zw - tu) \to \IA^1_t.\]
	The only singular fiber is $M = Y_0$, and the log singular locus is $\itSigma_Y = Z(x) \cup Z(y) \cup Z(w)$, as was shown in \cite[Ex.~2.19, Rem.~2.20]{Felten-gls-gr}. By writing down the comparison \'etale morphisms with the local model explicitly, one may show that every point $0 \not= p \in Z(x)$ and $0 \not= p \in Z(y)$ is a toroidal log singularity of type $T(1;1)$. Similarly, every point $0 \not= p \in Z(w)$ is a toroidal log singularity of type $T(2;1)$. The point $p = 0$ is a non-toroidal log singularity. The morphism of partial log schemes $\mu\colon Y^\partial \to A_\NN^\dagger$ is expected to lie in $\scC^\m{adm}$. For more information and a diagram, see \cite[Ex.~2.19]{Felten-gls-gr} and \cite[Ex.~1.158]{FeltenGLDT}. \qedex
\end{ex}

To construct the smoothing $\phi\colon X \to \IA^1_t$, we would like to apply a logarithmic Bogomolov--Tian--Todorov theorem to the morphism of partial log schemes $\phi_0\colon V^\partial_\m{adm} \to S_0^\dagger$ obtained from $(Z_0^\m{adm},U_0^\m{adm},\lambda^\m{adm})$. Now, $\phi_0\colon V^\partial_\m{adm} \to S_0^\dagger$ is not log Calabi--Yau since $\cW^d_{X_0/S_0} \cong \omega_V^\circ$, the relative dualizing sheaf, which is ample. The local models in $\scC^\m{adm}_0$ have non-toroidal log singularities only in isolated points, so the same holds for $\phi_0\colon V^\partial_\m{adm} \to S_0^\dagger$. We hope to find a \emph{log transversal section}  $s \in \Gamma(V,\omega_V^\circ)$ in the sense of \cite[Defn.~18.39]{FeltenGLDT} such that $E_0 = \m{div}(s) \subset V$ is disjoint from the locus $Z_0^\m{adm} \cap \cC_2V$ of (possibly) non-toroidal log singularities. This log transversal section yields, via the construction in \cite[\S18.1]{FeltenGLDT}, a new morphism of partial log schemes $\phi_0'\colon V^\partial_\m{CY} \to S_0^\dagger$ which is log Calabi--Yau, and all of whose log singularities are either toroidal (in a slightly more general sense than Definition~\ref{defn:toroidal-log-sing}) or of class $\scC_0^\m{adm}$.

If it is possible to prove an unobstructedness theorem for $\phi_0'\colon V_\m{CY}^\partial \to S_0^\dagger$ with the strategy that we have used for the log smooth and log toroidal Bogomolov--Tian--Todorov theorems, relying on the degeneration of the log Hodge--de Rham spectral sequence and the flatness of the universal log de Rham cohomology, then we can also expect that $\scC = \scC^\m{adm} \cup \scC^\m{tor}$ is cffc-admissible, where $\scC^\m{tor}$ is the class of toric local models. As a negative result of this article, we will see below that this is (in all likelihood) not the case.

\subsection{Cohomology in the non-toroidal case}
\label{sec:CohomologyInTheNonToroidalCase}

Let $\scC$ be a class of log singularities, for example $\scC = \scC^\m{adm} \cup \scC^\m{tor}$ from the previous section. Then, we can test the hypothesis that $\scC$ is cffc-admissible by writing down an \emph{example} of such a family and determining if $\cW^\bullet_{X/S}$ is cffc. The first original contribution of this article is that we develop a numerical criterion to determine if $\cW^\bullet_{X/S}$ is cffc.

For a non-toroidal log singularity, $\phi\colon X^\partial \to S^\dagger$ does not always have the base change property. Therefore, we introduce notation 
\[\cA^\bullet_{X/S} \coloneq \cW^\bullet_{X/S}\quad \text{and} \quad \cA^\bullet_{Y/T} \coloneq c^*\cA^\bullet_{X/S} = c^*\cW^\bullet_{X/S};
\]
see Diagram (\ref{eqn:BaseChangeDiagram}).
In particular, $\cA^\bullet_{X_0/S_0}$ is the restriction of $\cA^\bullet_{X/S}$ to the
central fiber $X_0$. 

\begin{prop}\label{prop:numerical-reduction}\note{prop:numerical-reduction}
	In Situation~\ref{sitn:basic}, choose some $Z_X \supset \itSigma_X$ and assume that $\phi\colon X^\partial \to S^\dagger$ is generically log smooth. Let $0 \not= s \in S^\star$ be a closed point, and suppose that
	\[h^j(X_0,\cA^i_{X_0/S_0}) = h^j(X_s,\Omega_{X_s/\kk}^i)\]
	for all $i,j \in \ZZ$. Then, $\cA^\bullet_{X/S}$ is cffc.
\end{prop}

\noindent The reader can find the proof on page~\pageref{pf:numerical-reduction}.

\subsection{The examples and their tags}

We have applied Proposition~\ref{prop:numerical-reduction} to a number of explicit examples of Situation~\ref{sitn:basic}. Every example that we consider has a description in this article and typically two files in \cite{FZ26} with code that we used to compute the geometric properties and the Hodge numbers with the computer algebra systems Macaulay2 \cite{M2} and OSCAR \cite{OSCAR}. To see more easily which file corresponds to which example, every example has a unique tag, which is always included in the filename; for instance, the tag of Example~\ref{ex:C21-FQA} is C21-FQA.

\begin{rem}
	We explain the meaning of the tags. Since we are considering only complete intersections, we have $\omega_{X_t} \cong \cO_X(\delta)$ for $\delta = - n - 1 + \sum_{j = 1}^r \m{deg}(f_j)$. The first symbol indicates if $\delta$ is negative (N), zero (C for Calabi--Yau), or positive (P). The second symbol indicates the dimension of $X_t$, and the third symbol the number of equations. The last three symbols can be either letters or numbers; they are de facto arbitrary, but we have tried to make them somewhat suggestive. For instance, ``FQA'' stands for ``Fermat quartic.'' We explain the abbreviations in the code files.	The tags of \emph{affine} examples start with A, for instance A32-ADM for Example~\ref{ex:adm-model}. \qedex
\end{rem}

\subsection{Experimental observations: surfaces}\label{sect:experimental-2}

In \cite[Defn.~9.7]{FeltenGLDT}, the first author introduced the notion of a scheme with \emph{toroidal crossing singularities}. This is a generalization of a normal crossing scheme, and the main property is that $X_0$ has toroidal crossing singularities if $\phi_0\colon X_0^\dagger \to S_0^\dagger$ is log smooth, saturated, and vertical over the standard log point. For example, $\{xyz = 0\}$, $\{xy = zw = 0\}$, and $\{xy - w^2 = zw = 0\}$ all have toroidal crossing singularities. To our great surprise, we have found for surfaces $X_0$ with such singularities:

\begin{obs}\label{obs:two-dim-cffc}\note{obs:two-dim-cffc}
	In Situation~\ref{sitn:basic}, assume that $\widebar X \subset \PP^n \times \widebar S$ is a complete intersection. Choose $Z_X \supset \itSigma_X$ and assume that $\phi\colon X^\partial \to S^\dagger$ is generically log smooth. Assume furthermore that $d = 2$ so that $\phi\colon X \to S$ is a degeneration of surfaces. Assume that $X_0$ has toroidal crossing singularities. In all instances of this setup that we tested, our computations either show or strongly suggest that $\cA^\bullet_{X/S}$ is cffc. Many of them have log singularities which are not toroidal.
\end{obs}

We discuss these examples in detail in Sections~\ref{sect:surfaces-IP3}, \ref{sect:surfaces-IP4}, and \ref{sect:surfaces-IP5}. There are $13$ examples; in ten examples, the Hodge numbers in Proposition~\ref{prop:numerical-reduction} match over $\QQ$, and in three examples, we could only verify that they match over a finite field.

On the one hand, the observation suggests that $\cA^\bullet_{X/S}$ is cffc if $X_0$ has only toroidal crossing singularities. On the other hand, one also has to keep in mind that all families considered are complete intersections in $\PP^n \times \IA^1_t$, and that we have $h^0(X_s,\Omega^1_{X_s/\kk}) = 0$ and $h^2(X_s,\Omega^1_{X_s/\kk}) = 0$ for a smooth fiber $X_s = \phi^{-1}(s)$, so if $\cA^\bullet_{X/S}$ is cffc, not only the local structure of the singularities but also the global geometry could play a role. This concern is mitigated by Example~\ref{ex:N21-SSE} (N21-SSE), Example~\ref{ex:C21-BE6} (C21-BE6), and Example~\ref{ex:C21-SBP} (C21-SBP), where $X_0$ is a singular normal surface and $\cA^\bullet_{X/S}$ is not cffc.

Interestingly, all these latter three examples have at least one non-canonical singularity. On the other hand, Example~\ref{ex:N21-SE6} (N21-SE6) has only a single $E_6$-singularity, and $\cA^\bullet_{X/S}$ is cffc. Moreover, in Example~\ref{ex:umbrella-smoothing} (N21-SWU), $X_0$ is not normal, $\cA^\bullet_{X/S}$ is cffc, and $X_0$ has semi-smooth singularities (in the sense of \cite[Defn.~4.2]{KSB1988}) but not toroidal crossing singularities. As a joint generalization of the two cases, this suggests that $\cA^\bullet_{X/S}$ might be cffc if $X_0$ has only \emph{semi-canonical} singularities.\footnote{For the definition of semi-canonical singularities, see \cite[Defn.~4.17]{KSB1988}. For a proof that both semi-smooth and canonical singularities are semi-canonical, see \cite[Thm.~4.21]{KSB1988}.}

All toroidal crossing singularities are \emph{semi--log canonical} \cite[Lem.~9.8]{FeltenGLDT}, but already $\{xyz = 0\} \subset \IA^3$ is not semi-canonical. The reader might therefore think that $X_0$ having semi--log canonical singularities is the correct condition. Moreover, the pair $(X,X_0)$ is log canonical if $X_0$ has semi--log canonical singularities in our setup. However, $X_0$ is a normal surface with log canonical singularities in Example~\ref{ex:N21-SSE} (N21-SSE), and $\cA^\bullet_{X/S}$ is not cffc.

We currently cannot assess if it is true that, in those examples where $X_0$ has toroidal crossing singularities (Observation~\ref{obs:two-dim-cffc}), the log singularities of $\phi\colon X^\partial \to S^\dagger$ arise from the two-dimensional Gross--Hacking--Keel intrinsic mirror construction (see Section~\ref{sect:intrinsic-mirror}). The reason is that not enough is known about which singularities arise from this construction.

\begin{rem}\label{rem:Jorder}
	Let $X/\CC$ be a normal projective variety, and let $\Omega^{[\bullet]}_{X/\CC}$ be the reflexive de Rham complex. As far as we see, it follows from \cite[Prop.~5.1]{Jorder2014} that the Hodge--de Rham spectral sequence of $\Omega^{[\bullet]}_{X/\CC}$ degenerates at $E_1$ if $X$ is a surface with rational singularities.\footnote{See also \cite{HuberJorder2014} for related results.} This result is related to but slightly different from our discussion of $\cA^\bullet_{X/S}$, even after adjusting for the fact that we consider a family while [\textit{loc.\,cit.}] considers a single variety: In Example~\ref{ex:N21-SE6} (N21-SE6), the central fiber $X_0$ has du Val and hence rational singularities, but $\cA^1_{X_0/S_0}$ is not reflexive at $p_w$, so our computation shows the degeneration of a new de Rham complex different from $\Omega^{[\bullet]}_{X_0/\CC}$. \qedex
\end{rem}

\subsection{Experimental observations: threefolds}

Maybe more in accordance with a sober expectation, we have found for threefolds:

\begin{obs}
	In Situation~\ref{sitn:basic}, assume that $\widebar X \subset \PP^n \times \widebar S$ is a complete intersection. Choose $Z_X \supset \itSigma_X$ and assume that $\phi\colon X^\partial \to S^\dagger$ is generically log smooth. Assume that $d = 3$ so that $\phi\colon X \to S$ is a degeneration of threefolds. Assume furthermore that $X_0$ has toroidal crossing singularities. In every instance of this setup with non-toroidal log singularities that we tested, our computations either show or strongly suggest that $\cA^\bullet_{X/S}$ is not cffc.
\end{obs}

We discuss these examples briefly in Section~\ref{sect:threefolds}. There are five examples in total, and four of them have non-toroidal log singularities. Our identification of non-matching Hodge numbers is only in Example~\ref{ex:N32-ADM} (N32-ADM) conditional.

Example~\ref{ex:C31-FQI} (C31-FQI) has toroidal log singularities, and indeed, we found that $\cA^\bullet_{X/S}$ is cffc (the Hodge numbers match over $\QQ$). This can be considered as experimental evidence that Theorem~\ref{thm:log-toroidal-cffc} is indeed correct.

Example~\ref{ex:N32-ADM} (N32-ADM) has two non-toroidal log singularities $q_0$, $q_\infty$ modeled on Example~\ref{ex:adm-model}. We believe that it has otherwise only toroidal log singularities (though we did not prove this), so we assume that its log singularities are of the conjectural class $\scC = \scC^\m{adm} \cup \scC^\m{tor}$ introduced in Section~\ref{sect:Fano-smoothing}. If our assumptions on the log singularities other than $q_0$, $q_\infty$ are correct, our computations show that $\cA^\bullet_{X/S}$ is not cffc. Therefore, $\scC = \scC^\m{adm} \cup \scC^\m{tor}$ is presumably not cffc-admissible.

\subsection{Logarithmic deformations in the non-toroidal case}

In Situation~\ref{sitn:basic}, let us only assume that we have some $Z_X \subset X$ such that $\phi\colon X^\partial \to S^\dagger$ is generically log smooth, and that $\cW^d_{X/S}$ is a line bundle, but not that $\phi\colon X^\partial \to S^\dagger$ would have toroidal log singularities. Then, $\phi\colon X^\partial \to S^\dagger$ gives rise to a deformation functor, which we now denote by 
\[\m{ELD}_{X_0/S_0}^\scD\colon \enspace \m{Art}_{\kk\lsem t\rsem\sslash\kk} \to \m{Set}.\]
For more details on its construction, we refer the reader to Section~\ref{sect:log-defo-functor} and the references therein. If $\cA^d_{X_0/S_0} \cong \cO_{X_0}$, and if $\cA^\bullet_{X/S}$ is cffc, then this deformation functor is unobstructed, see Theorem~\ref{thm:ELD-unobstructed}.

Many of our examples of non-toroidal degenerations of surfaces below satisfy the conditions of Theorem~\ref{thm:ELD-unobstructed}. All of them are new examples of log singular surfaces with a smooth local moduli space.

\subsection{The computations}

Every computation has two parts. In the first part, we use Macaulay2 \cite{M2} and the tools of \cite{Felten-gls-gr,FeltenZenodo} to establish that the example satisfies the geometric conditions laid out in Section~\ref{sect:log-geo-setup}. These computations can be easily carried out on a personal computer. In \cite{FZ26}, we have a separate file for each of these computations. Furthermore, \cite{FZ26} contains additional methods that complement the methods of \cite{Felten-gls-gr,FeltenZenodo}, see Section~\ref{sect:non-toroidal-log-sing}.

In the second part, we use OSCAR \cite{OSCAR} to compute $h^j(X_0,\cA_{X_0/S_0}^i)$ in order to apply Proposition~\ref{prop:numerical-reduction}. With a suitable strategy, most of these computations can also be carried out on a personal computer. In addition, we used the compute servers of the Algebra, Geometry, and Computer Algebra Group (AGAG) at the Department of Mathematics at the RPTU Kaiserslautern for cases that needed a larger amount of resources; among the examples that we finally included in this article, this was necessary only for Example~\ref{ex:C23-ZW3} (C23-ZW3).

For the second part, a detailed tabulation of computation times and memory assignments is included in the code file in \cite{FZ26}.

\subsection{Structure of the article}

The main text is divided into two parts. Part~\ref{part:theory} contains our theoretical results. In Section~\ref{sect:log-geo-setup}, we explain the precise setup of our examples, and how to verify, for instance, that $\phi\colon X^\partial \to S^\dagger$ is generically log smooth after constructing $Z_X \supset \itSigma_X$. In Section~\ref{sect:homological-algebra}, we discuss spectral sequences, the cffc property, and the proof of Proposition~\ref{prop:numerical-reduction}. In Section~\ref{sect:log-defo-functor}, we discuss the logarithmic deformation functor $\m{ELD}^\scD_{X_0/S_0}$ and its Bogomolov--Tian--Todorov theorem. In Section~\ref{sect:algorithm}, we discuss the details of our computations.

Part~\ref{part:examples} contains the examples which we computed, most of them with a description of their geometry. In the first three sections of Part~\ref{part:examples}, we consider degenerations of complete-intersection surfaces in $\PP^n$: in Section~\ref{sect:surfaces-IP3} inside $\PP^3$, in Section~\ref{sect:surfaces-IP4} inside $\PP^4$, and in Section~\ref{sect:surfaces-IP5} inside $\PP^5$. Finally, in Section~\ref{sect:threefolds}, we consider degenerations of complete-intersection threefolds in $\PP^n$.

\subsection{Acknowledgments}

The first author thanks Bernd Siebert and Mark Gross for discussions about intrinsic mirror symmetry, in particular Bernd Siebert for first suggesting to study the log singularities that occur in the Gross--Hacking--Keel mirror construction.

\medskip

\noindent\textbf{Funding.} The first author was financially supported by the German Research Foundation (DFG) through grant FE 2102/1-1 and by the Royal Society through Lukas Brantner's grant URF\textbackslash R1\textbackslash 211075. The first author thanks Columbia University and the University of Oxford for their hospitality. Further, he thanks the Deutsche Forschungsgemeinschaft 
(DFG, German Research Foundation) – Project-ID 286237555 – TRR 195 for support during a visit at the Universit\"at des Saarlandes in 2024. 

The second author has been funded by the Deutsche Forschungsgemeinschaft 
(DFG, German Research Foundation) – Project-ID 286237555 – TRR 195. 

\medskip

\noindent\textbf{AI disclosure.} The first author has obtained complimentary access to ChatGPT and Claude. With the exception of the instances specified in the main text, he used AI only for literature research, checking the correctness of the manuscript, improving the exposition, and similar tasks. The second author has not used any AI-driven tools at all for the work on this article. In particular, we wrote all accompanying code \cite{FZ26} ourselves.

\newpage

\part{Theory}\label{part:theory}

\section{The log-geometric setup}\label{sect:log-geo-setup}\note{sect:log-geo-setup}

We describe the log geometry of our examples in a general form. This form, a variant of \cite[Sitn.~2.18]{Felten-gls-gr}, is such that the results of \cite{Felten-gls-gr} allow us to verify computationally that the constructed family $\phi\colon (X,U,\cM_U) \to (S,\cM_S)$ is generically log smooth.

\subsection{The degenerating family}

We fix a field $\kq$ as the base field for our computations. For theoretical results, we assume $\m{char}(\kq) = 0$, and that we have a field extension $\kq \to \kk$, where $\kk$ is the algebraically closed field from Situation~\ref{sitn:basic}. When actually performing a computation, $\kq/\QQ$ should be a finite field extension for practical reasons. All examples considered in this article have $\kq = \QQ$. We also use $\kq = \FF_p$ as a less costly heuristic in very resource-intense calculations. Our examples start with an input datum of the following kind:

\begin{dat}\label{dat:basic-input}\note{dat:basic-input}
	We have $R_{\widebar\scS} = \kq[t]$ and a graded polynomial ring
	\[R_\PP^\bullet = \kq[t,x_0,\ldots,x_n]\]
	with $\m{deg}(t) = 0$ and $\m{deg}(x_i) = 1$ for $0 \leq i \leq n$. For some $r \geq 0$, we have homogeneous polynomials $f_1,\ldots,f_r \in R_\PP^\bullet$. We have an integer $d \geq 1$ which serves as the \emph{intended relative dimension} of our degeneration. \qedex
\end{dat}

In this situation, let $R_{\widebar S} = R_{\widebar\scS} \otimes_\kq \kk$, let $\widebar\scS = \Spec R_{\widebar\scS}$, and let $\widebar S = \Spec R_{\widebar S}$. We have the graded quotient ring
\[R_{\widebar\scX}^\bullet = R_\PP^\bullet/(f_1,\ldots,f_r)\]
and a closed subscheme
\[\widebar\scX = \operatorname{Proj} R_{\widebar\scX}^\bullet \subset \PP^n \times \widebar\scS.\]
After setting $R_{\widebar X}^\bullet = R_{\widebar\scX}^\bullet \otimes_\kq \kk$, we have $\widebar X = \operatorname{Proj} R_{\widebar X}^\bullet \subset \PP^n \times \widebar S$. Let $\phi\colon \widebar X \to \widebar S$ and $\phi\colon \widebar\scX \to \widebar\scS$ be the compositions of these embeddings with the projections to $\widebar S$ respectively $\widebar\scS$; they are proper. The family $\phi\colon \widebar X \to \widebar S$ is our primary object of interest.

Let
\[\widebar X_\vartriangle = \Spec R_{\widebar X}^\bullet \subset \IA^{n + 1} \times \widebar S\]
be the affine cone over $\widebar X$, and let $\phi_\vartriangle\colon \widebar X_\vartriangle \to \widebar S$ be the composition with the second projection. We define $\widebar\scX_\vartriangle$ and $\phi_\vartriangle\colon \widebar\scX_\vartriangle \to \widebar\scS$ analogously.

\begin{ass}\label{ass:apgls}\note{ass:apgls}
	We assume that $\phi_\vartriangle\colon \widebar\scX_\vartriangle \to \widebar\scS$ is \emph{algorithmically pre--generically log smooth} in the sense of \cite[Defn.~2.11]{Felten-gls-gr}. Concretely, we assume:
	\begin{enumerate}[label=\textnormal{(\alph*)}]
		\item\label{cond:integral} the ring $R_{\widebar\scX}^\bullet$ is integral, Cohen--Macaulay, and of dimension $d + 2 \geq 3$;
		\item\label{cond:injective} the ring homomorphism $R_{\widebar\scS} \to R_{\widebar\scX}^\bullet$ is injective;
		\item\label{cond:singular-locus} we have $\m{dim}\ \m{Sing}(\widebar\scX_\vartriangle/\widebar\scS;d + 1) \leq d$, where \[\m{Sing}(\widebar\scX_\vartriangle/\widebar\scS;d + 1) \subset \widebar\scX_\vartriangle\]
		is the reduced closed subset cut out by the radical of the $(d + 1)$-th Fitting ideal $\m{Fitt}_{d + 1}(\Omega^1_{\widebar\scX_\vartriangle/\widebar\scS})$ of the sheaf of relative K\"ahler differential forms $\Omega^1_{\widebar\scX_\vartriangle/\widebar\scS}$. \qedex
	\end{enumerate}
\end{ass}

By \cite[Lemma~2.10]{Felten-gls-gr}, this implies that $\phi_\vartriangle\colon \widebar\scX_\vartriangle \to \widebar\scS$ is flat, and that its fibers are Cohen--Macaulay, pure of dimension $d + 1$, and geometrically reduced. Furthermore, $\phi_\vartriangle$ is dominant.

Note that Condition~\ref{cond:integral} guarantees that $\widebar\scX_\vartriangle$ is integral, but it does not guarantee that $\widebar X_\vartriangle$ is integral.

Condition~\ref{cond:injective} implies that $R_{\widebar\scS} \to R_{\widebar\scX}^0$ is injective so that it is an isomorphism because $R_\PP^0 = \kq[t]$. Therefore, the augmentation homomorphism $R_{\widebar\scX}^\bullet \to R_{\widebar\scX}^0 = R_{\widebar\scS}$ defines a closed immersion $\eps\colon \widebar\scS \to \widebar\scX_\vartriangle$. We set 
\[\widebar\scX_\punct = \widebar\scX_\vartriangle \setminus \eps(\widebar\scS) = \widebar\scX_\vartriangle \setminus (\{0\} \times \widebar\scS),\]
the pointed affine cone, and denote its inclusion by $\iota\colon \widebar\scX_\punct \to \widebar\scX_\vartriangle$. Since $\kq \to \kk$ is flat, also $R_{\widebar S} \to R_{\widebar X}^\bullet$ is injective, and $R_{\widebar S} \to R_{\widebar X}^0$ is an isomorphism. We obtain a closed immersion $\eps\colon \widebar S \to \widebar X_\vartriangle$ and set $\widebar X_\punct = \widebar X_\vartriangle \setminus \eps(\widebar S)$. We write $\iota\colon \widebar X_\punct \to \widebar X_\vartriangle$ for the open immersion.

As in \cite[Sitn.~2.18]{Felten-gls-gr}, we have affine smooth surjective morphisms $\pi\colon \widebar\scX_\punct \to \widebar\scX$ and $\pi\colon \widebar X_\punct \to \widebar X$ of relative dimension $1$. Then, $\phi\colon \widebar\scX \to \widebar\scS$ and $\phi\colon \widebar X \to \widebar S$ are flat and surjective, and their fibers are Cohen--Macaulay, pure of dimension $d$, and geometrically reduced.

\begin{ass}\label{ass:smooth-fiber}\note{ass:smooth-fiber}
	We assume that $\m{Sing}(\widebar\scX_\vartriangle/\widebar\scS;d + 1) \cap \widebar\scX_\punct \cap \phi_\vartriangle^{-1}(s) = \varnothing$ for some closed point $s \in \widebar\scS$. \qedex
\end{ass}

After setting $\scX_s = \phi^{-1}(s) \subset \widebar\scX$, the fiber $\scX_s \to \Spec \kappa(s)$ is smooth. Since $\phi\colon \widebar\scX \to \widebar\scS$ is proper,
\[\m{Disc}(\phi) \coloneq \phi(\m{Sing}(\widebar\scX/\widebar\scS;d)) \subset \widebar\scS\]
is closed, and since $s \notin \m{Disc}(\phi)$, it is a proper closed subset.

We write $\scS \coloneq \widebar\scS \setminus (\m{Disc}(\phi) \setminus \{0\})$ and $\scS^\star = \scS \setminus \{0\}$, and we set $\scX = \phi^{-1}(\scS)$ and $\scX^\star = \phi^{-1}(\scS^\star)$. Then, $\scV \coloneq \scX_0$ is the only singular fiber of $\phi\colon \scX \to \scS$. Similarly, $V \coloneq X_0$ is the only singular fiber of $\phi\colon X \to S$, where $S$ is the base change of $\scS$ along $\Spec \kk \to \Spec \kq$. 

The purpose of assuming that $\phi\colon X^\star \to S^\star$ is smooth is that this implies that the Hodge--de Rham spectral sequence of $\phi\colon X^\star \to S^\star$ degenerates on the page $E_1$. We will use this in our analysis of the Hodge theory in the central fiber, cf.~Corollary~\ref{cor:cohom-flat-criterion}.

\subsection{The geometry of the central fiber $\scV_\vartriangle = \scX_{\vartriangle,0}$}\label{sect:V-geometry}\note{sect:V-geometry}

We continue in the setup of the previous section. We write $\scV^\circ_\vartriangle \subseteq \scV_\vartriangle$ for the smooth locus, and $\scD_{\vartriangle} \subset \scV_{\vartriangle}$ for its complement with the reduced induced scheme structure. This subset is cut out by the radical of the $(d + 1)$-th Fitting ideal $\m{Fitt}_{d + 1}(\Omega^1_{\scV_\vartriangle/\kq})$.

\begin{defn}[\protect{cf.~\cite[Defn.~6.6]{Felten-gls-gr}}]\label{defn:intrinsic-dnc}\note{defn:intrinsic-dnc}
	Let $\nu_\vartriangle\colon \widetilde\scV_\vartriangle \to \scV_\vartriangle$ be the normalization, let $\scC_\vartriangle \to \scV_\vartriangle$ be the conductor locus, and let $\mu_\vartriangle\colon \widetilde\scC_\vartriangle \to \scC_\vartriangle$ be the base change of $\nu_\vartriangle$ along $\scC_\vartriangle \to \scV_\vartriangle$. An open subset $\scW \subseteq \scV_\vartriangle$ is \emph{intrinsically double normal crossing} if
	\begin{enum:alph}
		\item $\widetilde\scW \coloneq \nu_\vartriangle^{-1}(\scW) \subseteq \widetilde\scV_\vartriangle$ is smooth;
		\item $\nu_\vartriangle\colon \widetilde\scW \to \scW$ is unramified;
		\item $\scC_\vartriangle \cap \scW$ is smooth and of dimension $d$;
		\item $\mu_{\vartriangle,*}\cO_{\widetilde\scC_\vartriangle}|_{\scW}$ is a free $\cO_{\scC_\vartriangle \cap \scW}$-module of rank $2$. \qeddef
	\end{enum:alph}
\end{defn}

We write 
\[\scV_\vartriangle^\times = \m{Dnc}(\scV_\vartriangle) \subseteq \scV_\vartriangle\]
for the \emph{(intrinsic) double normal crossing locus}, i.e., the largest open subset which is intrinsically double normal crossing. 
Note that, somewhat counterintuitive, every smooth point $p \in \scV_\vartriangle$ is 
intrinsically double normal crossing according to this definition. The characterizing 
property is hence to be understood as the ``worst case scenario.''
By the results of \cite[\S6]{Felten-gls-gr}, the subset $V_\vartriangle^\times = \scV_\vartriangle^\times \otimes_\kq \kk \subseteq V_\vartriangle$ is precisely the set of points where $V_\vartriangle$ is \'etale locally of the form $\Spec \kk[x,y,z_1,\ldots,z_d]/(xy)$, i.e., a double normal crossing scheme in the ordinary sense. Following \cite[Nota.~2.5]{Felten-gls-gr}, we write 
\[\scT_\vartriangle^\times = \m{nDnc}(\scV_\vartriangle) \coloneq \scV_\vartriangle \setminus \scV_\vartriangle^\times \subset \scV_\vartriangle \quad \text{and} \quad \scD_\vartriangle^\circ = \m{pDnc}(\scV_\vartriangle) \coloneq \scV_\vartriangle^\times \setminus \scV_\vartriangle^\circ.\]
When endowed with the reduced induced scheme structure, $\scD_\vartriangle^\circ$ is smooth of dimension $d$. We denote the set of irreducible (= connected) 
components of $\scD_\vartriangle^\circ$ by $[\scD_\vartriangle^\circ]$.

Under the projection $\pi_0\colon \scV_\punct \to \scV$, we have $\pi_0^{-1}(\scV^\circ) = \scV_\punct^\circ$ and $\pi_0^{-1}(\scD) = \scD_\punct$ and $\pi_0^{-1}(\scV^\times) = \scV^\times_\punct$ and $\pi_0^{-1}(\scD^\circ) = \scD_\punct^\circ$ and $\pi_0^{-1}(\scT^\times) = \scT_\punct^\times$ and $\pi_0^{-1}(\scT) = \scT_\punct$ for the analogous subsets of $\scV$ and $\scV_\punct$. The formation of all these subsets commutes with base change along $\Spec \kk \to \Spec \kq$.

We have bijections between $[\scD^\circ]$, $[\scD_\punct^\circ]$, and $[\scD_\vartriangle^\circ]$.
\medskip

In this article, we shall apply two strategies to compute $\scT_\vartriangle^\times$.

\begin{strat}\label{strat:Tx-sing}\note{strat:Tx-sing}
	For the first strategy, let $\scT_\vartriangle \subset \scD_\vartriangle$ be the closed subset cut out by the radical of the $d$-th Fitting ideal $\m{Fitt}_d(\Omega^1_{\scD_\vartriangle/\kq})$. When assuming that every irreducible component of $\scD_\vartriangle$ has dimension $d$, $\scT_\vartriangle$ is the reduced singular locus of $\scD_\vartriangle$. In particular, we have $\scT_\vartriangle \cap \scV_\vartriangle^\times = \varnothing$ since $(\scT_\vartriangle \cap \scV_\vartriangle^\times) \otimes_\kq \kk \cong \m{Sing}(D_\vartriangle) \cap V_\vartriangle^\times$, and since the reduced singular locus of $\Spec \kk[x,y,z_1,\ldots,z_{d - 1}]/(xy)$ is itself smooth. Therefore, we have $\scT_\vartriangle \subset \scT_\vartriangle^\times$. Now, if $\scV_\vartriangle \setminus \scT_\vartriangle$ is intrinsically double normal crossing, then $\scT_\vartriangle = \scT_\vartriangle^\times$. \qedex
\end{strat}

\begin{strat}\label{strat:Tx-unramified}\note{strat:Tx-unramified}
	For the second strategy, let $\scA_\vartriangle$ be the support of $(\nu_\vartriangle)_*\Omega^1_{\widetilde\scV_\vartriangle/\scV_\vartriangle}$ with its reduced induced scheme structure. Again, we have $\scA_\vartriangle \cap \scV_\vartriangle^\times = \varnothing$ so that $\scA_\vartriangle \subset \scT_\vartriangle^\times$. If $\scV_\vartriangle \setminus \scA_\vartriangle$ is intrinsically double normal crossing, then $\scA_\vartriangle = \scT_\vartriangle^\times$. \qedex
\end{strat}

After computing $\scT_\vartriangle^\times$ with any of the two strategies, the following assumption is computationally verifiable.

\begin{ass}\label{ass:Tx-small}\note{ass:Tx-small}
	We assume $\m{dim}(\scT_\vartriangle^\times) \leq d - 1$. \qedex
\end{ass}

\noindent Under this assumption, we can compute the set of irreducible components of $\scD_\vartriangle^\circ$ as
\[[\scD_\vartriangle^\circ] = \{\scE \setminus \scT^\times_\vartriangle \ | \ \scE \in \m{Irr}(\scD_\vartriangle), \ \m{dim}(\scE) = d\}.\]

\subsection{The log structure}

We continue in the setup of the previous sections. In this section, we enhance $\phi\colon X \to S$ to a morphism of log schemes $\phi\colon X^\dagger \to S^\dagger$ and introduce Assumption~\ref{ass:agls}, which allows us to find an open subset $U \subseteq X$ such that, after restricting the log structure to $U$, the induced morphism of partial log schemes $\phi\colon X^\partial \to S^\dagger$ is generically log smooth. This is achieved via an algebraic criterion that (when satisfied at a point $p \in X_0$) guarantees that $\phi\colon X \to S$ is \'etale locally of the form $xy = t^\ell$ (near the point $p \in X_0$).

We construct the log structures. The divisor $\{0\} = \{t = 0\} \subset S$ defines a log structure $\cM_S$ in the \'etale topology of $S$. The resulting log scheme $(S,\cM_S) = (S|\{t = 0\})$ is canonically an open subset of $\Spec(\NN \to \kk[t])$. The divisor $V = \{\phi^\sharp(t) = 0\} \subset X$ defines a log structure $\cM_X$ in the \'etale topology of $X$. We obtain an induced morphism of log schemes $\phi\colon (X,\cM_X) \to (S,\cM_S)$, which we normally write as $\phi\colon X^\dagger \to S^\dagger$.

To find an open subset $U \subseteq X$ such that $\phi\colon (X,U,\cM_U) \to (S,\cM_S)$ is generically log smooth, we use the algebraic criterion developed in \cite{Felten-gls-gr}. To this end, we need some notation. For $k \geq 0$, we set $\scS_k = \Spec \kq[t]/(t^{k + 1})$ and $\scX_{\vartriangle,k} = \scX_\vartriangle \times_\scS \scS_k$. We write \[\varrho\colon \enspace \Omega^1_{\widebar\scX_\vartriangle} \to \Omega^1_{\widebar\scX_\vartriangle/\widebar\scS}\]
for the canonical quotient map, and we set
\[\xi_k \coloneq \cE xt^1(\varrho|_{\scX_{\vartriangle,k}},\cO_{\scX_{\vartriangle,k}})\colon \qquad \cE xt^1(\Omega^1_{\scX_{\vartriangle,k}/\scS_{k}},\cO_{\scX_{\vartriangle,k}}) \to \cE xt^1(\Omega^1_{\widebar\scX_\vartriangle}|_{\scX_{\vartriangle,k}},\cO_{\scX_{\vartriangle,k}}).\]
We denote the set-theoretic support of its kernel by $\scK_{\vartriangle}^{(k)} \subset \scV_\vartriangle$, considered as a reduced closed subscheme.

For $\ell \geq 1$, we write 
\[\sigma_{\ell - 1}\colon \quad \cE xt^1(\Omega^1_{\widebar\scX_\vartriangle}|_{\scX_{\vartriangle,\ell - 1}},\cO_{\scX_{\vartriangle,\ell - 1}}) \xrightarrow{(-)\cdot \phi^\sharp(t)^{\ell - 1}}\cE xt^1(\Omega^1_{\widebar\scX_\vartriangle}|_{\scX_{\vartriangle,\ell - 1}},\cO_{\scX_{\vartriangle,\ell - 1}})\]
for the multiplication map. We denote the set-theoretic support of its image by $\scZ_{\vartriangle}^{(\ell)} \subset \scV_\vartriangle$, considered as a reduced closed subscheme. Now, the following assumption guarantees that we can use the criterion from \cite{Felten-gls-gr}:

\begin{ass}\label{ass:agls}\note{ass:agls}
	We assume for every $\scE \in [\scD_\vartriangle^\circ]$ that there is some $\ell = \ell(\scE) \geq 1$ with:
	\begin{enum:alph}
		\item $\scE \cap \scK_{\vartriangle}^{(k)} = \varnothing$ for all $0 \leq k \leq \ell - 2$ (this condition is empty if $\ell = 1$);
		\item $\m{dim}(\scE \cap \scZ_{\vartriangle}^{(\ell)}) \leq d - 1$. \qedex
	\end{enum:alph}
\end{ass}

The number $\ell = \ell(\scE)$ is unique if it exists, and it is called the \emph{kink} of $\scE$. For $\scE \in [\scD_\vartriangle^\circ]$, we set $Z_{\vartriangle}^{\scE} \coloneq (\scE \otimes_\kq \kk) \cap (\scZ_{\vartriangle}^{(\ell)} \otimes_\kq \kk)$ with $\ell = \ell(\scE)$. Then, we set
\begin{equation}\label{eqn:clubsuit}\note{eqn:clubsuit}
	Z_{\vartriangle} = T^\times_\vartriangle \cup Z_{\vartriangle}^\circ \quad \text{with} \quad Z_{\vartriangle}^\circ = \bigsqcup_{\scE \in [\scD_\vartriangle^\circ]} Z_{\vartriangle}^{\scE} \tag{$\clubsuit$}.
\end{equation} 
This is a closed subset of $V_\vartriangle$ of codimension $\geq 2$. We write $Z_\punct = Z_\vartriangle \cap V_\punct$. Since $\phi_\punct\colon X_\punct^\star \to S^\star$ is smooth, it is a straightforward consequence of \cite[Prop.~2.19]{Felten-gls-gr} that there is a (necessarily unique) closed subset $Z \subset V$ of codimension $\geq 2$ with $\pi_0^{-1}(Z) = Z_\punct$, and such that $\phi\colon (X,U,\cM_U) \to (S,\cM_S)$ is generically log smooth with $U = X \setminus Z$.

We normally write $\phi\colon X^\partial \to S^\dagger$ for this generically log smooth family.

Let $E \in [D^\circ]$, and let $p \in E \setminus Z$. Then, an \'etale local model for $\phi\colon X \to S$ near $p$ is given by
\[\Spec \kk[t,x,y,z_1,\ldots,z_{d - 1}]/(xy - t^\ell) \to \Spec \kk[t]\] for $\ell = \ell(E)$. When $p \in Z \cap V^\times$, one can show that $\phi\colon X^\dagger \to S^\dagger$ is indeed not log smooth at $p$ \cite[Thm.~2.7]{Felten-gls-gr}, but if $p \in T^\times$, then $\phi\colon X^\dagger \to S^\dagger$ might actually be log smooth at $p$. For a typical example, see \cite[Rem.~2.17]{Felten-gls-gr}.

\subsection{The log de Rham complex}\label{sect:log-de-Rham}\note{sect:log-de-Rham}

We continue in the setup of the previous sections. Let $k\colon \widebar\scX^\circ \to \widebar\scX$ be the inclusion of the smooth locus of $\phi\colon \widebar\scX \to \widebar\scS$. Then, we define the \emph{algorithmic log de Rham complex} as
\begin{equation}\label{eqn:A-bullet-scX}\note{eqn:A-bullet-scX}
	\cA^{\bullet,\m{algo}}_{\widebar\scX/\widebar\scS} \coloneq k_*\Omega^\bullet_{\widebar\scX^\circ/\widebar\scS},
\end{equation}
the direct image of the ordinary relative de Rham complex of $\widebar\scX^\circ \to \widebar\scS$. This is a first-order-differential-operator complex in the sense of Definition~\ref{defn:coh-flat}.

Since the fibers of $\phi\colon \widebar\scX \to \widebar\scS$ are pure of dimension $d$, we have $\widebar\scX \setminus \widebar\scX^\circ = \m{Sing}(\widebar\scX/\widebar\scS;d)$, so Part~\ref{cond:singular-locus} of Assumption~\ref{ass:apgls} together with $\widebar\scX$ being Cohen--Macaulay guarantee that $\cO_{\widebar\scX} \to k_*\cO_{\widebar\scX^\circ}$ is an isomorphism. Therefore, we have
\[\cA^{i,\m{algo}}_{\widebar\scX/\widebar\scS} \cong \cH om(\cH om(\Omega^i_{\widebar\scX/\widebar\scS},\cO_{\widebar\scX}),\cO_{\widebar\scX}),\]
and $\cA^{i,\m{algo}}_{\widebar\scX/\widebar\scS}$ is coherent and reflexive.

Since $\kq[t]$ is a Dedekind domain, a $\kq[t]$-module is flat if and only if it is torsion-free. Let $0 \not= f \in \kq[t]$, and consider the multiplication map $\mu_f\colon \cA^{i,\m{algo}}_{\widebar\scX/\widebar\scS} \to \cA^{i,\m{algo}}_{\widebar\scX/\widebar\scS}$. Since $\phi\colon \widebar\scX \to \widebar\scS$ is flat and $\Omega^i_{\widebar\scX^\circ/\widebar\scS}$ is finite locally free, $\mu_f$ is injective on $\widebar\scX^\circ$. But then, $\mu_f$ is injective, and $\cA^{i,\m{algo}}_{\widebar\scX/\widebar\scS}$ is flat over $\widebar\scS$.

Since $\Omega^1_{\widebar\scX^\circ/\widebar\scS}$ is locally free of rank $d$, the complex $\cA^{\bullet,\m{algo}}_{\widebar\scX/\widebar\scS}$ lies in degrees $0 \leq i \leq d$.

Let $\cA^{\bullet,\m{algo}}_{\widebar X/\widebar S} \coloneq \cA^{\bullet,\m{algo}}_{\widebar\scX/\widebar\scS} \otimes_\kq \kk$. By \cite[Prop.~C.19]{FeltenGLDT}, we can in fact perform the base change of the differential of $\cA^{\bullet,\m{algo}}_{\widebar\scX/\widebar\scS}$. Since $\Spec \kk \to \Spec \kq$ is flat, we have $\cA^{\bullet,\m{algo}}_{\widebar X/\widebar S} \cong k_*\Omega^\bullet_{\widebar X^\circ/\widebar S}$ for the inclusion $k\colon \widebar X^\circ \to \widebar X$ of the smooth locus of $\phi\colon \widebar X \to \widebar S$. After restriction to $S \subseteq \widebar S$, we have $X^\circ \subseteq U$, and the natural map $\Omega^\bullet_{X^\circ/S} \to \Omega^{\bullet,\m{log}}_{X^\circ/S}$ is an isomorphism because $(X^\circ,\cM_{X^\circ}) \to (S,\cM_S)$ is strict and smooth. But then, we have
\[\cA^{\bullet,\m{algo}}_{X/S} \cong k_*\Omega^\bullet_{X^\circ/S} \cong k_*\Omega^{\bullet,\m{log}}_{X^\circ/S} \cong j_*\Omega^{\bullet,\m{log}}_{U/S} \eqcolon \cW^\bullet_{X/S}\]
for the inclusion $j\colon U \to X$. Indeed, we have $k_*\Omega^{i,\m{log}}_{X^\circ/S} \cong j_*\Omega^{i,\m{log}}_{U/S}$ because $\Omega^{i,\m{log}}_{U/S}$ is finite locally free (on $U$).

\subsection{The log canonical bundle}\label{sect:log-canonical}\note{sect:log-canonical}

We continue in the setup of the previous sections. Recall that, for every separated morphism of finite type $f\colon Y \to T$ between Noetherian schemes, we have a functor $f^!\colon \m{D}^+_\m{qc}(\cO_T) \to \m{D}^+_\m{qc}(\cO_Y)$. Following \cite{Verdier1969,Sastry2004}, we have the \emph{Verdier isomorphism}
\begin{equation*}\label{eqn:Verdier-isom}\note{eqn:Verdier-isom}
	\Omega^d_{\widebar\scX^\circ/\widebar\scS} \cong \cH^{-d}(\phi'^!\cO_{\widebar\scS})
\end{equation*}
for the smooth morphism $\phi'\colon \widebar\scX^\circ \to \widebar\scS$. Note that, unlike the more standard version which can be found for instance in \cite[Tag~0AU3]{stacks}, we do not need properness of $\phi'\colon \widebar\scX^\circ \to \widebar\scS$.

\begin{ass}\label{ass:Gorenstein}\note{ass:Gorenstein}
	The scheme $\widebar\scX$ is Gorenstein. \qedex
\end{ass}

Then, $\phi\colon \widebar\scX \to \widebar\scS$ is a Gorenstein morphism in the sense of \cite[Tag~0C04]{stacks}, and $\cH^{-d}(\phi^!\cO_{\widebar\scS})$ is a line bundle on $\widebar\scX$ which extends $\Omega^d_{\widebar\scX^\circ/\widebar\scS}$. In particular, we have $\cA^{d,\m{algo}}_{\widebar\scX/\widebar\scS} \cong \cH^{-d}(\phi^!\cO_{\widebar\scS})$, and $\cW^d_{X/S}$ is a line bundle. In other words, $\phi\colon X^\partial \to S^\dagger$ is log Gorenstein in the sense of \cite[Defn.~8.4]{FeltenGLDT}.

In all our examples, we shall use the following result to verify Assumption~\ref{ass:Gorenstein}.

\begin{lem}\label{lem:log-canonical-ci}\note{lem:log-canonical-ci}
	Suppose that Assumptions~\ref{ass:apgls}, \ref{ass:smooth-fiber}, \ref{ass:Tx-small}, and \ref{ass:agls} are satisfied. Assume that $n = d + r$ for $n,d,r$ as in Datum~\ref{dat:basic-input}. Then, $\widebar\scX$ is Gorenstein, and moreover, we have an isomorphism 
	\[\cA^{d,\m{algo}}_{\widebar\scX/\widebar\scS} \cong \cO_{\widebar\scX}(-n - 1 + \delta)\]
	for $\delta = \sum_{j = 1}^r \delta_j = \sum_{j = 1}^r \m{deg}(f_j)$.
\end{lem}
\begin{proof}
	Let $i\colon \widebar\scX \to \PP^n \times \widebar\scS$ be the closed immersion, and let $\cI \subset \cO_{\PP^n \times \widebar\scS}$ be the corresponding ideal. At every closed point $p \in \widebar\scX$, the ideal $\scI_p$ is generated by $r$ elements, and $\m{dim}(\cO_{\widebar\scX,p}) = d + 1$. Thus, by \cite[Tag~02JN]{stacks}, the generators of $\cI_p$ are a regular sequence, and $i\colon \widebar\scX \to \PP^n \times \widebar\scS$ is a regular immersion in the sense of \cite[Tag~063J]{stacks}. By \cite[Tag~0DWA]{stacks}, $\widebar\scX_\vartriangle$ is Gorenstein because $\PP^n \times \widebar\scS$ is regular, hence Gorenstein.
	
	We consider the conormal sheaf $\cC_{\widebar\scX/\PP^n \times \widebar\scS} \coloneq \cI/\cI^2$. We have a factorization of $i$ into a sequence of closed immersions
	\[\widebar\scX = \widebar\scX_{(r)} = \{f_1 = \ldots = f_r = 0\} \subset \ldots \subset \widebar\scX_{(1)} = \{f_1 = 0\} \subset \widebar\scX_{(0)} = \PP^n \times \widebar\scS.\]
	For $p \in \widebar\scX$, each of these closed immersions is regular at $p$. Let $\delta_i = \m{deg}(f_i)$, and let $\cI_i$ be the ideal corresponding to the closed immersion $\widebar\scX_{(i)} \subset \widebar\scX_{(i - 1)}$. We have $\cN_{\widebar\scX_{(i)}/\widebar\scX_{(i - 1)}} = \cI_i/\cI_i^2 \cong \cI_i \otimes \cO_{\widebar\scX_i}$, and so the canonical surjection $\cO_{\PP^n \times \widebar\scS}(-\delta_i) \otimes \cO_{\widebar\scX_{(i - 1)}} \to \cI_i$ yields a canonical surjection $\cO_{\PP^n \times \widebar\scS}(-\delta_i) \otimes \cO_{\widebar\scX_{(i)}} \to \cN_{\widebar\scX_{(i)}/\widebar\scX_{(i-1)}}$. At $p \in \widebar\scX$, this is an isomorphism. Now, \cite[Tag~063N]{stacks} applied inductively allows us to conclude that $\cC_{\widebar\scX/\PP^n \times \widebar\scS}$ is locally free of rank $r$, and that $\bigwedge^r\cC_{\widebar\scX/\PP^n \times \widebar\scS} \cong \cO_{\widebar\scX}(-\delta)$, where $\delta = \sum_{i = 1}^r\delta_i$.
	
	We have $\cA^{d,\m{algo}}_{\widebar\scX/\widebar\scS} \cong \cH^{-d}(\phi^!\cO_{\widebar\scS})$. Let $\pi\colon \PP^n \times \widebar\scS \to \widebar\scS$ be the second projection, and let $\cN_{\widebar\scX/\PP^n \times \widebar\scS} = \cH om(\cC_{\widebar\scX/\PP^n \times \widebar\scS},\cO_{\widebar\scX})$ be the normal sheaf of the inclusion. Then, we have
	\begin{equation*}
		\phi^!\cO_{\widebar\scS} \cong i^!\pi^!\cO_{\widebar\scS} \cong Li^*(\Omega^n_{\PP^n \times \widebar\scS/\widebar\scS}[n]) \otimes \bigwedge^r \cN_{\widebar\scX/\PP^n \times \widebar\scS}[-r] \cong \cO_{\widebar\scX}(- n - 1 + \delta)[d]
	\end{equation*}
	by \cite[Tag~0AU3]{stacks}. In particular, $\cA^{d,\m{algo}}_{\widebar\scX/\widebar\scS} \cong \cO_{\widebar\scX}(-n - 1 + \delta)$.
\end{proof}

\section{Homological algebra}\label{sect:homological-algebra}\note{sect:homological-algebra}

In this section, we develop the homological algebra necessary to prove Proposition~\ref{prop:numerical-reduction}.

\subsection{Prespectral sequences}

Let $S$ be a scheme. When $E$ is a spectral sequence in $\cO_S$-modules, and $b\colon T \to S$ is a morphism of schemes, then it is not immediately clear what the \emph{pullback spectral sequence} $b^*E$ should be. To make sense of the notation $b^*E$, we introduce the following slight weakening of the notion of a spectral sequence.

\begin{defn}
	Let $\bA$ be an Abelian category. A (right) \emph{prespectral sequence} is a tuple 
	\[(E_r^{p,q},d_r^{p,q},\phi_r^{p,q})\] 
	consisting of:
	\begin{enum:alph}
		\item for every $r \geq 1$ and $p,q \in \ZZ$, an object $E_r^{p,q}$ of $\bA$;
		\item for every $r \geq 1$ and $p,q \in \ZZ$, a morphism $d_r^{p,q}\colon E_r^{p,q} \to E_r^{p + r,q +1 - r}$ in $\bA$ such that $d_r^{p,q} \circ d_r^{p - r, q + r - 1} = 0$;
		\item for every $r \geq 1$ and $p,q \in \ZZ$, a morphism 
		\[\phi_r^{p,q}\colon \enspace E_{r + 1}^{p,q} \to \frac{\m{ker}(d_r^{p,q}\colon E_r^{p,q} \to E_r^{p + r,q + 1 - r})}{\m{im}(d_r^{p - r, q + r - 1}\colon E_r^{p - r,q + r - 1} \to E_r^{p,q})}\]
		in $\bA$.
	\end{enum:alph}
	A \emph{spectral sequence} is a prespectral sequence such that every $\phi_r^{p,q}$ is an isomorphism. A prespectral sequence \emph{degenerates} at $E_1$ if $d_r^{p,q} = 0$ for all $r \geq 1$. A \emph{homomorphism} between prespectral sequences $D$ and $E$ consists of a morphism $h_r^{p,q}\colon D_r^{p,q} \to E_r^{p,q}$ in $\bA$ for every $r \geq 1$ and every $p,q \in \ZZ$. These must commute with $d_r$ and $\phi_r$. \qeddef
\end{defn}

\begin{lem}
	Let $F\colon \bA \to \bB$ be an additive and right exact functor between Abelian categories, and let $E$ be a prespectral sequence in $\bA$. Then:
	\begin{enum:arabic}
		\item $F(E)$ is a prespectral sequence in $\bB$; 
		\item if $E$ degenerates at $E_1$, then $F(E)$ degenerates at $E_1$;
		\item if $E$ is a spectral sequence and $F$ is exact, then $F(E)$ is a spectral sequence.
	\end{enum:arabic}
\end{lem}
\begin{proof}
	Let $\cp\colon A \to B$ and $\psi\colon B \to C$ be morphisms in $\bA$ with $\psi \circ \cp = 0$. Let $K \coloneq \m{ker}(\psi)$, and let $\chi\colon A \to K$ be the induced map. Let $H \coloneq \m{coker}(\chi)$ so that $H = \m{ker}(\psi)/\m{im}(\cp)$. Since $F$ is right exact, we have $F(H) = \m{coker}(F\chi)$. Because we have a natural morphism $F(K) \to \m{ker}(F\psi)$, there is a natural morphism $\m{coker}(F\chi) \to \m{coker}(F(A) \to \m{ker}(F\psi))$ so that we obtain a natural morphism
	\[ F(\m{ker}(\psi)/\m{im}(\cp))  \to \m{ker}(F\psi)/\m{im}(F\cp) .\]
	If $F$ is exact, then also $F(K) \to \m{ker}(F\psi)$ is an isomorphism so that the morphism in the display is an isomorphism.
\end{proof}

\subsection{Spectral sequences from double complexes}

Let $\bA$ be an Abelian category. Our convention for a \emph{double complex} $(I^{\bullet,\bullet},\partial,\bar\partial)$ in $\bA$ (with $\partial\colon I^{p,q} \to I^{p + 1,q}$ and $\bar\partial\colon I^{p,q} \to I^{p,q + 1}$) is that $\partial\bar\partial + \bar\partial\partial = 0$. When $I^{\bullet,\bullet}$ is bounded below in both indices, we filter $I^{\bullet,\bullet}$ as $F^iI^{p,q} = I^{p,q}$ if $p \geq i$ and $F^iI^{p,q} = 0$ if $p < i$. This gives rise to a spectral sequence ${'E} = {'E}(I^{\bullet,\bullet})$ associated with $
I^{\bullet,\bullet}$, for which we have 
\[{'E}_1^{p,q} = H^q(I^{p,\bullet},\bar\partial) \quad \text{and} \quad {'E}_1^{p,q} \Rightarrow H^{p + q}(\m{Tot}(I^{\bullet,\bullet}),\partial + \bar\partial).\] 
We are primarily interested in the case where $\bA = \m{Mod}(R_S)$ is the category of modules over a ring $R_S$. In this case, let us review the construction of the spectral sequence ${'E}(I^{\bullet,\bullet})$.

\begin{constr}
	Let $R_S$ be a ring, and let $(I^{\bullet,\bullet},\partial,\bar\partial)$ be a double complex in $\m{Mod}(R_S)$, bounded below in both indices. First, let \[Z_1^{p,q} = \{\alpha \in I^{p,q} \ | \ \bar\partial\alpha = 0\} \quad \text{and} \quad B_1^{p,q} = \{\bar\partial \beta \in I^{p,q} \ | \ \beta \in I^{p,q - 1}\}\] 
	so that ${'E}_1^{p,q} \coloneq Z_1^{p,q}/B_1^{p,q} = H^q(I^{p,\bullet},\bar\partial)$. We define an $R_S$-linear map $d_1\colon {'E}_1^{p,q} \to {'E}_1^{p + 1,q}$ by setting $d_1[\alpha] \coloneq [\partial\alpha]$ for $\alpha \in Z_1^{p,q}$ representing $[\alpha] \in {'E}_1^{p,q}$. This satisfies $d_1 \circ d_1 = 0$.
	
	Now, for $r \geq 1$ and $p,q \in \ZZ$, we define an \emph{$r$-zigzag starting at $(p,q)$} as a sequence $\alpha_\bullet = (\alpha_1,\ldots,\alpha_r)$ with $\alpha_i \in I^{p - 1 + i,q + 1 - i}$ such that $\bar\partial\alpha_1 = 0$ and $\bar\partial\alpha_{i + 1} + \partial\alpha_i = 0$ for $1 \leq i \leq r - 1$. Let $\m{Zig}_r^{p,q}$ be the $R_S$-module of $r$-zigzags starting with $\alpha_1 \in I^{p,q}$. With this definition, we set 
	\[Z_r^{p,q} \coloneq \{\alpha \in I^{p,q} \ | \ \exists\ r\text{-zigzag}\ (\alpha,\alpha_2,\ldots,\alpha_r)\} = \m{im}(\m{Zig}_r^{p,q} \to I^{p,q}, \: \alpha_\bullet \mapsto \alpha_1).\]
	By truncating zigzags, we see that $Z_{r + 1}^{p,q} \subseteq Z_r^{p,q}$. 
	
	For $r \geq 2$, we define $B_r^{p,q} \subseteq I^{p,q}$ such that $\alpha \in B_r^{p,q}$ if and only if we can find $\beta^0 \in I^{p,q - 1}$ and a sequence of zigzags $\beta^1_\bullet = (\beta^1_1)$, $\beta^2_\bullet = (\beta^2_1,\beta^2_2)$, \ldots, $\beta^{r - 1}_\bullet = (\beta^{r - 1}_1,\ldots,\beta^{r - 1}_{r - 1})$ with $\beta^1_1,\ldots,\beta^{r - 1}_{r - 1} \in I^{p - 1,q}$ such that 
	\[\alpha = \bar\partial\beta^0 + \partial\beta^1_1 + \partial\beta^2_2 + \ldots + \partial\beta^{r - 1}_{r - 1}.\]
	In other words, we have 
	\[B_r^{p,q} = \partial I^{p,q - 1} + \sum_{i = 1}^{r - 1} \m{im}(\m{Zig}_i^{p - i,q + i - 1} \to I^{p,q},\: \beta_\bullet \mapsto \partial\beta_i).\]
	We clearly have $B_r^{p,q} \subseteq B_{r + 1}^{p,q}$, and it is easy to see that also $B_r^{p,q} \subseteq Z_s^{p,q}$ for any $r,s \geq 1$.
	
	For $r \geq 2$, we set ${'E}_r^{p,q} \coloneq Z_r^{p,q}/B_r^{p,q}$. Then, we define an $R_S$-linear map 
	\[d_r\colon {'E}_r^{p,q} \to {'E}_r^{p + r,q + 1 - r}\] 
	as follows: For $\alpha \in Z_r^{p,q}$ representing $[\alpha] \in {'E}_r^{p,q}$, we choose an $r$-zigzag $\alpha_\bullet = (\alpha,\alpha_2,\ldots,\alpha_r)$, and then we set $d_r[\alpha] \coloneq [\partial\alpha_r]$. This is well-defined because the difference between two $r$-zigzags $(\alpha,\alpha_2,\ldots,\alpha_r)$ and $(\alpha,\alpha'_2,\ldots,\alpha_r')$ is an $(r - 1)$-zigzag. We have $\bar\partial\partial\alpha_r = -\partial\bar\partial\alpha_r = \partial^2\alpha_{r - 1} = 0$. Therefore, $(\partial\alpha_r,0,\ldots,0)$ is an $r$-zigzag, and we have $d_r \circ d_r = 0$.
	
	One may show that $d_r[\alpha] = 0$ if and only if $\alpha \in Z_{r + 1}^{p,q} \subseteq Z_r^{p,q}$, and that the image of $d_r\colon {'E}_r^{p - r,q + r - 1} \to {'E}_r^{p,q}$ is $B_{r + 1}^{p,q}/B_r^{p,q}$. This yields the isomorphism $\phi_r^{p,q}$ in the definition of a spectral sequence. 
	\qeddef
\end{constr}

The explicit description allows us to show that the construction of ${'E}(I^{\bullet,\bullet})$ behaves well with respect to base change along ring homomorphisms.

\begin{lem}\label{lem:spectral-seq-comp-map}
	Let $R_S$ be a ring, and let $(I^{\bullet,\bullet},\partial,\bar\partial)$ be a double complex in $\m{Mod}(R_S)$, bounded below in both indices. Let $\cp\colon R_S \to R_T$ be a ring homomorphism. Then, we have a comparison map
	\[\gamma_\cp\colon {'E}(I^{\bullet,\bullet}) \otimes_{R_S} R_T \to {'E}(I^{\bullet,\bullet} \otimes_{R_S} R_T)\]
	between the prespectral sequence on the left and the spectral sequence on the right. If $\cp$ is flat, then this is an isomorphism. When $\psi\colon R_T \to R_U$ is another ring homomorphism, then we have 
	\[\gamma_\psi \circ (\gamma_\cp \otimes_{R_T} R_U) = \gamma_{\psi \circ \cp}\colon \quad {'E}(I^{\bullet,\bullet}) \otimes_{R_S} R_U\to {'E}(I^{\bullet,\bullet} \otimes_{R_S} R_U).\] 
\end{lem}
\begin{proof}
	Let $\widetilde I^{\bullet,\bullet} \coloneq I^{\bullet,\bullet} \otimes_{R_S} R_T$, and let $\widetilde B_r^{p,q}$, $\widetilde Z_r^{p,q}$ etc.~be the $R_T$-modules constructed from $\widetilde I^{\bullet,\bullet}$. Then, we have $\cp$-homomorphisms $\m{Zig}_r^{p,q} \to \widetilde{\m{Zig}}_r^{p,q}$, and consequently $\cp$-homomorphisms $B_r^{p,q} \to \widetilde B_r^{p,q}$ and $Z_r^{p,q} \to \widetilde Z_r^{p,q}$ which are compatible with the inclusions. This induces $\cp$-homomorphisms ${'E}_r^{p,q} \to {'\widetilde E}_r^{p,q}$ which are compatible with $d_r$, $\tilde d_r$ and $\phi_r^{p,q}$, $\tilde\phi_r^{p,q}$. The desired homomorphism $\gamma_\cp$ of prespectral sequences is then given by adjunction and the homomorphism 
	\[H^n({'E}_r,d_r) \otimes_{R_S} R_T \to H^n({'E}_r \otimes_{R_S} R_T,d_r).\]
	Compatibility with composition of ring homomorphisms is clear from this construction. When $\cp\colon R_S \to R_T$ is flat, then the comparison homomorphism in the preceding display is an isomorphism so that ${'E}(I^{\bullet,\bullet}) \otimes_{R_S} R_T$ is a spectral sequence. Furthermore, $Z_1^{p,q} \otimes_{R_S} R_T \to \widetilde Z_1^{p,q}$ and $B_1^{p,q} \otimes_{R_S} R_T \to \widetilde B_1^{p,q}$ are isomorphisms so that $\gamma_\cp$ is an isomorphism on the page $E_1$. Isomorphy then propagates to all pages, and we have an isomorphism of spectral sequences.
\end{proof}

The explicit description also allows us to give a criterion for $E_1$-degeneration of the spectral sequence ${'E}(I^{\bullet,\bullet})$.

\begin{lem}\label{lem:degen-char}
	Let $R_S$ be a ring, and let $(I^{\bullet,\bullet},\partial,\bar\partial)$ be a double complex in $\m{Mod}(R_S)$, bounded below in both indices. Then the following statements are equivalent:
	\begin{enum:roman}
		\item\label{item:sseq-degen} the spectral sequence ${'E}(I^{\bullet,\bullet})$ degenerates at $E_1$;
		\item\label{item:hyper-surj} for all $p,q \in \ZZ$, the homomorphism $H^{p + q}(F^p\m{Tot}(I^{\bullet,\bullet}),\partial + \bar\partial) \to H^q(I^{p,\bullet},\bar\partial)$ is surjective;
		\item\label{item:hyper-inj} for all $n,p \in \ZZ$, the homomorphism $H^n(F^p\m{Tot}(I^{\bullet,\bullet}),\partial + \bar\partial) \to H^n(\m{Tot}(I^{\bullet,\bullet}),\partial + \bar\partial)$ is injective.
	\end{enum:roman}
\end{lem}
\begin{proof}
	First, we show the equivalence of \ref{item:sseq-degen} and \ref{item:hyper-surj}. Let us say that an $r$-zigzag is \emph{acyclic} if $\partial\alpha_r = 0$, and let $\m{AZig}^{p,q}_r \subseteq \m{Zig}^{p,q}_r$ be the subset of acyclic zigzags. Note that every acyclic $r$-zigzag can be arbitrarily extended by $\alpha_{r + 1} = \ldots = 0$, and let $\m{AZig}^{p,q}$ be the thus constructed union of acyclic $r$-zigzags for arbitrary $r$. Then, we have $\m{AZig}^{p,q} = Z^{p + q}(F^p\m{Tot}(I^{\bullet,\bullet}),\partial + \bar\partial)$ via $(\alpha_1,\alpha_2,\ldots) \mapsto \sum_{i \geq 1} \alpha_i$. The map $H^{p + q}(F^p\m{Tot}(I^{\bullet,\bullet}),\partial + \bar\partial) \to H^q(I^{p,\bullet},\bar\partial)$ in \ref{item:hyper-surj} is given by $(\alpha_1,\alpha_2,\ldots) \mapsto [\alpha_1]$.
	
	First, assume that the spectral sequence ${'E}(I^{\bullet,\bullet})$ degenerates at $E_1$. Then, we have $Z_r^{p,q} = Z_1^{p,q}$ for all $r \geq 1$. Thus, every $\alpha \in Z_1^{p,q}$ can be extended to an arbitrarily long $r$-zigzag. Since $I^{\bullet,\bullet}$ is bounded below in both indices, this implies that $\alpha$ can be extended to an acyclic zigzag $(\alpha,\alpha_2,\ldots)$. This acyclic zigzag represents a class in $H^{p + q}(F^p\m{Tot}(I^{\bullet,\bullet}),\partial + \bar\partial)$ which maps to $[\alpha] \in H^q(I^{p,\bullet},\bar\partial)$. Conversely, assume that the map in \ref{item:hyper-surj} is surjective. Then, for every $\alpha \in Z^q(I^{p,\bullet},\bar\partial) = Z_1^{p,q}$, we can find an acyclic zigzag $(\alpha_1,\alpha_2,\ldots)$ such that $\alpha_1 - \alpha \in B_1^{p,q}$. Thus, there is some $\beta \in I^{p,q- 1}$ with $\bar\partial\beta = \alpha_1 - \alpha$, and then $(\alpha_1 - \alpha,\partial\beta,0,\ldots)$ is an acyclic zigzag. In other words, we can find an acyclic zigzag $(\alpha,\alpha_2,\ldots)$ extending $\alpha$. Thus, we have $\alpha \in Z_r^{p,q}$ for every $r \geq 1$. Since the kernel of $d_r$ is given by $Z_{r + 1}^{p,q}/B_r^{p,q}$, this implies $d_r = 0$, and the spectral sequence degenerates at $E_1$.
	
	Now, we show the equivalence of \ref{item:hyper-surj} and \ref{item:hyper-inj}. To this end, we consider the short exact sequence 
	\[0 \to (F^{p + 1}\m{Tot}(I^{\bullet,\bullet}),\partial + \bar\partial) \to (F^p\m{Tot}(I^{\bullet,\bullet}),\partial + \bar\partial) \to (I^{p,\bullet}[-p],\bar\partial) \to 0\]
	of complexes, where $(I^{p,\bullet}[-p])^k \coloneq I^{p,k - p}$. The associated long exact sequence in cohomology shows that the map in \ref{item:hyper-surj} is surjective for all $p,q$ if and only if \[H^n(F^{p + 1}\m{Tot}(I^{\bullet,\bullet}),\partial + \bar\partial) \to H^n(F^p\m{Tot}(I^{\bullet,\bullet}),\partial + \bar\partial)\] 
	is injective for all $n$ and $p$. Since $I^{\bullet,\bullet}$ is bounded below in both indices, this is clearly equivalent to statement \ref{item:hyper-inj}.
\end{proof}

\subsection{Cohomological flatness}

Let us first make the definition of cohomological flatness as a filtered complex precise.

\begin{defn}\label{defn:coh-flat}\note{defn:coh-flat}
	Let $S$ be a locally Noetherian scheme, and let $\phi\colon X \to S$ be flat and proper. For a morphism $b\colon T \to S$ from a locally Noetherian scheme $T$, we write $\psi\colon Y \to T$ for the base change of $\phi$ along $b$. Furthermore, we write $c\colon Y \to X$ for the induced map.
	\begin{enum:arabic}
		\item A quasicoherent sheaf $\cF$ on $X$ is \emph{cohomologically flat} if it is flat over $S$, if $R^j\phi_*\cF$ is a flat $\cO_S$-module for every $j \geq 0$, and if, for every $b\colon T \to S$ and every $j \geq 0$ (with $T$ locally Noetherian), the canonical homomorphism 
		\[b^*R^j\phi_*\cF \to R^j\psi_*c^*\cF\]
		is an isomorphism.
		\item A \emph{first-order-differential-operator complex} $(\cA^\bullet,\partial)$ consists of a $\ZZ$-graded $\cO_X$-module $\cA^\bullet$ with $\cA^i = 0$ for $i \ll 0$, and where every $\cA^i$ is quasicoherent, and of $\phi^{-1}(\cO_S)$-linear maps $\partial\colon \cA^i \to \cA^{i + 1}$ with $\partial^2 = 0$, and which are differential operators of order $\leq 1$. Concretely, this means that we have 
		\[\partial(fga) - f\partial(ga) - g\partial(fa) + fg\partial(a) = 0\]
		for $f,g\in \cO_X$ and $a \in \cA^i$.
		
		Let $b\colon T \to S$ be a morphism from a locally Noetherian scheme $T$. Although the operators $\partial$ are not $\cO_X$-linear, they can be pulled back along $c^*$ by \cite[Prop.~C.19]{FeltenGLDT}, yielding a first-order-differential-operator complex $c^*(\cA^\bullet,\partial)$ on $Y = X \times_S T$.
		\item\label{item:fodo-cohflat} A first-order-differential-operator complex $\cA^\bullet$ is \emph{cohomologically flat} if each $\cA^i$ is flat over $S$, if $R^j\phi_*\cA^\bullet$ is a flat $\cO_S$-module, and if the canonical homomorphism 
		\[b^*R^j\phi_*\cA^\bullet \to R^j\psi_*c^*\cA^\bullet\]
		is an isomorphism for every $b\colon T \to S$ and every $j \in \ZZ$ (with $T$ locally Noetherian).
		\item A first-order-differential-operator complex $\cA^\bullet$ is \emph{cohomologically flat as a filtered complex} (\emph{cffc}) if each $\cA^i$ is quasicoherent and cohomologically flat, and if
		\[R^{i + j}\phi_*(F^i\cA^\bullet) \to R^j\phi_*\cA^i\]
		is surjective, where $F^\bullet\cA^\bullet$ is the stupid filtration characterized by $F^j\cA^i = \cA^i$ if $i \geq j$ and $F^j\cA^i = 0$ if $i < j$. \qeddef
	\end{enum:arabic}
\end{defn}
\begin{rem}
	The notion of cohomological flatness for quasicoherent sheaves (and complexes thereof) goes back at least to \cite[\S7.8]{EGAIII-2}. Although what we call \emph{cohomological flatness as a filtered complex} appears to be a well-known notion, we are not aware of any work introducing a name for it. \qedex
\end{rem}

In the rest of this section, we work in the following situation.

\begin{sitn}\label{sitn:cohomology-of-complexes}\note{sitn:cohomology-of-complexes}
	We fix a Noetherian ring $R_S$ and set $S = \Spec R_S$. We have a flat and proper morphism of schemes $\phi\colon X \to S$, as well as a first-order-differential-operator complex $\cA^\bullet_{X/S}$ in the sense of Definition~\ref{defn:coh-flat}. For every $i \in \ZZ$, we assume that $\cA^i_{X/S}$ is coherent and flat over $S$.
	
	For every affine Noetherian scheme $T = \Spec R_T$ and every morphism $b\colon T \to S$, we form the base change diagram
	\[
	\xymatrix{
		Y \ar[r]^c \ar[d]^\psi & X \ar[d]^\phi \\
		T \ar[r]^b & S. \\
	}
	\]
	We set $\cA^i_{Y/T} \coloneq c^*\cA^i_{X/S}$. Now, $\partial\colon \cA^i_{X/S} \to \cA^{i + 1}_{X/S}$ is an $\phi^{-1}(\cO_S)$-linear differential operator of order $\leq 1$, so \cite[Prop.~C.19]{FeltenGLDT} provides induced $\psi^{-1}(\cO_T)$-linear differential operators $\partial\colon \cA^i_{Y/T} \to \cA^{i + 1}_{Y/T}$ of order $\leq 1$; by \cite[Rem.~C.20]{FeltenGLDT}, we have $\partial^2 = 0$.\footnote{We do not have a base change result for arbitrary $\phi^{-1}(\cO_S)$-linear maps. The issue is localization in local sections of $\cO_Y$.}
	
	We fix an affine open cover $\cU = \{U_a\}_{a \in I}$ of $X$, indexed by the finite set $I = \{0,1,\ldots,n\}$ for some $n$. 
	On every $Y$, we consider the affine open cover $\cU \times_S T = \{U_a \times_S T\}_{a \in I}$. \qeddef
\end{sitn}

The prototypical example is the relative de Rham complex $\Omega^\bullet_{X/S}$ of a smooth and proper morphism $\phi\colon X \to S$. This example shows that we cannot assume $\partial$ to be $\cO_X$-linear. For us, the situation where $\cA^\bullet_{X/S}$ is some variant of a relative log de Rham complex is the most relevant.

Affinity of $S$ is not essential but simplifies the situation, especially when we also require $T$ to be affine, since then $\cU \times_S T$ is an affine cover of $Y$.

For every $b\colon T \to S$, we obtain a \v Cech double complex given by
\[\Gamma(W,\check\cC^j(\cU \times_S T,\cA_{Y/T}^i)) = \prod_{J \subseteq I, |J| = j + 1} \Gamma\left(W \cap \bigcap_{a \in J} (U_a \times_S T),\: \cA_{Y/T}^i\right)\]
for open subsets $W \subseteq Y$. For $\bar\partial\colon \check\cC^j(\cU \times_S T,\cA^i_{Y/T}) \to \check\cC^{j + 1}(\cU \times_S T,\cA^i_{Y/T})$, we use the usual \v Cech differential, but for the other differential, we set $\partial((\alpha_J)_J) = (-1)^j \cdot (\partial \alpha_J)_J$ in order to ensure $\partial\bar\partial + \bar\partial\partial = 0$.

We have maps $\cA^i_{Y/T} \to \check\cC^0(\cU \times_S T,\cA^i_{Y/T})$ which turn $\cA_{Y/T}^i \to \check\cC^\bullet(\cU \times_S T,\cA^i_{Y/T})$ into an $R\psi_*$-acyclic resolution. Similarly, $\cA^\bullet_{Y/T} \to \check\cC^\bullet(\cU \times_S T,\cA^\bullet_{Y/T})$ is an $R\psi_*$-acyclic resolution.

The double complex $\psi_*\check\cC^\bullet(\cU \times_S T,\cA^\bullet_{Y/T})$ defines a spectral sequence
\[E(T)\colon \enspace R^j\psi_*\cA^i_{Y/T} \Rightarrow R^{i + j}\psi_*\cA^\bullet_{Y/T}.\]
We have a canonical comparison map
\[c^*\check\cC^\bullet(\cU,\cA^\bullet_{X/S}) \to \check\cC^\bullet(\cU \times_S T,\cA^\bullet_{Y/T}),\]
which is an isomorphism. Furthermore,
\[b^*\phi_*\check\cC^\bullet(\cU,\cA^\bullet_{X/S}) \to \psi_*c^*\check\cC^\bullet(\cU,\cA^\bullet_{X/S}) \to \psi_*\check\cC^\bullet(\cU \times_S T,\cA^\bullet_{Y/T})\]
is an isomorphism. Therefore, the two associated spectral sequences are isomorphic, and we obtain a map of prespectral sequences 
\begin{equation}\label{eqn:spectral-sequence-comp}
	b^*(R^j\phi_*\cA^i_{X/S} \Rightarrow R^{i + j}\phi_*\cA^\bullet_{X/S}) \to (R^j\psi_*\cA^i_{Y/T} \Rightarrow R^{i + j}\psi_*\cA^\bullet_{Y/T})
\end{equation}
by Lemma~\ref{lem:spectral-seq-comp-map}. When $b\colon T \to S$ is flat, then Lemma~\ref{lem:spectral-seq-comp-map} guarantees that this map is an isomorphism.

\begin{lem}\label{lem:cffc-criterion}
	In Situation~\ref{sitn:cohomology-of-complexes}, the following statements are equivalent:
	\begin{enum:roman}
		\item\label{item:cffc} $\cA^\bullet_{X/S}$ is cffc in the sense of Definition~\ref{defn:coh-flat}.
		\item\label{item:cf-degen} Each $\cA^i_{X/S}$ is cohomologically flat, and the corresponding spectral sequence $E(S)$ degenerates at $E_1$.
		\item\label{item:cffc-max} For every affine Noetherian $T = \Spec R_T$ and every morphism $b\colon T \to S$, the spectral sequence $E(T)$ degenerates at $E_1$, each term $E^{i,j}_r(T)$ is a finite locally free $\cO_T$-module for every $r \geq 1$, the map \eqref{eqn:spectral-sequence-comp} is an isomorphism, each $R^k\psi_*\cA_{Y/T}^\bullet$ is a finite locally free $\cO_T$-module, and $b^*R^k\phi_*\cA_{X/S}^\bullet \to R^k\psi_*\cA_{Y/T}^\bullet$ is an isomorphism.
	\end{enum:roman}
\end{lem}
\begin{proof}
	Statements~\ref{item:cffc} and \ref{item:cf-degen} are equivalent by Lemma~\ref{lem:degen-char}. Now, assume that \ref{item:cffc-max} holds. By our assumptions in Situation~\ref{sitn:cohomology-of-complexes}, each $\cA^i_{X/S}$ is flat over $S$. For all Noetherian $T = \Spec R_T$ and morphisms $b\colon T \to S$, the sheaf $R^j\psi_*\cA^i_{Y/T}$ is finite locally free, and $b^*R^j\phi_*\cA^i_{X/S} \to R^j\psi_*\cA^i_{Y/S}$ is an isomorphism. This clearly extends to a morphism $b\colon T \to S$ from an arbitrary locally Noetherian scheme $T$, so $\cA^i_{X/S}$ is cohomologically flat in the sense of Definition~\ref{defn:coh-flat}. The spectral sequence $E(S)$ degenerates at $E_1$, so \ref{item:cf-degen} holds.
	
	Conversely, assume that \ref{item:cffc} holds. Then, we have a commutative diagram
	\[
		\xymatrix{
			b^*R^{i + j}\phi_*F^i\cA_{X/S}^\bullet \ar@{->>}[r] \ar[d] & b^*R^j\phi_*\cA_{X/S}^i \ar[d]^\cong \\
			R^{i + j}\psi_*F^i\cA_{Y/T}^\bullet \ar[r] & R^j\psi_*\cA_{Y/T}^i \\
		}
	\]
	so that the lower horizontal map is surjective. In particular, $E(T)$ degenerates at $E_1$ by Lemma~\ref{lem:degen-char}. Since $R^j\phi_*\cA_{X/S}^i$ is flat and $\phi\colon X \to S$ is proper, $R^j\phi_*\cA_{X/S}^i$ is finite locally free. Since the right vertical map in the above diagram is an isomorphism, also $E_1^{i,j}(T) = R^j\psi_*\cA_{Y/T}^i$ is finite locally free. The $E_1$-degeneration of $E(T)$ now implies that every $E_r^{i,j}$ is finite locally free, and that \eqref{eqn:spectral-sequence-comp} is an isomorphism. Now $R^{i + j}\psi_*\cA_{Y/T}^\bullet$ has a finite and exhaustive filtration with subquotients $R^j\psi_*\cA^i_{Y/T}$ so that also $R^{i + j}\psi_*\cA_{Y/T}^\bullet$ is finite locally free. Since also $R^{i + j}\phi_*\cA_{X/S}^\bullet$ has a finite exhaustive filtration with finite locally free subquotients, its pullback $b^*R^{i + j}\phi_*\cA_{X/S}^\bullet$ has such a filtration. This filtration is respected by the map $b^*R^{i + j}\phi_*\cA_{X/S}^\bullet \to R^{i + j}\psi_*\cA_{Y/T}^\bullet$ so that it being an isomorphism follows from the map \eqref{eqn:spectral-sequence-comp} being an isomorphism on the page $E_1$.
\end{proof}

\subsection{The generalized Hodge--de Rham spectral sequence}

We now specialize to the following situation.

\begin{sitn}\label{sitn:generalized-Hodge-de-Rham}
	In Situation~\ref{sitn:cohomology-of-complexes}, suppose additionally that $0 \in S \subseteq \IA^1_t = \Spec \kk[t]$, and assume that $X_0 = \phi^{-1}(0)$ is the only singular fiber. Let $S^\star = S \setminus \{0\}$ and $X^\star = \phi^{-1}(S^\star)$. Assume that $\cA^\bullet_{X^\star/S^\star} \cong \Omega^\bullet_{X^\star/S^\star}$, the relative de Rham complex of $\phi\colon X^\star \to S^\star$. Assume that $\phi\colon X \to S$ is projective.\footnote{Properness is enough, but assuming projectivity---which is given in our examples---allows us to simply cite Hartshorne's book \cite{Hartshorne1977} for cohomology and base change.} Let $b_0\colon S_0 = \{0\} \to S$ and $c_0\colon X_0 \to X$ be the closed immersions, and let $\phi_0\colon X_0 \to S_0$ be the central fiber. \qeddef
\end{sitn}

\begin{lem}
	In Situation~\ref{sitn:generalized-Hodge-de-Rham}, assume that \[b_0^*R^j\phi_*\cA_{X/S}^i \to R^j(\phi_0)_*\cA_{X_0/S_0}^i\] 
	is an isomorphism for all $i,j \in \ZZ$. Then, $\cA^\bullet_{X/S}$ is cffc.
\end{lem}
\begin{proof}
	On $S^\star$, the sheaf $R^j\phi_*\cA^i_{X^\star/S^\star} \cong R^j\phi_*\Omega^i_{X^\star/S^\star}$ is finite locally free because $\phi\colon X^\star \to S^\star$ is smooth and proper. Under the assumption of the claim, it follows from the Cohomology and Base Change Theorem \cite[Thm.~12.11]{Hartshorne1977} that $R^j\phi_*\cA_{X/S}^i$ is also finite locally free at $0 \in S$.\footnote{Here, we use our assumption that $\cA^i_{X/S}$ is flat over $S$.} Furthermore, \cite[Prop.~12.10]{Hartshorne1977} shows that $b^*R^j\phi_*\cA_{X/S}^i \to R^j\psi_*\cA_{Y/T}^i$ is an isomorphism for all affine morphisms $b\colon T \to S$. Therefore, $\cA_{X/S}^i$ is cohomologically flat.
	
	It remains to show that $E = E(S)$ degenerates at $E_1$. We prove by induction over $r \geq 1$ that $d_r^{i,j} = 0$ and that $E_{r + 1}^{i,j}$ is finite locally free. For $r = 1$, we already have seen that $E_r^{i,j}$ is finite locally free. Since $E(S^\star)$ degenerates at $E_1$, we have $d_1^{i,j}|_{S^\star} = 0$. However, since the target of $d_1^{i,j}$ is finite locally free on the integral scheme $S$, this implies $d_1^{i,j} = 0$. But then, we have $E_2^{i,j} \cong E_1^{i,j}$ so that $E_2^{i,j}$ is finite locally free. The induction step is analogous.
\end{proof}

In practice, we can easily verify this criterion as follows.

\begin{cor}\label{cor:cohom-flat-criterion}\note{cor:cohom-flat-criterion}
	In Situation~\ref{sitn:generalized-Hodge-de-Rham}, choose some $s \in S^\star$, and let $X_s = \phi^{-1}(s)$ be the fiber. Assume that
	\[h^j(X_0,\cA_{X_0/S_0}^i) = h^j(X_s,\Omega_{X_s/\kk}^i)\]
	for all $i,j \in \ZZ$. Then $\cA^\bullet_{X/S}$ is cffc.
\end{cor}
\begin{proof}
	This follows from \cite[Cor.~12.9]{Hartshorne1977} since $S$ is integral.
\end{proof}

\begin{proof}[Proof of Proposition~\ref{prop:numerical-reduction}]\label{pf:numerical-reduction}
	Recall that the setup of Proposition~\ref{prop:numerical-reduction} is an open subset $0 \in S \subseteq \IA^1_t$ and a flat projective morphism $\phi\colon X \to S$ with $d$-dimensional fibers such that $X_0$ is the only singular fiber, and such that $\phi\colon X^\partial \to S^\dagger$ is generically log smooth for a chosen $Z_X \subset X$. The reflexive log de Rham complex $\cA^\bullet_{X/S} = \cW^\bullet_{X/S}$ is a first-order-differential-operator complex. To show that $\cA^i_{X/S}$ is flat over $S$, it is sufficient to show that it is torsion-free as a $\kk[t]$-module on every affine open subset $W \subseteq X$. Since $\Omega^{i,\m{log}}_{U/S}$ is $\kk[t]$-torsion-free, also $\Gamma(W,\cA^i_{X/S}) = \Gamma(U \cap W,\Omega^{i,\m{log}}_{U/S})$ is $\kk[t]$-torsion-free. After fixing an affine open cover $\cU$ of $X$, we are in Situation~\ref{sitn:cohomology-of-complexes}.
	
	Since $\phi\colon (X^\dagger)^\star \to (S^\dagger)^\star$ is strict, we have $\cA^\bullet_{X/S}|_{U^\star} = \Omega^{\bullet,\m{log}}_{U^\star/S^\star} \cong \Omega^\bullet_{U^\star/S^\star}$ for $U^\star = X^\star \cap U$. Since $\phi\colon X^\star \to S^\star$ is smooth and $j_*\cO_{U^\star} = \cO_{X^\star}$, we have $\Omega^\bullet_{X^\star/S^\star} \cong j_*\Omega^\bullet_{U^\star/S^\star}$, and therefore $\cA^\bullet_{X/S}|_{X^\star} \cong \Omega^\bullet_{X^\star/S^\star}$. Thus, we are in Situation~\ref{sitn:generalized-Hodge-de-Rham}. Now, $\cA^\bullet_{X/S}$ is cffc by Corollary~\ref{cor:cohom-flat-criterion}.
\end{proof}

\section{The logarithmic deformation functor}\label{sect:log-defo-functor}\note{sect:log-defo-functor}

In this section, we review the construction of a logarithmic deformation functor associated with the setup of Section~\ref{sect:log-geo-setup} and show that it is unobstructed if $\cA^d_{X_0/S_0} \cong \cO_{X_0}$ and $\cA^\bullet_{X/S}$ is cffc.

In analogy with \cite[Constr.~10.27]{FeltenGLDT}, we turn $\phi\colon X^\partial \to S^\dagger$ into an \emph{enhanced generically log smooth family} $\phi\colon X^{\m{en}} \to S^\partial$ in the sense of \cite[Defn.~10.8]{FeltenGLDT}. As a part of the structure of an enhanced generically log smooth family, we endow it with the distinguished log de Rham complex
\[\cA^\bullet_{X/S} \coloneq \cW^\bullet_{X/S}\] 
and the distinguished log polyvector fields
\[\cG^{-p}_{X/S} \coloneq \cV^{-p}_{X/S} \coloneq j_*\Theta^{p,\m{log}}_{U/S}\]
for $p \geq 0$, where 
\[\Theta^{p,\m{log}}_{U/S} = \bigwedge^p_{\cO_U} \Theta^{1,\m{log}}_{U/S}\]
are the (usual) log polyvector fields.\footnote{Note our convention that $\cG^\bullet_{X/S}$ lies in negative degrees.} Together, they form a two-sided Gerstenhaber calculus in $\cO_X$-modules and $\phi^{-1}(\cO_S)$-multilinear maps. For the definition of a two-sided Gerstenhaber calculus, see \cite[Defn.~3.22]{FeltenGLDT}.

We have seen in Section~\ref{sect:log-canonical} that $\phi\colon X^\partial \to S^\dagger$ is log Gorenstein, as required in the definition of an enhanced generically log smooth family. We have established flatness of $\cA^i_{X/S} \cong \cA^{i,\m{algo}}_{\scX/\scS} \otimes_\kq \kk$ in Section~\ref{sect:log-de-Rham}, and since
\[\cG^{-p}_{X/S} \otimes \cA^d_{X/S} \cong \cA^{-p + d}_{X/S},\]
also $\cG^{-p}_{X/S}$ is flat over $S$. It is now clear that we have indeed defined an enhanced generically log smooth family. It follows from \cite[Lemma~10.26]{FeltenGLDT} and \cite[Lemma~B.4]{FeltenGLDT} that $\phi\colon X^\m{en} \to S^\dagger$ is \emph{torsionless} in the sense of \cite[Defn.~10.12]{FeltenGLDT}.

After base change along $S_0^\dagger \to S^\dagger$, we obtain an enhanced generically log smooth family $\phi_0\colon X_0^\m{en} \to S_0^\dagger$. When $c_0\colon X_0 \to X$ is the inclusion of the central fiber, the distinguished log de Rham complex is now $\cA_{X_0/S_0}^\bullet = c_0^*\cA^\bullet_{X/S}$, and the distinguished log polyvector fields are $\cG^{-p}_{X_0/S_0} = c_0^*\cG^{-p}_{X/S}$. We have canonical maps \[\varpi\colon \enspace \cA^i_{X_0/S_0} \to \cW^i_{X_0/S_0}\quad \text{and} \quad \varpi\colon \enspace \cG^{-p}_{X_0/S_0} \to \cV^{-p}_{X_0/S_0},\]
which are not necessarily isomorphisms. Here, $\cV^{-p}_{X_0/S_0} = j_*\Theta^{p,\m{log}}_{U_0/S_0}$ denotes the reflexive log polyvector fields of $\phi_0\colon X_0^\partial \to S_0^\dagger$.

We wish to construct a well-behaved deformation functor which classifies some (not necessarily all) deformations of $\phi_0\colon X_0^\m{en} \to S_0^\dagger$ over $S_A^\dagger$ for $A \in \m{Art}_{\kk\lsem t\rsem\sslash\kk}$. As a starting point, we take the base change $\phi_A\colon X_A^\m{en} \to S_A^\dagger$ of $\phi\colon X^\m{en} \to S^\dagger$ along $S_A^\dagger \to S^\dagger$. The main idea is then to classify deformations $\psi_A\colon Y_A^\m{en} \to S_A^\dagger$ which are \emph{locally} isomorphic to $\phi_A\colon X_A^\m{en} \to S_A^\dagger$.

From a technical perspective, we first construct an \emph{enhanced system of deformations} $\scD$ in the sense of \cite[Defn.~10.21(1)]{FeltenGLDT}; this can be done in analogy with \cite[Constr.~10.27]{FeltenGLDT}. In particular, we choose an open affine cover $\cV = \{V_\alpha\}_\alpha$ of $V = X_0$, and we set $V_{\alpha;A}^\m{en} \coloneq X_A^\m{en}|_{V_\alpha}$.\footnote{The primary purpose of having $\cV$ in the definition of an enhanced system of deformations is to accommodate the common situation that we do not have a global model $\phi_A\colon X_A^\m{en} \to S_A^\dagger$ for the deformations to be classified but can construct models $V_{\alpha;A}^\m{en} \to S_A^\dagger$ locally.}

An \emph{enhanced generically log smooth deformation of type $\scD$} is a deformation $\psi_A\colon Y_A^\partial \to S_A^\dagger$ of $\phi_0\colon X_0^\partial \to S_0^\dagger$ as a morphism of partial log schemes, together with isomorphisms of deformations $\chi_\alpha\colon Y_A^\partial|_{V_\alpha} \cong V_{\alpha;A}^\partial$, subject to the condition that the transition maps
\[V_{\alpha;A}^\partial|_{V_\alpha \cap V_\beta} \xrightarrow{\chi_\alpha^{-1}} Y_A^\partial|_{V_\alpha \cap V_\beta} \xrightarrow{\chi_\beta} V_{\beta;A}^\partial|_{V_\alpha \cap V_\beta} \xrightarrow{\cong} V_{\alpha;A}^\partial|_{V_\alpha \cap V_\beta}\]
arise from the action of $\cG^{-1}_{X_A/S_A}$ on the two-sided Gerstenhaber calculus via the Lie bracket and the Lie derivative. This last condition ensures that $\psi_A\colon Y_A^\partial \to S_A^\dagger$ can be endowed with the structure of an enhanced generically log smooth family $\psi_A\colon Y_A^\m{en} \to S_A^\dagger$. For more details, see \cite[Defn.~10.21(2)]{FeltenGLDT}.

We have a notion of \emph{equivalence} between enhanced generically log smooth deformations of type $\scD$ which takes into account both isomorphisms of the underlying deformations as morphisms of partial log schemes, and ``inessential'' ways of changing the isomorphisms $\chi_\alpha$. Then, equivalence classes of enhanced generically log smooth deformations of type $\scD$ form a deformation functor 
\[\m{ELD}^\scD_{X_0/S_0}\colon \enspace \m{Art}_{\kk\lsem t\rsem\sslash\kk} \to \m{Set}.\]

The whole point of introducing such a subtle deformation functor is that, unlike for instance deformations of $\phi_0\colon X_0^\partial \to S_0^\dagger$ as a morphism of partial log schemes, we have the necessary theory to study this deformation functor. For example, \cite[Cor.~13.38]{FeltenGLDT} gives a $\kk\lsem t\rsem$-linear curved Lie algebra $L^\bullet_{X_0/\kk\lsem t\rsem}$ which controls the deformation functor. For the definition of a curved Lie algebra in this setting, see \cite[Defn.~3.3]{FeltenGLDT}. For the definition of the deformation functor associated with a curved Lie algebra, see \cite[\S5.1.1]{FeltenGLDT}.

In the Calabi--Yau case, we have the following logarithmic Bogomolov--Tian--Todorov theorem.

\begin{thm}\label{thm:ELD-unobstructed}\note{thm:ELD-unobstructed}
	In the setup of Section~\ref{sect:log-geo-setup}, assume that $\delta = n + 1$ so that $\cA^d_{X_0/S_0} \cong \cO_{X_0}$. If $\cA^\bullet_{X/S}$ is cffc, then $\m{ELD}_{X_0/S_0}^\scD$ is unobstructed.
\end{thm}
\begin{proof}
	For this proof, we assume that the reader is familiar with \cite{FeltenGLDT}. Choose a second affine open cover $\cU$ of $X_0$. By \cite[Cor.~13.37]{FeltenGLDT}, we have the $\kk\lsem t\rsem$-linear curved two-sided Gerstenhaber calculus
	\[GC_{X_0/\kk\lsem t\rsem} = (PV_{X_0/\kk\lsem t\rsem}^{\bullet,\bullet},DR_{X_0/\kk\lsem t\rsem}^{\bullet,\bullet}),\]
	and for $A \in \m{Art}_{\kk\lsem t\rsem\sslash\kk}$, the \emph{universal log de Rham cohomology} is 
	\[\HH^k(X_0/A) \coloneq H^k(\m{Tot}^\bullet(DR_{X_0/\kk\lsem t\rsem}^{\bullet,\bullet} \otimes_{\kk\lsem t\rsem} A),\ \partial + \bar\partial + \cL_\ell(-)).\]
	It can be computed as 
	\[\HH^k(X_0/A) \cong R^k(\phi_A)_*\cA^\bullet_{X_A/S_A},\]
	so $\cA^\bullet_{X/S}$ being cffc implies that $\HH^k(X_0/A)$ is a finite free $A$-module, and that
	\[\HH^k(X_0/B) \otimes_B A \cong \HH^k(X_0/A)\]
	for any surjection $B \to A$ in $\m{Art}_{\kk\lsem t\rsem\sslash\kk}$.
	
	The first spectral sequence of the double complex $DR^{\bullet,\bullet}_{X_0/\kk\lsem t\rsem} \otimes_{\kk\lsem t\rsem} \kk$ can be computed as 
	\[H^j(X_0,\cA^i_{X_0/S_0}) \Rightarrow H^{i + j}(X_0,\cA^\bullet_{X_0/S_0}),\]
	which degenerates at $E_1$ because $\cA^\bullet_{X/S}$ is cffc. Thus,  $GC_{X_0/\kk\lsem t\rsem}$ is perfect in the sense of \cite[Defn.~5.18]{FeltenGLDT}. We conclude by \cite[Thm.~1.101]{FeltenGLDT} and the observation that
	\[\m{MC}(PV^{\bullet,\bullet}_{X_0/\kk\lsem t\rsem},-) \Rightarrow \m{Def}(L^\bullet_{X_0/\kk\lsem t\rsem},-)\]
	is a smooth natural transformation of deformation functors.
\end{proof}

\section{The algorithm}\label{sect:algorithm}\note{sect:algorithm}

We explain our computations. We set $\kq = \QQ$ but leave the algebraically closed extension $\kk$ unspecified. We assume that we have $n$, $d$, and $f_1,\ldots,f_r \in \QQ[t,x_0,\ldots,x_n]$ as in Datum~\ref{dat:basic-input}.

\subsection{The log-geometric setup}\label{sect:log-geo-setup-computation}\note{sect:log-geo-setup-computation}

Let $\mathtt{h} = (\QQ[t] \to R_\PP^\bullet)$ be the ring homomorphism, and let $\mathtt{ff} = (f_1,\ldots,f_r)$ be the list of functions. The method \texttt{projectiveGlsFamily(h,ff,d)} from \cite{Felten-gls-gr,FeltenZenodo}, written in Macaulay2 \cite{M2}, first tries to verify Assumption~\ref{ass:apgls}. If this fails, it throws an error. If this is correct, it tries to verify Assumption~\ref{ass:Tx-small} via Strategy~\ref{strat:Tx-sing}. If this fails, it throws an error. Otherwise, it tries to verify Assumption~\ref{ass:agls}. Again, if this fails, it throws an error. If this is correct, it attempts to find a smooth fiber by checking $t \in \{-1,1,2\}$. If none of them is smooth, it simply states that it did not find a smooth fiber. For the examples considered in this article, this is more than sufficient to verify Assumption~\ref{ass:smooth-fiber}.

When Strategy~\ref{strat:Tx-sing} fails to verify Assumption~\ref{ass:Tx-small} and we want to use Strategy~\ref{strat:Tx-unramified}, we cannot use the method \texttt{projectiveGlsFamily(h,ff,d)}. However, we can still use more specialized methods from \cite{Felten-gls-gr,FeltenZenodo} to verify Assumptions~\ref{ass:apgls}, \ref{ass:smooth-fiber}, \ref{ass:Tx-small}, and \ref{ass:agls}. This is only relevant for Example~\ref{ex:umbrella-smoothing}, where the reader can find the particulars.

In our examples, we always have $n = d + r$ so that Lemma~\ref{lem:log-canonical-ci} guarantees that Assumption~\ref{ass:Gorenstein} is satisfied.

\subsection{The Hodge numbers of the smooth fiber}

The Hodge numbers of smooth complete intersections in $\PP^n$ can be fairly easily computed with the method developed by Brückmann in \cite{Brueckmann1973,Brueckmann1976,Brueckmann1977}. This was implemented by Belmans, Best, and Chapoton \cite{BelmansZenodo}. For our examples, we read off the Hodge numbers from Belmans' website \cite{BelmansHodge}.

\subsection{The log Hodge numbers of the central fiber}\label{sect:how-to-apply}\note{sect:how-to-apply}

We explain how to compute the numbers $h^j(X_0,\cA^i_{X_0/S_0})$ which appear in Proposition~\ref{prop:numerical-reduction}. We have implemented our strategy in OSCAR \cite{OSCAR}, and the code files (one for each example) along with their results and some metadata can be found in \cite{FZ26}.

\begin{enum:arabic}
	\item To simplify our notation, we write $P^\bullet \coloneq R_\PP^\bullet$ and $R^\bullet \coloneq R^\bullet_{\widebar\scX}$. Furthermore, we write $F^\bullet$ for $R^\bullet$ considered as an $R^\bullet$-module, and $F^\bullet_\m{amb}$ for $R^\bullet$ considered as a $P^\bullet$-module via the quotient map $P^\bullet \to R^\bullet$.
	
	\item We compute a graded $R^\bullet$-module $M^{1,\bullet}$ such that $\widetilde M^{1,\bullet} \cong \Omega^1_{\widebar\scX/\widebar\scS}$ for the associated quasi-coherent sheaf on $\Proj R^\bullet = \widebar\scX$. We use the standard procedure of first computing $\Omega^1_{\PP^n \times \widebar\scS/\widebar\scS}|_{\widebar\scX}$ by means of the Euler sequence, and then dividing by the relations coming from the defining ideal of $\widebar\scX \subset \PP^n \times \widebar\scS$.
	
	\item We set 
	\[N^{1,\bullet} \coloneq (M^{1,\bullet})^\vee \coloneq  \m{GrHom}^\bullet_{R^\bullet}(M^{1,\bullet},F^\bullet).\]
	By \cite[Prop.~2.5.13]{EGAII}, we then have $\widetilde N^{1,\bullet} \cong \cH om(\Omega^1_{\widebar\scX/\widebar\scS},\cO_{\widebar\scX})$ because $R^\bullet$ is generated by $R^1$, and hence $\widetilde N^{1,\bullet} \cong \cT_{\widebar\scX/\widebar\scS}$, the classical tangent sheaf.
	
	\item We write $M^{1,\bullet}_\m{amb}$ for $M^{1,\bullet}$ considered as a $P^\bullet$-module via the quotient map $P^\bullet \to R^\bullet$. We write $N^{1,\bullet}_\m{amb}$ for $N^{1,\bullet}$ considered as a $P^\bullet$-module, and this can be computed as 
	\[N^{1,\bullet}_\m{amb} = \m{GrHom}^\bullet_{P^\bullet}(M^{1,\bullet}_\m{amb},F^\bullet_\m{amb}).\]
	In practice, computing $N^{1,\bullet}_\m{amb}$ this way tends to be faster than via computing $N^{1,\bullet}$ first.
	
	\item Let $\Theta^{1,\m{log}}_{X/S} \coloneq j_*\Theta^{1,\m{log}}_{U/S}$ be the sheaf of relative log derivations of $\phi\colon X \to S$. Then, the forgetful map $\Theta^{1,\m{log}}_{X/S} \to \cT_{X/S}$ (which maps a log derivation $(D,\Delta)$ to the derivation $D$) is an isomorphism on the strict and smooth locus\footnote{Note that this locus is smaller than the log smooth locus. The complement of the strict and smooth locus has codimension $1$ in the central fiber and codimension $\geq 2$ only in the total space while the log smooth locus has codimension $\geq 2$ already in every fiber.} of $\phi\colon X \to S$. Since the complement of this locus has codimension $\geq 2$ in the normal variety $X$, and since both $\Theta^{1,\m{log}}_{X/S}$ and $\cT_{X/S}$ are reflexive, the forgetful map is an isomorphism everywhere. Thus, we have
	\[\Theta^{1,\m{log}}_{X/S} \cong (\widetilde N^{1,\bullet} \otimes_\QQ \kk)|_X.\]
	
      \item The perfect pairing $\Omega^{1,\m{log}}_{U/S} \otimes \Omega^{d - 1,\m{log}}_{U/S} \to \Omega^{d,\m{log}}_{U/S} \cong \cO_U(-n - 1 + \delta)$ on $U \subseteq X$ induces an isomorphism $\cW^{d - 1}_{X/S} \cong \Theta^{1,\m{log}}_{X/S} \otimes \cO_X(-n - 1 + \delta)$ in the notation of Section \ref{sec:CohomologyInTheNonToroidalCase}.
	
	\item We have 
	\begin{align*}
		h^j(X_0,\cA^{d - 1}_{X_0/S_0}) &= h^j(X_0,(\Theta^{1,\m{log}}_{X/S} \otimes \cO_X(-n - 1 + \delta))|_{X_0}) \\
		&= h^j(X_0, ((\widetilde N^{1,\bullet} \otimes_\QQ \kk) \otimes \cO_{\widebar X}(-n - 1 + \delta))|_{X_0}) \\
		&= h^j(\scX_0, (\widetilde N^{1,\bullet} \otimes \cO_{\widebar\scX}(-n - 1 + \delta))|_{\scX_0})
	\end{align*}
	since cohomology commutes with flat base change along $\Spec \kk \to \Spec \QQ$.
	
	\item By \cite[Tag~01MT]{stacks}, $\widetilde N^{1,\bullet} \otimes \cO_{\widebar\scX}(-n - 1 + \delta)$ is isomorphic to the quasi-coherent sheaf associated with $N^{1,\bullet}[-n - 1 + \delta]$. By \cite[Prop.~3.5.3]{EGAII}, the restriction to $\scX_0$ is the quasi-coherent sheaf associated with $N^{1,\bullet}[-n - 1 + \delta] \otimes_{\QQ[t]} \QQ$ (for $\QQ[t] \to \QQ$ given by $t \mapsto 0$).
	
	\item Sheaf cohomology can be computed on the ambient $\PP^n_\QQ$ so that we have
	\[h^j(X_0,\cA^{d - 1}_{X_0/S_0}) = h^j(\PP^n_\QQ,N^{1,\bullet}_\m{amb}[-n - 1 + \delta] \otimes_{\QQ[t]} \QQ).\]
	Here, we suppress the passage to the associated quasi-coherent sheaf in the notation.
	
	\item To compute $h^j(\PP^n_\QQ,N^{1,\bullet}_\m{amb}[-n - 1 + \delta] \otimes_{\QQ[t]} \QQ)$, we use the Bernstein--Gelfand--Gelfand correspondence as developed by Eisenbud, Fl\o ystad, and Schreyer in \cite{EisenbudBGG2003}.
	
	\item Since $\cA^0_{X_0/S_0} \cong \cO_{X_0}$ and $\cA^d_{X_0/S_0} \cong \cO_{X_0}(-n -1 + \delta)$, we compute 
	\begin{align*}
		h^j(X_0,\cA^0_{X_0/S_0}) &= h^j(\PP^n_\QQ,F^\bullet_\m{amb} \otimes_{\QQ[t]} \QQ) \\		
		\text{and} \quad h^j(X_0,\cA^d_{X_0/S_0}) &= h^j(\PP^n_\QQ,F^\bullet_\m{amb}[-n -1 + \delta] \otimes_{\QQ[t]} \QQ) 
	\end{align*}
	similarly.
	
	\item In the case $d = 2$ of surfaces, we are finished.
	
	\item If $d \geq 3$, we set 
	\[M_\m{amb}^{[1],\bullet} \coloneq \m{GrHom}^\bullet_{P^\bullet}(N^{1,\bullet}_\m{amb},F^\bullet_\m{amb}),\]
	corresponding to $\cW^1_{X/S}$. Then, we compute $h^j(X_0,\cA^1_{X_0/S_0})$ as $h^j(\PP^n_\QQ,M^{[1],\bullet}_\m{amb} \otimes_{\QQ[t]} \QQ)$.
\end{enum:arabic}

\begin{rem}
	Two strategies to compute $N^{1,\bullet}_\m{amb}$ are possible. We can first compute $N^{1,\bullet}$ as the $R^\bullet$-dual of $M^{1,\bullet}$ and then consider it as a $P^\bullet$-module (\emph{subscheme strategy}), or we can first pass to the $P^\bullet$-module $M^{1,\bullet}_\m{amb}$ and then compute $N^{1,\bullet}_\m{amb}$ as a homomorphism module over $P^\bullet$ (\emph{ambient-space strategy}). We have compared the computation time for one implementation of each strategy in the following files \cite{FZ26}:
	\begin{quote}
		\texttt{Theta1-comparison-C21-FQA.jl} \\
		\texttt{Theta1-comparison-N22-OP5.jl} \\
		\texttt{Theta1-comparison-C22-BT4.jl}
	\end{quote}
	In Example~\ref{ex:C21-FQA} (C21-FQA), the subscheme strategy is slightly faster. In Example~\ref{ex:N22-OP5} (N22-OP5), the ambient-space strategy is several times faster than the subscheme strategy. Finally, in Example~\ref{ex:C22-BT4} (C22-BT4), the ambient-space strategy was complete after 13 seconds while we aborted the subscheme strategy after $> 45$ minutes. In view of these tests, we always use the ambient-space strategy in our computations. \qedex
\end{rem}

\begin{rem}
	Computing $M^{[1],\bullet}_\m{amb}$ directly can be expensive, and we often did not succeed in computing it. For instance, the code file \cite{FZ26}
	\begin{quote}
		\texttt{Hodge-numbers-N21-WS3-surface-A1-failed.jl}
	\end{quote}
	contains an attempt at computing $M^{[1],\bullet}_\m{amb}$ in Example~\ref{ex:N21-WS3} (N21-WS3) over $\kq = \FF_{4021}$ which we aborted after $> 1$ hour of runtime.\footnote{In the seemingly only slightly simpler Example~\ref{ex:N21-WS2} (N21-WS2) and Example~\ref{ex:N21-WS1} (N21-WS1), it took far less than one second to compute $M^{[1],\bullet}_\m{amb}$. The code files can be found in \cite{FZ26}.} Another example is the code file \cite{FZ26}
	\begin{quote}
		\texttt{Hodge-numbers-N32-ADM-FF4021-cluster.jl},
	\end{quote}
	which contains an attempt at computing $M^{[1],\bullet}_\m{amb}$ in Example~\ref{ex:N32-ADM} (N32-ADM) over $\kq = \FF_{4021}$ which we aborted after more than 3 weeks of runtime.
	
	This is the reason why we compute $\cA^1_{X/S} = \cA^{d - 1}_{X/S}$ as $\Theta^{1,\m{log}}_{X/S} \otimes \cO_X(-n - 1 + \delta)$ in the surface case. In the threefold case, we can compute $h^j(X_0,\cA^2_{X_0/S_0})$ with this method; when we are unable to compute $M^{[1],\bullet}_\m{amb}$, it is sometimes possible to compute $h^j(X_0,\cA^1_{X_0/S_0})$ with a more elaborate strategy, for instance in Example~\ref{ex:N32-ADM} (N32-ADM).\qedex
\end{rem}

\begin{rem}
  \label{rem:NonReproducibleExample}
  In Example~\ref{ex:C23-ZW3} (C23-ZW3), we have not been able to compute the log Hodge numbers with the official Release 1.8 of OSCAR \cite{OSCAR}, but a modification of the development branch (which has a different binding for computation of Gr\"obner
  bases for modules over quotient rings) does work. We intend these improvements to be included in a future official release of OSCAR and will update this preprint accordingly. All other successful computations reported in this article can be completed with the current official release. \qedex
\end{rem}

\subsection{Non-toroidal log singularities}\label{sect:non-toroidal-log-sing}\note{sect:non-toroidal-log-sing}

In the setup of Section~\ref{sect:log-geo-setup}, assume that we have established that $\cW^\bullet_{X/S}$ is cffc. Now, if $\phi\colon X^\partial \to S^\dagger$ is log toroidal, it is already known by prior results that $\cW^\bullet_{X/S}$ is cffc. If $\phi_0\colon X_0^\partial \to S_0^\dagger$ is additionally log Calabi--Yau, it is also known by prior results that the logarithmic deformation functor is unobstructed. In this section, we discuss how we establish that $\phi\colon X^\partial \to S^\dagger$ is \emph{not} log toroidal in (some of) our examples, so that these examples are indeed new.

Suppose that $\phi\colon X^\partial \to S^\dagger$ is log toroidal. By \cite[Thm.~8.2]{FFR2021}, $\cA^{d - 1}_{X_0/S_0} = \cW^{d - 1}_{X/S} \otimes \cO_{X_0}$ is reflexive in this case. Since $\cW^d_{X/S}$ is a line bundle, this is equivalent to
\[\Theta^{1,\m{log}}_{X/S} \otimes \cO_{X_0}\]
being reflexive, and we obtain a necessary criterion for log toroidality of $\phi\colon X^\partial \to S^\dagger$.

Let $S_1 = \Spec(\NN \to \kk[t]/(t^2))$, and let $\phi_1\colon X_1^\partial \to S_1^\dagger$ be the base change of $\phi\colon X^\partial \to S^\dagger$ along $S_1^\dagger \to S^\dagger$. If $\phi\colon X^\partial \to S^\dagger$ is log toroidal, then $\phi_1\colon X_1^\partial \to S_1^\dagger$ is log toroidal as well. In our examples, we use a stronger criterion that (if satisfied) shows not only that $\phi\colon X^\partial \to S^\dagger$ is not log toroidal but that $\phi_1\colon X_1^\partial \to S_1^\dagger$ is not log toroidal. We fall back to reflexivity of $\Theta^{1,\m{log}}_{X/S} \otimes \cO_{X_0}$ only if the stronger criterion is too costly to evaluate. Before explaining the criterion, we explain why we want this stronger criterion.

Let $\scD$ be the enhanced system of deformations \cite[Defn.~10.21]{FeltenGLDT} as discussed in Section~\ref{sect:log-defo-functor}. By forgetting some of the data, we obtain what we may call its underlying (plain, non-enhanced) system of deformations $\scD'$ \cite[Defn.~10.15]{FeltenGLDT}. The structure $\scD'$ has the same underlying data as a system of deformations in the sense of [\textit{loc.\,cit.}], but some of the conditions in the definition are possibly not satisfied. Concretely, in the notation of [\textit{loc.\,cit.}], $\cW^i_{V_{\alpha;A}/S_A}$ may fail to be flat over $S_A$, and its formation may fail to commute with base change.

First, suppose that $\phi\colon X^\partial \to S^\dagger$ is log toroidal. Then, $\scD'$ is a system of deformations (with all conditions satisfied), and the forgetful map 
\begin{equation}\label{eqn:ELD-LD}
	\m{ELD}^\scD_{X_0/S_0} \Rightarrow \m{LD}^{\scD'}_{X_0/S_0}
\end{equation}
is an isomorphism. If $\phi_0\colon X_0^\partial \to S_0^\dagger$ is log Calabi--Yau, the log toroidal Bogomolov--Tian--Todorov theorem guarantees unobstructedness.

Now, suppose that $\Theta^{1,\m{log}}_{X/S} \otimes \cO_{X_0}$ is indeed not reflexive at a point $p \in X_0$. Then, it could still happen that $\phi_0\colon X_0^\partial \to S_0^\dagger$ is log toroidal, and that indeed all base changes $\phi_k\colon X_k^\partial \to S_k^\dagger$ are log toroidal, for $S_k^\dagger = \Spec(\NN \to \kk[t]/(t^{k + 1}))$.\footnote{It seems unlikely that this can happen, but we did not prove that it cannot. It certainly can happen that $\phi_k$ is log toroidal for a fixed $k \geq 0$ while $\phi$ is not log toroidal.} In this case, all local models in $\scD'$ are log toroidal, $\scD'$ is once again a system of deformations, and we have the natural transformation in Display~\eqref{eqn:ELD-LD}. If $\phi_0\colon X_0^\partial \to S_0^\dagger$ is log Calabi--Yau, the log toroidal Bogomolov--Tian--Todorov theorem guarantees that $\m{LD}^{\scD'}_{X_0/S_0}$ is unobstructed. The meaning of the condition that $\Theta^{1,\m{log}}_{X/S} \otimes \cO_{X_0}$ fails to be reflexive is as follows: It implies that the Gerstenhaber calculus of the enhanced generically log smooth family $\phi_0\colon X_0^\m{en} \to S_0^\dagger$ is not isomorphic to the reflexive Gerstenhaber calculus of $\phi_0\colon X_0^\partial \to S_0^\dagger$, and that the canonical enhanced system of deformations obtained from $\scD'$ does not coincide with $\scD$. The natural transformation in Display~\eqref{eqn:ELD-LD} is possibly not an isomorphism, and the only new insight (of our computation compared to applying the log toroidal Bogomolov--Tian--Todorov theorem) is that not only $\m{LD}^{\scD'}_{X_0/S_0}$ but also $\m{ELD}^\scD_{X_0/S_0}$ is unobstructed.

To establish that we have found a genuinely new example of unobstructed deformations, we would ideally show that $\phi_0\colon X_0^\partial \to S_0^\dagger$ is not log toroidal. However, it seems to be difficult to achieve this computationally, so we show instead that at least $\phi_1\colon X_1^\partial \to S_1^\dagger$ is not log toroidal. This excludes the scenario described in the previous paragraph. We use the following strategy:

\begin{strat}\label{strat:E1}
	We would like to show that $\phi_1\colon X_1^\partial \to S_1^\dagger$ is not log toroidal. Recall that $\pi \colon |X_{0, \punct}| \to |X_0|$ is the projection from the pointed affine cone to its projective quotient. Let $W \subseteq \lvert X_0\rvert$ be an open subset, and assume that $\phi_1$ is log toroidal on $W$. Then, also $\phi_{1,\vartriangle}\colon X_{1,\vartriangle}^\partial \to S_1^\dagger$ is log toroidal on $\pi^{-1}(W) \subseteq \lvert X_{0,\punct}\rvert \subseteq \lvert X_{0,\vartriangle}\lvert$.
	
	Let
	\[\cE_1 = \cW^1_{X_{1,\vartriangle}/S_1} \otimes_{\cO_{X_{1,\vartriangle}}} \cO_{X_{0,\vartriangle}}.\]
	By \cite[Thm.~8.2]{FFR2021}, the canonical homomorphism to the bidual $\cE_1 \to \cE^{\vee\vee}_1$ is an isomorphism on $\pi^{-1}(W)$.
	
	We have $\Theta^{1,\m{log}}_{X_\vartriangle/S} \cong \cT_{X_\vartriangle/S} = \cH om(\Omega^1_{X_\vartriangle/S},\cO_{X_\vartriangle})$, and the reflexive log differential $1$-forms can be computed as 
	\[\cW^1_{X_{1,\vartriangle}/S_1} \cong \cH om(\Theta^{1,\m{log}}_{X_\vartriangle/S} \otimes_{\cO_{X_\vartriangle}} \cO_{X_{1,\vartriangle}},\cO_{X_{1,\vartriangle}}).\]
	After setting $A_1 = \QQ[t]/(t^2)$, we have $\cE_1 = \widetilde{E_1 \otimes_\QQ \kk}$ for 
	\[E_1 = \m{Hom}_{R \otimes_{\QQ[t]} A_1}\left(\m{Hom}(\Omega^1_{R/\QQ[t]},R) \otimes_{\QQ[t]} A_1, A_1\right) \otimes_{A_1} \QQ,\]
	where we write $R = R^\bullet_{\widebar\scX}$, considered as an ungraded ring.
	
	Let $K_1 = \m{ker}(E_1 \to E^{\vee\vee}_1)$ and $Q_1 = \m{coker}(E_1 \to E^{\vee\vee}_1)$. Then, we must have
	\[\pi^{-1}(W) \cap \m{Supp}(\widetilde K_1 \otimes_\QQ \kk) = \varnothing\quad \text{and} \quad \pi^{-1}(W) \cap \m{Supp}(\widetilde Q_1 \otimes_\QQ \kk) = \varnothing.\]
	This gives us a computational means to verify that a  certain locus in $X_1^\partial$ consists of non-toroidal log singularities. \qedex
\end{strat}

In Example~\ref{ex:C31-US1} and Example~\ref{ex:C32-BCU}, Strategy~\ref{strat:E1} is not sufficient to show that $\phi_1\colon X_1^\partial \to S_1^\dagger$ is not log toroidal. For these examples, we use the following alternative strategy.

\begin{strat}
	We would like to show that $\phi_1\colon X_1^\partial \to S_1^\dagger$ is not log toroidal. We compute
	\[\cE_1' \coloneq \Theta^{1,\m{log}}_{X_{1,\vartriangle}/S_1} \otimes_{\cO_{X_{1,\vartriangle}}} \cO_{X_{0,\vartriangle}} \cong \widetilde{E_1' \otimes_\QQ \kk},\]
	$K_1' = \m{ker}(E_1' \to E_1'^{\vee})$, and $Q_1' = \m{coker}(E_1' \to E_1'^{\vee\vee})$. Now, $\phi_\punct\colon X_\punct \to S$ is a Gorenstein morphism, so the normalized relative dualizing sheaf $\omega_{X_\punct/S}^\circ$ is a line bundle. Furthermore, we have $\cW^{d + 1}_{X^\circ_\punct/S} \cong \omega_{X_\punct^\circ/S}^\circ$ on the \emph{smooth} locus $X_\punct^\circ$ of $\phi_\punct\colon X_\punct \to S$, so $\cW^{d + 1}_{X_\punct/S}$ is a line bundle because the complement of $X_\punct^\circ \subseteq X_\punct$ has (absolute) codimension $\geq 2$. Now, $\Theta^{1,\m{log}}_{X_{1,\punct}/S_1}$ is locally isomorphic to $\cW^d_{X_{1,\punct}/S_1}$. In particular, we have 
	\[\pi^{-1}(W) \cap \m{Supp}(\widetilde K_1' \otimes_\QQ \kk) = \varnothing \quad \text{and} \quad \pi^{-1}(W) \cap \m{Supp}(\widetilde Q_1' \otimes_\QQ \kk) = \varnothing\]
	if $\phi_1\colon X_1^\partial \to S_1^\dagger$ is log toroidal on $W \subseteq X_1$. \qedex
\end{strat}

We provide implementations 
% We would like this to be inline, but it tends to produce overflows.
%\verb|diffNonToroidal'W1'inf(-)| and \verb|diffNonToroidal'Th1'inf(-)|
\begin{quote}
	\texttt{diffNonToroidal'W1'inf(-)} \quad and \\
	\texttt{diffNonToroidal'Th1'inf(-)} 
\end{quote}
of these strategies in the code file \texttt{new-diff-non-toroidal-methods.m2} in \cite{FZ26}. The file also contains implementations of the non-infinitesimal versions, relying on the reflexivity of $\cW^1_{X_\vartriangle/S} \otimes \cO_{X_{0,\vartriangle}}$ and $\Theta^1_{X_\vartriangle/S} \otimes \cO_{X_{0,\vartriangle}}$.

\newpage

\part{Examples}\label{part:examples}

\section{Surfaces in $\PP^3$}\label{sect:surfaces-IP3}\note{sect:surfaces-IP3}

In this section, we discuss instances of the log-geometric setup from Section~\ref{sect:log-geo-setup} with $n = 3$ and $r = 1$, i.e., families of surfaces in $\PP^3$. We set $\kq = \QQ$ but leave the algebraically closed extension $\kk$ unspecified.

\subsection{Calabi--Yau examples}

In these examples, $\cA^\bullet_{X/S}$ is cffc, and we obtain an unobstructed deformation functor.

\begin{ex}[C21-FQA]\label{ex:C21-FQA}\note{ex:C21-FQA}
	We consider $R_\PP^\bullet = \QQ[t,x,y,z,w]$,
	\[f_1 = xyzw - t(x^4 + y^4 + z^4 + w^4),\]
	and $d = 2$. The fiber $\scX_{t = 1} \subset \PP^3$ is a smooth quartic K3 surface so that Assumption~\ref{ass:smooth-fiber} is satisfied. The algorithmic tools from \cite{Felten-gls-gr,FeltenZenodo} verify that $\phi\colon \widebar\scX_\vartriangle \to \widebar\scS$ satisfies Assumptions~\ref{ass:apgls}, \ref{ass:Tx-small} via Strategy~\ref{strat:Tx-sing}, and \ref{ass:agls}. Since $d + r = 3 = n$, Lemma~\ref{lem:log-canonical-ci} yields that Assumption~\ref{ass:Gorenstein} is satisfied, and that $\cA^{d,\m{algo}}_{\widebar\scX/\widebar\scS} \cong \cO_{\widebar\scX}$.
	
	The central fiber $V = \{xyzw = 0\} \subset \PP^3$ has four irreducible components
	\begin{align*}
		V(x,y,z) &= \{w = 0\} & &\cong \Proj \kk[x,y,z], \\
		V(x,y,w) &= \{z = 0\} & &\cong \Proj \kk[x,y,w], \\
		V(x,z,w) &= \{y = 0\} & &\cong \Proj \kk[x,z,w], \\
		V(y,z,w) &= \{x = 0\} & &\cong \Proj \kk[y,z,w].
	\end{align*}
	Their pairwise intersections are 
	\begin{align*}
		D(x,y) &= V(x,y,z) \cap V(x,y,w) = \{z = w = 0\} \cong \Proj \kk[x,y], \\
		D(x,z) &= V(x,y,z) \cap V(x,z,w) = \{y = w = 0\} \cong \Proj \kk[x,z], \\
		D(x,w) &= V(x,y,w) \cap V(x,z,w) = \{y = z = 0\} \cong \Proj \kk[x,w], \\
		D(y,z) &= V(x,y,z) \cap V(y,z,w) = \{x = w = 0\} \cong \Proj \kk[y,z], \\
		D(y,w) &= V(x,y,w) \cap V(y,z,w) = \{x = z = 0\} \cong \Proj \kk[y,w], \\
		D(z,w) &= V(x,z,w) \cap V(y,z,w) = \{x = y = 0\} \cong \Proj \kk[z,w],
	\end{align*}
	and the union of these six copies of $\PP^1$ coincides with the reduced singular locus $D \subset V$. The locus $T \subset D$ coincides with the four points where three irreducible components intersect. The same description holds over $\QQ$.
	
	The scheme $V$ has normal crossings and thus toroidal crossing singularities.
	
	The kink of every $\scE \in [\scD^\circ]$ is $\ell = 1$. For $\scD(x,y)^\circ = \scD(x,y) \setminus \scT$, we have $\scZ^{(1)} \cap \scD(x,y)^\circ = \{x^4 + y^4 = 0\} \cap \scD(x,y)^\circ$, and the other $\scE \in [\scD^\circ]$ are similar. Therefore, after base change along $\Spec \kk \to \Spec \QQ$, the projective counterpart $Z^\circ$ of $Z^\circ_\vartriangle$ from Equation~\eqref{eqn:clubsuit} consists of $24$ points, four on each irreducible component of $D$.
	
	In \cite[Ex.~1.80]{FeltenGLDT}, the first author showed that each $p \in Z^\circ$ is a toroidal log singularity of type $T(1;1)$. Furthermore, the family $\phi\colon X \to S$ is semistable at the four points of $T$. Namely, $V$ is clearly a normal crossing scheme, and $X$ is smooth at any $p \in T$, so we conclude with \cite[Lemma~1.44]{FeltenGLDT}. Thus, $\phi\colon X^\partial \to S^\dagger$ is log toroidal, and we already know that $\m{LD}_{X_0/S_0}^\scD$ is unobstructed. The explicit computation yields that the (log) Hodge numbers are
	\begin{center}
		\renewcommand{\arraystretch}{1.2}
		\begin{tabular}{l|ccc}
			$t = 0$ & $h^0$ & $h^1$ & $h^2$ \\ \hline
			$\cA^2_{X/S}|_{X_0} \cong \cO_{X_0}(0)$ & 1 & 0 & 1 \\ \hline
			$\cA^1_{X/S}|_{X_0} \cong \Theta^{1,\m{log}}_{X/S}|_{X_0} \otimes \cO_{X_0}(0)$ & 0 & 20 & 0 \\ \hline
			$\cA^0_{X/S}|_{X_0} \cong \cO_{X_0}$ & 1 & 0 & 1 \\
		\end{tabular}
		\hspace{1.6cm}
		\begin{tabular}{l|ccc}
			$t = 1$ & $h^0$ & $h^1$ & $h^2$ \\ \hline
			$\Omega^2_{X_{t = 1}}$& 1 & 0 & 1 \\ \hline
			$\Omega^1_{X_{t = 1}}$ & 0 & 20 & 0 \\ \hline
			$\cO_{X_{t = 1}}$ & 1 & 0 & 1 \\
		\end{tabular}
	\end{center}
	for $t = 0$ and $t = 1$; we may interpret this as a computational check of \cite[Thm.~1.9]{FFR2021}, which asserts that the log Hodge--de Rham spectral sequence of a log toroidal family over a log point degenerates at $E_1$. \qedex
\end{ex}

\begin{figure}[t]
	\begin{mdframed}
		\includegraphics[scale=.66]{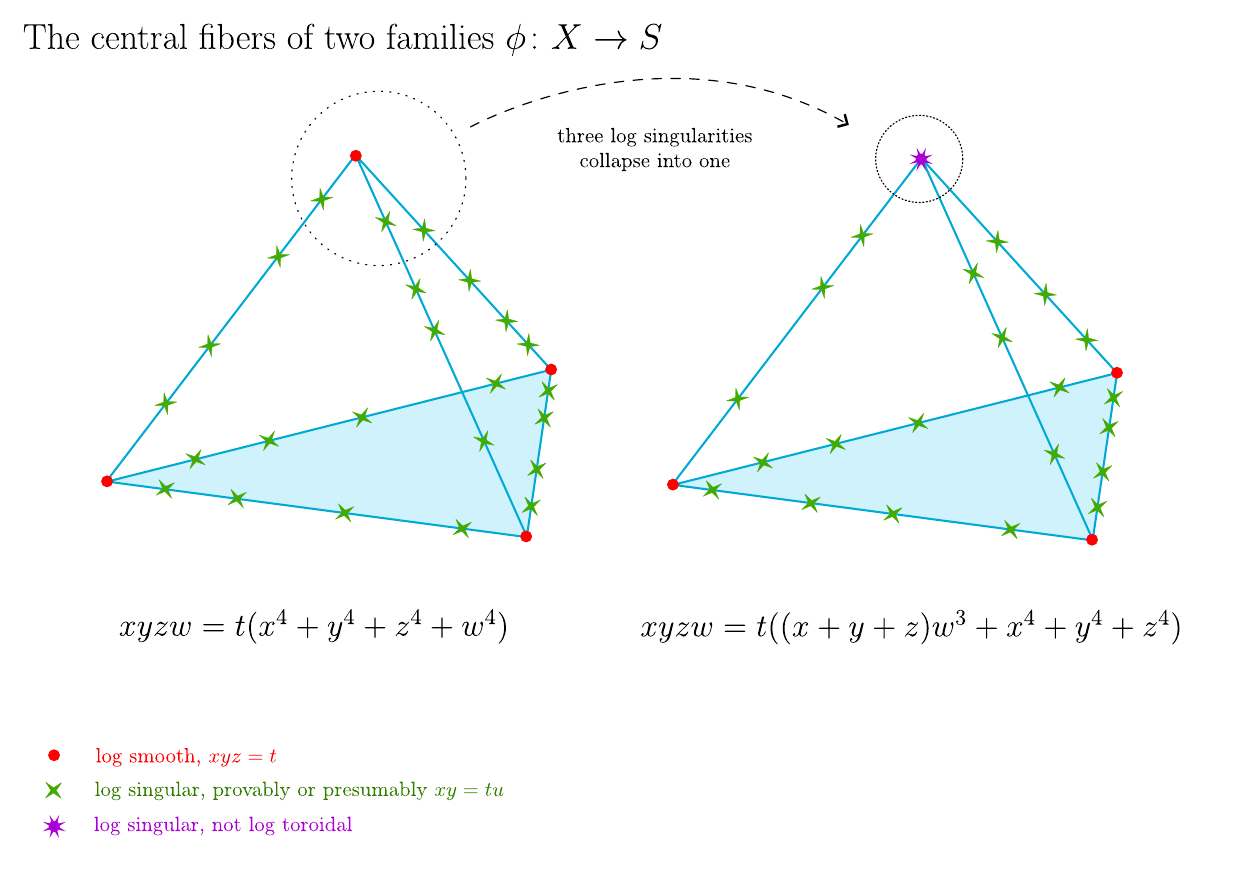}
		\caption{Example~\ref{ex:C21-FQA} and Example~\ref{ex:C21-WS3}}
	\end{mdframed}
\end{figure}

\begin{ex}[C21-WS3]\label{ex:C21-WS3}\note{ex:C21-WS3}
	We consider $R^\bullet_\PP = \QQ[t,x,y,z,w]$,
	\[f_1 = xyzw - t((x + y + z)w^3 + x^4 + y^4 + z^4),\]
	and $d = 2$. The fiber $\scX_{t = 1} = \phi^{-1}(1) \subset \PP^3$ is a smooth quartic K3 surface so that Assumption~\ref{ass:smooth-fiber} is satisfied. The algorithmic tools from \cite{Felten-gls-gr,FeltenZenodo} verify that $\phi\colon \widebar\scX_\vartriangle \to \widebar\scS$ satisfies Assumptions~\ref{ass:apgls}, \ref{ass:Tx-small} via Strategy~\ref{strat:Tx-sing}, and \ref{ass:agls}. Since $d + r = 3 = n$, Lemma~\ref{lem:log-canonical-ci} yields that Assumption~\ref{ass:Gorenstein} is satisfied, and that $\cA^{d,\m{algo}}_{\widebar\scX/\widebar\scS} \cong \cO_{\widebar\scX}$.
	
	The central fiber $V$ is the same as in the previous Example~\ref{ex:C21-FQA}.
	
	For every $\scE \in [\scD^\circ]$, the kink is $\ell(\scE) = 1$. On the one hand, we have 
	\[\scZ(x,y)^\circ = \{x^4 + y^4 = 0\}, \quad \scZ(y,z)^\circ = \{y^4 + z^4 = 0\}, \quad \scZ(x,z)^\circ = \{x^4 + z^4 = 0\}.\]
	On the other hand, we have
	\begin{align*}
		\scZ(x,w)^\circ &= \{x(x^3 + w^3) = 0\}, \\
		\scZ(y,w)^\circ &= \{y(y^3 + w^3) = 0\}, \\
		\scZ(z,w)^\circ &= \{z(z^3 + w^3) = 0\}.
	\end{align*}
	Thus, $Z^\circ$ consists of $21$ points, four on each of the three lines of $D$ in $\{w = 0\}$, and three on each of the three lines of $D$ which contain the point $p_a$ given by $\{p_a\} = \{t = x = y = z = 0\}$.
	
	Any $p \in Z^\circ$ is presumably a toroidal log singularity of type $T(1;1)$, but we did not check this. At the three points of $T \cap \{w = 0\}$, the total space $X$ is smooth, so $\phi\colon X \to S$ is semistable at these points, and hence log smooth.
	
	In the notation of Section~\ref{sect:non-toroidal-log-sing}, the radical of the annihilator of $Q_1 = \m{coker}(E_1 \to E_1^{\vee\vee})$ is $(x,y,z)$. Therefore, $\phi_1\colon X_1^\partial \to S_1^\dagger$ is not log toroidal at the fourth point $p_a$ of $T$, and we do not know unobstructedness by prior results.
	
	The (log) Hodge numbers are
	\begin{center}
		\renewcommand{\arraystretch}{1.2}
		\begin{tabular}{l|ccc}
			$t = 0$ & $h^0$ & $h^1$ & $h^2$ \\ \hline
			$\cA^2_{X/S}|_{X_0} \cong \cO_{X_0}(0)$ & 1 & 0 & 1 \\ \hline
			$\cA^1_{X/S}|_{X_0} \cong \Theta^{1,\m{log}}_{X/S}|_{X_0} \otimes \cO_{X_0}(0)$ & 0 & 20 & 0 \\ \hline
			$\cA^0_{X/S}|_{X_0} \cong \cO_{X_0}$ & 1 & 0 & 1 \\
		\end{tabular}
		\hspace{1.6cm}
		\begin{tabular}{l|ccc}
			$t = 1$ & $h^0$ & $h^1$ & $h^2$ \\ \hline
			$\Omega^2_{X_{t = 1}}$& 1 & 0 & 1 \\ \hline
			$\Omega^1_{X_{t = 1}}$ & 0 & 20 & 0 \\ \hline
			$\cO_{X_{t = 1}}$ & 1 & 0 & 1 \\
		\end{tabular}
	\end{center}
	for $t = 0$ and $t = 1$. Therefore, $\phi_0\colon X_0^\m{en} \to S_0^\dagger$ is a new example of an enhanced generically log smooth family with unobstructed deformations.\qedex
\end{ex}

\subsection{More log singularities, not Calabi--Yau}

The following families are not Calabi--Yau, but they are evidence for $\cA_{X/S}^\bullet$ being cffc under the presence of various log singularities.

\begin{ex}[N21-WS3]\label{ex:N21-WS3}\note{ex:N21-WS3}
	We consider $R^\bullet_\PP = \QQ[t,x,y,z,w]$,
	\[f_1 = xyz - t((x + y + z)w^2 + x^3 + y^3 + z^3),\]
	and $d = 2$. This is a variant of Example~\ref{ex:C21-WS3}. The fiber $\scX_{t = 1} \subset \PP^3$ is a smooth cubic surface so that Assumption~\ref{ass:smooth-fiber} is satisfied. The tools of \cite{Felten-gls-gr,FeltenZenodo} show that Assumptions~\ref{ass:apgls}, \ref{ass:Tx-small} via Strategy~\ref{strat:Tx-sing}, and \ref{ass:agls} are satisfied. Since $d + r = 3 = n$, Lemma~\ref{lem:log-canonical-ci} yields that Assumption~\ref{ass:Gorenstein} is satisfied, and that $\cA^{d,\m{algo}}_{\widebar\scX/\widebar\scS} \cong \cO_{\widebar\scX}(-1)$.
	
	The central fiber $V = \{xyz = 0\} \subset \PP^3$ is a normal crossing scheme with three irreducible components 
	\begin{align*}
		V(x,y,w) &= \{z = 0\} & &\cong \Proj \kk[x,y,w],\\
		V(x,z,w) &= \{y = 0\} & &\cong \Proj \kk[x,z,w],\\
		V(y,z,w) &= \{x = 0\} & &\cong \Proj \kk[y,z,w].
	\end{align*}
	Their pairwise intersections are 
	\begin{align*}
		D(x,w) &= V(x,y,w) \cap V(x,z,w) = \{y = z = 0\} \cong \Proj \kk[x,w],\\
		D(y,w) &= V(x,y,w) \cap V(y,z,w) = \{x = z = 0\} \cong \Proj \kk[y,w],\\
		D(z,w) &= V(x,z,w) \cap V(y,z,w) = \{x = y = 0\} \cong \Proj\kk[z,w],
	\end{align*}
	and their union is the reduced singular locus $D$. The triple locus $T$ consists of a single point $\{p_a\} = \{t = x = y = z = 0\}$.
	
	The scheme $V$ has normal crossings and thus toroidal crossing singularities.
	
	We have 
	\begin{align*}
		\scZ(x,w)^\circ &= \{x(x^2 + w^2) = 0\}, \\
		\scZ(y,w)^\circ &= \{y(y^2 + w^2) = 0\}, \\
		\scZ(z,w)^\circ &= \{z(z^2 + w^2) = 0\}
	\end{align*}
	so that $Z^\circ$ consists of $6$ points, two on each connected component of $D^\circ$. Each $p \in Z^\circ$ is presumably a toroidal log singularity of type $T(1;1)$, but we did not check this.
	
	In the notation of Section~\ref{sect:non-toroidal-log-sing}, the radical of the annihilator of $Q_1 = \m{coker}(E_1 \to E_1^{\vee\vee})$ is $(x,y,z)$. Therefore, $\phi_1\colon X_1^\partial \to S_1^\dagger$ is not log toroidal at the (unique) triple-intersection point $p_a \in T$.
	
	We obtain the (log) Hodge numbers
	\begin{center}
		\renewcommand{\arraystretch}{1.2}
		\begin{tabular}{l|ccc}
			$t = 0$ & $h^0$ & $h^1$ & $h^2$ \\ \hline
			$\cA^2_{X/S}|_{X_0} \cong \cO_{X_0}(-1)$ & 0 & 0 & 1 \\ \hline
			$\cA^1_{X/S}|_{X_0} \cong \Theta^{1,\m{log}}_{X/S}|_{X_0} \otimes \cO_{X_0}(-1)$ & 0 & 7 & 0 \\ \hline
			$\cA^0_{X/S}|_{X_0} \cong \cO_{X_0}$ & 1 & 0 & 0 \\
		\end{tabular}
		\hspace{1.6cm}
		\begin{tabular}{l|ccc}
			$t = 1$ & $h^0$ & $h^1$ & $h^2$ \\ \hline
			$\Omega^2_{X_{t = 1}}$& 0 & 0 & 1 \\ \hline
			$\Omega^1_{X_{t = 1}}$ & 0 & 7 & 0 \\ \hline
			$\cO_{X_{t = 1}}$ & 1 & 0 & 0 \\
		\end{tabular}
	\end{center}
	for $t = 0$ and $t = 1$. Thus, $\cA^\bullet_{X/S}$ is cffc. \qedex
\end{ex}

\begin{ex}[N21-WS2]\label{ex:N21-WS2}
	We consider $R^\bullet_\PP = \QQ[t,x,y,z,w]$,
	\[f_1 = xyz - t((x + y)w^2 + x^3 + y^3 + z^3),\]
	and $d = 2$. The fiber $\scX_{t = 1} \subset \PP^3$ is a smooth cubic surface so that Assumption~\ref{ass:smooth-fiber} is satisfied. The tools of \cite{Felten-gls-gr,FeltenZenodo} show that Assumptions~\ref{ass:apgls}, \ref{ass:Tx-small} via Strategy~\ref{strat:Tx-sing}, and \ref{ass:agls} are satisfied. Since $d + r = 3 = n$, Lemma~\ref{lem:log-canonical-ci} yields that Assumption~\ref{ass:Gorenstein} is satisfied, and that $\cA^{d,\m{algo}}_{\widebar\scX/\widebar\scS} \cong \cO_{\widebar\scX}(-1)$.
	
	The central fiber $V = \{xyz = 0\} \subset \PP^3$ is the same as in the previous Example~\ref{ex:N21-WS3}. For every $\scE \in [\scD^\circ]$, the kink is $\ell(\scE) = 1$.
	
	The log singular locus is 
	\[\scZ(x,w)^\circ = \{x(x^2 + w^2) = 0\}, \quad \scZ(y,w)^\circ = \{y(y^2 + w^2) = 0\}, \quad \scZ(z,w)^\circ = \{z = 0\}.\]
	Thus, $D(x,w)^\circ$ and $D(y,w)^\circ$ contain two log singular points each, but $D(z,w)^\circ$ does not contain any log singular point.
	
	In the notation of Section~\ref{sect:non-toroidal-log-sing}, the radical of the annihilator of $Q_1 = \m{coker}(E_1 \to E_1^{\vee\vee})$ is $(x,y,z)$. Therefore, $\phi_1\colon X_1^\partial \to S_1^\dagger$ is not log toroidal at the (unique) triple-intersection point $p_a \in T$.

	We obtain the following (log) Hodge numbers:
	\begin{center}
		\renewcommand{\arraystretch}{1.2}
		\begin{tabular}{l|ccc}
			$t = 0$ & $h^0$ & $h^1$ & $h^2$ \\ \hline
			$\cA^2_{X/S}|_{X_0} \cong \cO_{X_0}(-1)$ & 0 & 0 & 1 \\ \hline
			$\cA^1_{X/S}|_{X_0} \cong \Theta^{1,\m{log}}_{X/S}|_{X_0} \otimes \cO_{X_0}(-1)$ & 0 & 7 & 0 \\ \hline
			$\cA^0_{X/S}|_{X_0} \cong \cO_{X_0}$ & 1 & 0 & 0 \\
		\end{tabular}
		\hspace{1.6cm}
		\begin{tabular}{l|ccc}
			$t = 1$ & $h^0$ & $h^1$ & $h^2$ \\ \hline
			$\Omega^2_{X_{t = 1}}$& 0 & 0 & 1 \\ \hline
			$\Omega^1_{X_{t = 1}}$ & 0 & 7 & 0 \\ \hline
			$\cO_{X_{t = 1}}$ & 1 & 0 & 0 \\
		\end{tabular}
	\end{center}
	In particular, $\cA^\bullet_{X/S}$ is cffc. \qedex
\end{ex}

\begin{ex}[N21-WS1]\label{ex:N21-WS1}
	We consider $R^\bullet_\PP = \QQ[t,x,y,z,w]$,
	\[f_1 = xyz - t(xw^2 + x^3 + y^3 + z^3),\]
	and $d = 2$. The fiber $\scX_{t = 1} \subset \PP^3$ is a smooth cubic surface so that Assumption~\ref{ass:smooth-fiber} is satisfied. The tools of \cite{Felten-gls-gr,FeltenZenodo} show that Assumptions~\ref{ass:apgls}, \ref{ass:Tx-small} via Strategy~\ref{strat:Tx-sing}, and \ref{ass:agls} are satisfied. Since $d + r = 3 = n$, Lemma~\ref{lem:log-canonical-ci} yields that Assumption~\ref{ass:Gorenstein} is satisfied, and that $\cA^{d,\m{algo}}_{\widebar\scX/\widebar\scS} \cong \cO_{\widebar\scX}(-1)$.
	
	The central fiber $V = \{xyz = 0\} \subset \PP^3$ is the same as in Example~\ref{ex:N21-WS3}. For every $\scE \in [\scD^\circ]$, the kink is $\ell(\scE) = 1$.
	
	The log singular locus is 
	\[\scZ(x,w)^\circ = \{x(x^2 + w^2) = 0\}, \quad \scZ(y,w)^\circ = \{y = 0\}, \quad \scZ(z,w)^\circ = \{z = 0\}.\]
	Thus, $D(x,w)^\circ$ contains two log singular points while $D(y,w)^\circ$ and $D(z,w)^\circ$ contain no log singular point.
	
	In the notation of Section~\ref{sect:non-toroidal-log-sing}, the radical of the annihilator of $Q_1 = \m{coker}(E_1 \to E_1^{\vee\vee})$ is $(x,y,z)$. Therefore, $\phi_1\colon X_1^\partial \to S_1^\dagger$ is not log toroidal at the (unique) triple-intersection point $p_a \in T$.
	
	We obtain the following (log) Hodge numbers:
	\begin{center}
		\renewcommand{\arraystretch}{1.2}
		\begin{tabular}{l|ccc}
			$t = 0$ & $h^0$ & $h^1$ & $h^2$ \\ \hline
			$\cA^2_{X/S}|_{X_0} \cong \cO_{X_0}(-1)$ & 0 & 0 & 1 \\ \hline
			$\cA^1_{X/S}|_{X_0} \cong \Theta^{1,\m{log}}_{X/S}|_{X_0} \otimes \cO_{X_0}(-1)$ & 0 & 7 & 0 \\ \hline
			$\cA^0_{X/S}|_{X_0} \cong \cO_{X_0}$ & 1 & 0 & 0 \\
		\end{tabular}
		\hspace{1.6cm}
		\begin{tabular}{l|ccc}
			$t = 1$ & $h^0$ & $h^1$ & $h^2$ \\ \hline
			$\Omega^2_{X_{t = 1}}$& 0 & 0 & 1 \\ \hline
			$\Omega^1_{X_{t = 1}}$ & 0 & 7 & 0 \\ \hline
			$\cO_{X_{t = 1}}$ & 1 & 0 & 0 \\
		\end{tabular}
	\end{center}
	Thus, $\cA^\bullet_{X/S}$ is cffc. \qedex
\end{ex}

\begin{ex}[P21-FQI]
  \label{ex:P21-FQI} We consider $R_\PP^\bullet = \QQ[t,x,y,z,w]$,
	\[f_1 = xyzw(x + y + z + w) - t(x^5 + y^5 + z^5 + w^5),\]
	and $d = 2$. While $\scX_{t = 1} = \phi^{-1}(1) \subset \PP^3$ is not smooth, $\scX_{t = -1} = \phi^{-1}(-1)$ is a smooth quintic surface; in particular, Assumption~\ref{ass:smooth-fiber} is satisfied. The tools of \cite{Felten-gls-gr,FeltenZenodo} show that Assumptions~\ref{ass:apgls}, \ref{ass:Tx-small} via Strategy~\ref{strat:Tx-sing}, and \ref{ass:agls} are satisfied. Since $d + r = 3 = n$, Lemma~\ref{lem:log-canonical-ci} yields that Assumption~\ref{ass:Gorenstein} is satisfied, and that $\cA^{d,\m{algo}}_{\widebar\scX/\widebar\scS} \cong \cO_{\widebar\scX}(1)$.
	
	Let $s = x + y + z + w$. The central fiber $V = \{xyzws = 0\}$ has five irreducible components $V(x,y,z) = \{w = 0\}$, $V(x,y,w) = \{z = 0\}$, $V(x,z,w) = \{y = 0\}$, $V(y,z,w) = \{x = 0\}$, and $H = \{s = 0\}$, each of which is isomorphic to $\PP^2$.
	
	There are ten intersections of two irreducible components, each of which is a copy of $\PP^1$. We divide them into two groups; the first group consists of $D(x,y) = V(x,y,z) \cap V(x,y,w)$ and similarly $D(x,z)$, $D(x,w)$, $D(y,z)$, $D(y,w)$, and $D(z,w)$; the second group consists of $H(x,y,z) = H \cap V(x,y,z)$ and similarly $H(x,y,w)$, $H(x,z,w)$, $H(y,z,w)$. The union of these ten lines coincides with the reduced singular locus $D$.
	
	There are ten intersections of three irreducible components, each of which is a point. We denote the four triple-intersection points purely among the $V(x,y,z)$ by $T_0 = \{p_x, p_y, p_z, p_w\}$. We denote the other six triple-intersection points by $T_s = \{q_{xy}, q_{xz}, q_{xw}, q_{yz}, q_{yw}, q_{zw}\}$.
	
	There are no points where four irreducible components intersect. It is then easy to see that $V$ is a simple normal crossing scheme; in particular, it has toroidal crossing singularities.
	
	For every $\scE \in [\scD^\circ]$, the kink is $\ell(\scE) = 1$. For $u \not= v \in \{x,y,z,w\}$, the log singular locus on $\scD(u,v)^\circ$ is
	\[\scZ(u,v)^\circ = \{u^5 + v^5 = 0\}.\]
	After base change to $\kk$, there are four log singular points on $D(u,v)^\circ$ while the fifth lies on the intersection with $H$. For pairwise distinct $u,v,r \in \{x,y,z,w\}$, the log singular locus on $\scH(u,v,r)^\circ$ is 
	\[\scZ_H(u,v,r)^\circ = \{uvr(r^2 - uv) = 0\}.\]
	The reader may check that the closed subset cut out by this equation coincides with the one cut out by $uvr(u^2 - vr)$ and with the one cut out by $uvr(v^2 - ur)$. After base change to $\kk$, there are five distinct solutions, three of which lie in the triple locus $T$. Thus, $H(u,v,r)^\circ$ contains two log singular points.
	
	One may check that none of the four points of $T_0 = \{p_x,p_y,p_z,p_w\}$ lie in the singular locus of the total space $X$. Since $V$ is a normal crossing scheme, this implies that $\phi\colon X \to S$ is semistable at every $p \in T_0$, and hence log smooth.
	
	We obtain the (log) Hodge numbers
	\begin{center}
		\renewcommand{\arraystretch}{1.2}
		\begin{tabular}{l|ccc}
			$t = 0$ & $h^0$ & $h^1$ & $h^2$ \\ \hline
			$\cA^2_{X/S}|_{X_0} \cong \cO_{X_0}(1)$ & 4 & 0 & 1 \\ \hline
			$\cA^1_{X/S}|_{X_0} \cong \Theta^{1,\m{log}}_{X/S}|_{X_0} \otimes \cO_{X_0}(1)$ & 0 & 45 & 0 \\ \hline
			$\cA^0_{X/S}|_{X_0} \cong \cO_{X_0}$ & 1 & 0 & 4 \\
		\end{tabular}
		\hspace{1.6cm}
		\begin{tabular}{l|ccc}
			$t = -1$ & $h^0$ & $h^1$ & $h^2$ \\ \hline
			$\Omega^2_{X_{t = -1}}$& 4 & 0 & 1 \\ \hline
			$\Omega^1_{X_{t = -1}}$ & 0 & 45 & 0 \\ \hline
			$\cO_{X_{t = -1}}$ & 1 & 0 & 4 \\
		\end{tabular}
	\end{center}
	for $t = 0$ and $t = -1$, computed over $\kq = \FF_{4021}$ in the case of $t = 0$. We have not been able to compute them over $\kq = \QQ$.\footnote{The reader can find the code files for two failed attempts in \cite{FZ26}:
	\begin{quote}
		\texttt{Hodge-numbers-P21-FQI-QQ-ambient-space-strategy-cluster.jl} \\
		\texttt{Hodge-numbers-P21-FQI-QQ-subscheme-strategy-cluster.jl}
	\end{quote}
	For the ambient-space strategy, we aborted the attempt after $> 3$ weeks of runtime, and for the subscheme strategy, we aborted the attempt after 94 hours of runtime. In both cases, computing the dual of $\Omega^1_{\widebar\scX/\widebar\scS}$ is the point of failure. More information can be found in the files.}
	If the numbers are the same over $\QQ$, then $\cA^\bullet_{X/S}$ is cffc. \qedex
\end{ex}

\subsection{Smoothing of the Whitney umbrella}

Here, we see a smoothing of a projective version of the Whitney umbrella and find again that $\cA_{X/S}^\bullet$ is cffc.

\begin{ex}[N21-SWU]\label{ex:umbrella-smoothing}\note{ex:umbrella-smoothing}
	We consider $R^\bullet_\PP = \QQ[t,x,y,z,w]$,
	\[f_1 = x^2w - y^2z - t(x^3 + y^3 + z^3 + w^3),\]
	and $d = 2$. The central fiber $\scV = \{x^2w = y^2z\} \subset \PP^3$ is a projective variant of the Whitney umbrella, and the fiber $\scX_{t = 1} = \phi^{-1}(1) \subset \PP^3$ is a smooth cubic surface. In particular, Assumption~\ref{ass:smooth-fiber} is satisfied.
	
	The method \texttt{projectiveGlsFamily(h,ff,d)} throws an error. (Below, we will see why.) Therefore, we have to use lower-level methods from \cite{Felten-gls-gr,FeltenZenodo}.
	
	The method \texttt{isAlgoPreGLS(h,ff,d + 1)} from \cite{Felten-gls-gr,FeltenZenodo} confirms that Assumption~\ref{ass:apgls} is satisfied.
	
	The central fiber $\scV_\vartriangle$ is an integral Gorenstein scheme of dimension $3$, and its reduced singular locus is $\scD = \{x = y = 0\} \cong \Proj \QQ[z,w]$. When we follow Strategy~\ref{strat:Tx-sing} to verify Assumption~\ref{ass:Tx-small}, we therefore obtain $\scT_\vartriangle = \varnothing$, but $\scV_\vartriangle \setminus \scT_\vartriangle$ is not a double normal crossing scheme. As described in Section~\ref{sect:log-geo-setup-computation}, \texttt{projectiveGlsFamily(h,ff,d)} throws an error in this situation. Since Strategy~\ref{strat:Tx-sing} fails to confirm Assumption~\ref{ass:Tx-small}, we employ Strategy~\ref{strat:Tx-unramified}.
	
	Let $\nu\colon \widetilde\scV_\vartriangle \to \scV_\vartriangle$ be the normalization, let $\scV_\vartriangle^\times \subseteq \scV_\vartriangle$ be the (intrinsic) double normal crossing locus as in Section~\ref{sect:V-geometry}, and let $\scT_\vartriangle^\times = \scV_\vartriangle \setminus \scV_\vartriangle^\times$. Since $\nu\colon \widetilde\scV_\vartriangle \to \scV_\vartriangle$ is unramified over $\scV_\vartriangle^\times$, we have 
	\[\scT_\vartriangle^\times \supset \m{supp}(\nu_*\Omega^1_{\widetilde\scV_\vartriangle/\scV_\vartriangle}) = \{x = y = zw = 0\} \subset \scV_\vartriangle,\]
	the union of two lines. We make an ansatz $\scA_\vartriangle = \{x = y = zw = 0\}$ for $\scT_\vartriangle^\times$. Now, we check computationally that $\scV_\vartriangle \setminus \scA_\vartriangle$ is intrinsically double normal crossing, and we find that $\scT_\vartriangle^\times = \scA_\vartriangle$. In particular, we have $\m{dim}(\scT_\vartriangle^\times) = 1 \leq (d + 1) - 2$, and Assumption~\ref{ass:Tx-small} is satisfied.
	
	The method \texttt{Z'locus()} from \cite{Felten-gls-gr,FeltenZenodo} allows us to compute
	\[\scZ^{(1)}_\vartriangle = \{x = y = z^3 + w^3 = 0\} \subset \scD_\vartriangle.\]
	Thus, for the one irreducible component of $\scD^\circ_\vartriangle$, we have $\m{dim}(\scD^\circ_\vartriangle \cap \scZ^{(1)}_\vartriangle) = 1 \leq (d + 1) - 2$, and the kink is $\ell(\scD^\circ_\vartriangle) = 1$. Assumption~\ref{ass:agls} is satisfied, and $Z$ consists of the three points of $Z^\circ$ and the two points of $T^\times$.
	
	We have $d + r = 3 = n$, and Lemma~\ref{lem:log-canonical-ci} shows that Assumption~\ref{ass:Gorenstein} is satisfied, and that $\scA_{\widebar\scX/\widebar\scS}^{d,\m{algo}} \cong \cO_{\widebar\scX}(-1)$.
	
	The three points of $Z^\circ$ are presumably toroidal log singularities of type $T(1;1)$, but we did not check this.
	
	In the notation of Section~\ref{sect:non-toroidal-log-sing}, the radical of the annihilator of $Q_1 = \m{coker}(E_1 \to E_1^{\vee\vee})$ is $(x,y,zw)$. Therefore, $\phi_1\colon X_1^\partial \to S_1^\dagger$ is not log toroidal at the two points in $T^\times$.
	
	The (log) Hodge numbers are
	\begin{center}
		\renewcommand{\arraystretch}{1.2}
		\begin{tabular}{l|ccc}
			$t = 0$ & $h^0$ & $h^1$ & $h^2$ \\ \hline
			$\cA^2_{X/S}|_{X_0} \cong \cO_{X_0}(-1)$ & 0 & 0 & 1 \\ \hline
			$\cA^1_{X/S}|_{X_0} \cong \Theta^{1,\m{log}}_{X/S}|_{X_0} \otimes \cO_{X_0}(-1)$ & 0 & 7 & 0 \\ \hline
			$\cA^0_{X/S}|_{X_0} \cong \cO_{X_0}$ & 1 & 0 & 0 \\
		\end{tabular}
		\hspace{1.6cm}
		\begin{tabular}{l|ccc}
			$t = 1$ & $h^0$ & $h^1$ & $h^2$ \\ \hline
			$\Omega^2_{X_{t = 1}}$& 0 & 0 & 1 \\ \hline
			$\Omega^1_{X_{t = 1}}$ & 0 & 7 & 0 \\ \hline
			$\cO_{X_{t = 1}}$ & 1 & 0 & 0 \\
		\end{tabular}
	\end{center}
	for $t = 0$ and $t = 1$. In particular, $\cA_{X/S}^\bullet$ is cffc.
	
	The central fiber $V = \scV \otimes_\QQ \kk$ does \emph{not} have toroidal crossing singularities. Namely, if a scheme $V$ has toroidal crossing singularities and $D = \m{Sing}(V)$ is smooth, then $V$ is a double normal crossing scheme. On the other hand, $V$ is semi-smooth in the sense of \cite[Defn.~4.2]{KSB1988}, and in particular semi-canonical (cf.~Section~\ref{sect:experimental-2}). \qedex
\end{ex}

\subsection{Degenerations with normal singularities}\label{sect:normal-surface}

In this section, we consider examples where the central fiber $X_0$ has normal singularities. These examples are not immediately relevant to our main log-geometric interest, but it allows us to find examples where $\cA^\bullet_{X/S}$ is not cffc. In this section, we assume $\kk = \CC$ to use results from complex geometry to discuss the singularities of our examples.

\begin{ex}[N21-SE6]\label{ex:N21-SE6}
	We consider $R_\PP^\bullet = \QQ[t,x,y,z,w]$,
	\[f_1 = zx^2 - y^3 + wz^2 - t(x^3 + y^3 + z^3 + w^3),\]
	and $d = 2$. The equation for the central fiber $X_0$ was suggested by ChatGPT 5.6 as a surface with a single $E_6$-singularity and no other singularities. Indeed, one verifies easily that $\{p_w\} = \{x = y = z = 0\}$ is the only singularity of $X_0$. In the next step, the Tjurina algebra at $p_w$ is isomorphic to $\CC[x,y]/(x^3,y^2)$, so the theorem of Mather--Yau \cite[Thm.~I.2.26]{Greuel2007} shows that $p$ is an $E_6$-singularity.
	
	The methods of \cite{Felten-gls-gr} show that Assumption~\ref{ass:apgls} and Assumption~\ref{ass:smooth-fiber} are satisfied. A smooth fiber is given by $\scX_{t = 1}$. We have $\scT^\times_\vartriangle = \{x = y = z = 0\} \cong \IA^1_w$, so Assumption~\ref{ass:Tx-small} is satisfied. Next, $\scD^\circ_\vartriangle$ is empty, so Assumption~\ref{ass:agls} is satisfied. Finally, Assumption~\ref{ass:Gorenstein} is satisfied via Lemma~\ref{lem:log-canonical-ci}. We have $\cA^{d,\m{algo}}_{\widebar\scX/\widebar\scS} \cong \cO_{\widebar\scX}(-1)$.
	
	A computation in Macaulay2 shows that $\cA^1_{X_0/S_0}$ is not reflexive at $p_w$. This is used in Remark~\ref{rem:Jorder} to compare the present example with results of \cite{Jorder2014}.
	
	We obtain the (log) Hodge numbers
	\begin{center}
		\renewcommand{\arraystretch}{1.2}
		\begin{tabular}{l|ccc}
			$t = 0$ & $h^0$ & $h^1$ & $h^2$ \\ \hline
			$\cA^2_{X/S}|_{X_0} \cong \cO_{X_0}(-1)$ & 0 & 0 & 1 \\ \hline
			$\cA^1_{X/S}|_{X_0} \cong \Theta^{1,\m{log}}_{X/S}|_{X_0} \otimes \cO_{X_0}(-1)$ & 0 & 7 & 0 \\ \hline
			$\cA^0_{X/S}|_{X_0} \cong \cO_{X_0}$ & 1 & 0 & 0 \\
		\end{tabular}
		\hspace{1.6cm}
		\begin{tabular}{l|ccc}
			$t = 1$ & $h^0$ & $h^1$ & $h^2$ \\ \hline
			$\Omega^2_{X_{t = 1}}$& 0 & 0 & 1 \\ \hline
			$\Omega^1_{X_{t = 1}}$ & 0 & 7 & 0 \\ \hline
			$\cO_{X_{t = 1}}$ & 1 & 0 & 0 \\
		\end{tabular}
	\end{center}
	for $t = 0$ and $t = 1$. Thus, $\cA^\bullet_{X/S}$ is cffc.\qedex
\end{ex}

\begin{ex}[N21-SSE]\label{ex:N21-SSE}
	We consider $R_\PP^\bullet = \QQ[t,x,y,z,w]$,
	\[f_1 = x^3 + y^3 + z^3 + tw^3,\]
	and $d = 2$. This equation was suggested by ChatGPT 5.6 as a family of log canonical surfaces where $R^j\phi_*\cW^i_{X/S}$ is not locally free. Indeed, the only singular point of $X_0$ is $p_w = [0:0:0:1]$, whose affine equation is $x^3 + y^3 + z^3 = 0$, a simple elliptic singularity \cite{Saito1974} (hence log canonical but not canonical).
	
	The methods of \cite{Felten-gls-gr} show that Assumption~\ref{ass:apgls} and Assumption~\ref{ass:smooth-fiber} are satisfied. A smooth fiber is given by $\scX_{t = 1}$. We have $\scT^\times_\vartriangle = \{x = y = z = 0\} \cong \IA^1_w$, so Assumption~\ref{ass:Tx-small} is satisfied. Next, $\scD^\circ_\vartriangle$ is empty, so Assumption~\ref{ass:agls} is satisfied. Finally, Assumption~\ref{ass:Gorenstein} is satisfied via Lemma~\ref{lem:log-canonical-ci}. We have $\cA^{d,\m{algo}}_{\widebar\scX/\widebar\scS} \cong \cO_{\widebar\scX}(-1)$.
	
		We obtain the (log) Hodge numbers
	\begin{center}
		\renewcommand{\arraystretch}{1.2}
		\begin{tabular}{l|ccc}
			$t = 0$ & $h^0$ & $h^1$ & $h^2$ \\ \hline
			$\cA^2_{X/S}|_{X_0} \cong \cO_{X_0}(-1)$ & 0 & 0 & 1 \\ \hline
			$\cA^1_{X/S}|_{X_0} \cong \Theta^{1,\m{log}}_{X/S}|_{X_0} \otimes \cO_{X_0}(-1)$ & 0 & 8 & 1 \\ \hline
			$\cA^0_{X/S}|_{X_0} \cong \cO_{X_0}$ & 1 & 0 & 0 \\
		\end{tabular}
		\hspace{1.6cm}
		\begin{tabular}{l|ccc}
			$t = 1$ & $h^0$ & $h^1$ & $h^2$ \\ \hline
			$\Omega^2_{X_{t = 1}}$& 0 & 0 & 1 \\ \hline
			$\Omega^1_{X_{t = 1}}$ & 0 & 7 & 0 \\ \hline
			$\cO_{X_{t = 1}}$ & 1 & 0 & 0 \\
		\end{tabular}
	\end{center}
	for $t = 0$ and $t = 1$. Thus, $\cA^\bullet_{X/S}$ is not cffc. \qedex
\end{ex}

\begin{ex}[C21-BE6]\label{ex:C21-BE6}
	We consider $R_\PP^\bullet = \QQ[t,x,y,z,w]$,
	\[f_1 = x^2w^2 + y^3w + z^4 - t(x^4 + y^4 + z^4 + w^4),\]
	and $d = 2$. This is our naive attempt at constructing a smoothing of an $E_6$-singularity. Indeed, the central fiber has an $E_6$-singularity at $\{p_w\} = \{x = y = z = 0\}$; however, it also has a singularity at $\{p_x\} = \{y = z = w = 0\}$, and its Tjurina algebra is isomorphic to $\CC[y,z]/(y^5,z^3)$ so that it is a Brieskorn--Pham singularity of type $x^2 + y^4 + z^6$ (in particular not du Val) by the theorem of Mather--Yau \cite[Thm.~I.2.26]{Greuel2007}.
	
	The methods of \cite{Felten-gls-gr} show that Assumption~\ref{ass:apgls} and Assumption~\ref{ass:smooth-fiber} are satisfied. While $\scX_{t = 1}$ is not smooth, a smooth fiber is given by $\scX_{t = 2}$. We have $\scT^\times_\vartriangle = \{xw = y = z = 0\}$, so Assumption~\ref{ass:Tx-small} is satisfied. Next, $\scD^\circ_\vartriangle$ is empty, so Assumption~\ref{ass:agls} is satisfied. Finally, Assumption~\ref{ass:Gorenstein} is satisfied via Lemma~\ref{lem:log-canonical-ci}. We have $\cA^{d,\m{algo}}_{\widebar\scX/\widebar\scS} \cong \cO_{\widebar\scX}$.
	
	We obtain the (log) Hodge numbers
	\begin{center}
		\renewcommand{\arraystretch}{1.2}
		\begin{tabular}{l|ccc}
			$t = 0$ & $h^0$ & $h^1$ & $h^2$ \\ \hline
			$\cA^2_{X/S}|_{X_0} \cong \cO_{X_0}(0)$ & 1 & 0 & 1 \\ \hline
			$\cA^1_{X/S}|_{X_0} \cong \Theta^{1,\m{log}}_{X/S}|_{X_0} \otimes \cO_{X_0}(0)$ & 0 & 21 & 1 \\ \hline
			$\cA^0_{X/S}|_{X_0} \cong \cO_{X_0}$ & 1 & 0 & 1 \\
		\end{tabular}
		\hspace{1.6cm}
		\begin{tabular}{l|ccc}
			$t = 2$ & $h^0$ & $h^1$ & $h^2$ \\ \hline
			$\Omega^2_{X_{t = 2}}$& 1 & 0 & 1 \\ \hline
			$\Omega^1_{X_{t = 2}}$ & 0 & 20 & 0 \\ \hline
			$\cO_{X_{t = 2}}$ & 1 & 0 & 1 \\
		\end{tabular}
	\end{center}
	for $t = 0$ and $t = 2$. Thus, $\cA^\bullet_{X/S}$ is not cffc. \qedex
\end{ex}

\begin{ex}[C21-SBP]\label{ex:C21-SBP}
	We consider $R_\PP^\bullet = \QQ[t,x,y,z,w]$,
	\[f_1 = x^4 + y^4 + z^4 + tw^4,\]
	and $d = 2$. This equation was suggested by ChatGPT 5.6 after prompting for a degeneration of surfaces where $\cA^\bullet_{X/S}$ is not cffc. The central fiber has a single (Brieskorn--Pham) singularity at $\{p_w\} = \{x = y = z = 0\}$, which is not du Val.
	
	The methods of \cite{Felten-gls-gr} show that Assumption~\ref{ass:apgls} and Assumption~\ref{ass:smooth-fiber} are satisfied. A smooth fiber is given by $\scX_{t = 1}$. We have $\scT^\times_\vartriangle = \{x = y = z = 0\}$, so Assumption~\ref{ass:Tx-small} is satisfied. Next, $\scD^\circ_\vartriangle$ is empty, so Assumption~\ref{ass:agls} is satisfied. Finally, Assumption~\ref{ass:Gorenstein} is satisfied via Lemma~\ref{lem:log-canonical-ci}. We have $\cA^{d,\m{algo}}_{\widebar\scX/\widebar\scS} \cong \cO_{\widebar\scX}$.
	
	We obtain the (log) Hodge numbers
	\begin{center}
		\renewcommand{\arraystretch}{1.2}
		\begin{tabular}{l|ccc}
			$t = 0$ & $h^0$ & $h^1$ & $h^2$ \\ \hline
			$\cA^2_{X/S}|_{X_0} \cong \cO_{X_0}(0)$ & 1 & 0 & 1 \\ \hline
			$\cA^1_{X/S}|_{X_0} \cong \Theta^{1,\m{log}}_{X/S}|_{X_0} \otimes \cO_{X_0}(0)$ & 0 & 24 & 4 \\ \hline
			$\cA^0_{X/S}|_{X_0} \cong \cO_{X_0}$ & 1 & 0 & 1 \\
		\end{tabular}
		\hspace{1.6cm}
		\begin{tabular}{l|ccc}
			$t = 1$ & $h^0$ & $h^1$ & $h^2$ \\ \hline
			$\Omega^2_{X_{t = 1}}$& 1 & 0 & 1 \\ \hline
			$\Omega^1_{X_{t = 1}}$ & 0 & 20 & 0 \\ \hline
			$\cO_{X_{t = 1}}$ & 1 & 0 & 1 \\
		\end{tabular}
	\end{center}
	for $t = 0$ and $t = 1$. Thus, $\cA^\bullet_{X/S}$ is not cffc. \qedex
\end{ex}

\section{Surfaces in $\PP^4$}\label{sect:surfaces-IP4}\note{sect:surfaces-IP4}

In this section, we discuss instances of the log-geometric setup from Section~\ref{sect:log-geo-setup} with $n = 4$ and $r = 2$, i.e., families of complete-intersection surfaces in $\PP^4$. We set $\kq = \QQ$ but leave the algebraically closed extension $\kk$ unspecified.

\subsection{Calabi--Yau examples}

In these examples, $\cA^\bullet_{X/S}$ is cffc, and we obtain an unobstructed deformation functor.

\begin{figure}[t]
	\begin{mdframed}
		\includegraphics[scale=.64]{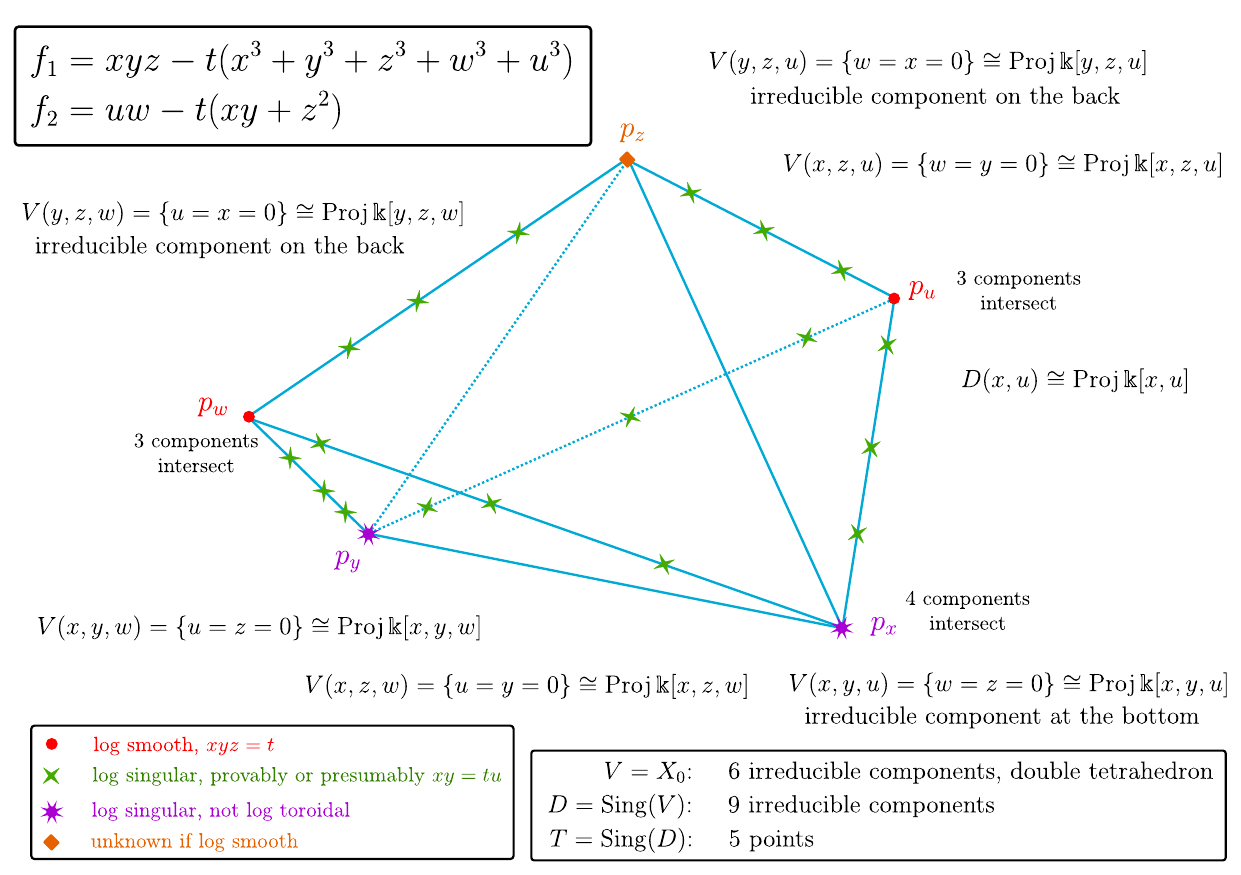}
		\caption{Example~\ref{ex:C22-BT5}}\label{fig:6-1}
	\end{mdframed}
\end{figure}

\begin{ex}[C22-BT5]\label{ex:C22-BT5}\note{ex:C22-BT5}
	We consider $R_\PP^\bullet = \QQ[t,x,y,z,w,u]$,
	\begin{align*}
		f_1 &= xyz - t(x^3 + y^3 + z^3 + w^3 + u^3), \\
		f_2 &= uw - t(xy + z^2), 
	\end{align*}
	and $d = 2$. The fiber $\scX_{t = 1} = \phi^{-1}(1) \subset \PP^4$ is a smooth sextic K3 surface so that Assumption~\ref{ass:smooth-fiber} is satisfied. The tools of \cite{Felten-gls-gr,FeltenZenodo} show that Assumptions~\ref{ass:apgls}, \ref{ass:Tx-small} via Strategy~\ref{strat:Tx-sing}, and \ref{ass:agls} are satisfied. Since $d + r = 4 = n$, Lemma~\ref{lem:log-canonical-ci} yields that Assumption~\ref{ass:Gorenstein} is satisfied, and that $\cA_{\widebar\scX/\widebar\scS}^{d,\m{algo}} \cong \cO_{\widebar\scX}$.
	
	The central fiber $V = \{xyz = uw = 0\} \subset \PP^4$ has six irreducible components, each of which is isomorphic to $\PP^2$. For their equations and intersection pattern, see Figure~\ref{fig:6-1}. For instance, one irreducible component is $V(x,y,w) = \{u = z = 0\} \cong \Proj \kk[x,y,w]$. The reduced singular locus $D = \m{Sing}(V)$ consists of nine copies of $\PP^1$, each of them the intersection of two irreducible components. For instance, one of them is $D(x,y) = V(x,y,w) \cap V(x,y,u) \cong \Proj \kk[x,y]$. Finally, the reduced singular locus $T = \m{Sing}(D)$ consists of five points $T = \{p_x,p_y,p_z,p_w,p_u\}$. At $p_w$ and $p_u$, three irreducible components of $V$ intersect, and $V$ is a normal crossing scheme at these points. In $p_x$, $p_y$, and $p_z$, four irreducible components of $V$ intersect, and the local form is $\Spec \kk[x,y]/(xy) \times \Spec \kk[u,w]/(uw)$. In particular, $V$ has toroidal crossing singularities.
	
	For every $\scE \in [\scD^\circ]$, the kink is $\ell(\scE) = 1$.
	
	The log singular locus is 
	\begin{gather*}
		\scZ(x,y)^\circ = \{xy = 0\}, \ \scZ(y,z)^\circ = \{yz = 0\}, \ \scZ(x,z)^\circ = \{xz = 0\}, \\
		\scZ(x,w)^\circ = \{w(x^3 + w^3) = 0\}, \ \scZ(y,w)^\circ = \{w(y^3 + w^3) = 0\}, \ \scZ(z,w)^\circ = \{w(z^3 + w^3) = 0\}, \\
		\scZ(x,u)^\circ = \{u(x^3 + u^3) = 0\}, \ \scZ(y,u)^\circ = \{u(y^3 + u^3) = 0\}, \ \scZ(z,u)^\circ = \{u(z^3 + u^3) = 0\}.
	\end{gather*}
	Thus, we have $\scZ(x,y)^\circ = \scZ(y,z)^\circ = \scZ(x,z)^\circ = \varnothing$, while the remaining six components of $D^\circ$ each carry three log singular points, presumably toroidal of type $T(1;1)$.
	
	The total space $X$ is smooth at the two points $p_w$ and $p_u$ so that $\phi\colon X^\dagger \to S^\dagger$ is log smooth of the form $xyz = t$ at these points.
	
	In the notation of Section~\ref{sect:non-toroidal-log-sing}, the radical of the annihilator of $Q_1 = \m{coker}(E_1 \to E_1^{\vee\vee})$ is $(xy,z,w,u)$. Therefore, $\phi_1\colon X_1^\partial \to S_1^\dagger$ is not log toroidal at $p_x$ and $p_y$. We do not know if $\phi\colon X^\dagger \to S^\dagger$ (or $\phi_0\colon X_0^\dagger \to S_0^\dagger$) is log smooth or log singular at the point $p_z$.
	
	The (log) Hodge numbers are
	\begin{center}
		\renewcommand{\arraystretch}{1.2}
		\begin{tabular}{l|ccc}
			$t = 0$ & $h^0$ & $h^1$ & $h^2$ \\ \hline
			$\cA^2_{X/S}|_{X_0} \cong \cO_{X_0}(0)$ & 1 & 0 & 1 \\ \hline
			$\cA^1_{X/S}|_{X_0} \cong \Theta^{1,\m{log}}_{X/S}|_{X_0} \otimes \cO_{X_0}(0)$ & 0 & 20 & 0 \\ \hline
			$\cA^0_{X/S}|_{X_0} \cong \cO_{X_0}$ & 1 & 0 & 1 \\
		\end{tabular}
		\hspace{1.6cm}
		\begin{tabular}{l|ccc}
			$t = 1$ & $h^0$ & $h^1$ & $h^2$ \\ \hline
			$\Omega^2_{X_{t = 1}}$& 1 & 0 & 1 \\ \hline
			$\Omega^1_{X_{t = 1}}$ & 0 & 20 & 0 \\ \hline
			$\cO_{X_{t = 1}}$ & 1 & 0 & 1 \\
		\end{tabular}
	\end{center}
	for $t = 0$ and $t = 1$. Therefore, $\phi_0\colon X_0^\m{en} \to S_0^\dagger$ is a new example of an unobstructed log Calabi--Yau surface. \qedex
\end{ex}

\begin{figure}[t]
	\begin{mdframed}
		\includegraphics[scale=.64]{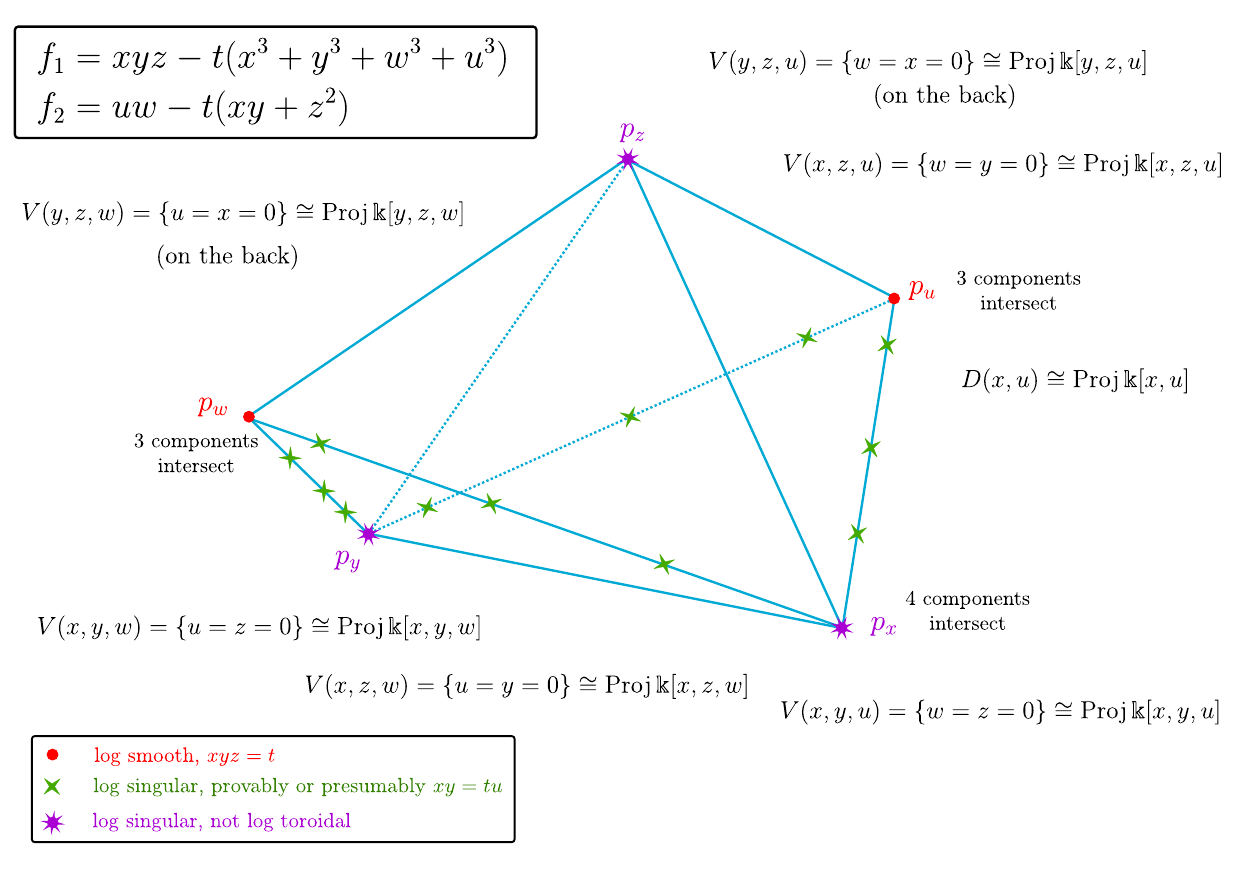}
		\caption{Example~\ref{ex:C22-BT4}}\label{fig:6-3}
	\end{mdframed}
\end{figure}

\begin{ex}[C22-BT4]\label{ex:C22-BT4}\note{ex:C22-BT4}
	We consider $R_\PP^\bullet = \QQ[t,x,y,z,w,u]$,
	\begin{align*}
		f_1 &= xyz - t(x^3 + y^3 + w^3 + u^3), \\
		f_2 &= uw - t(xy + z^2),
	\end{align*}
	and $d = 2$. The fiber $\scX_{t = 1} = \phi^{-1}(1) \subset \PP^4$ is a smooth sextic K3 surface so that Assumption~\ref{ass:smooth-fiber} is satisfied. The tools of \cite{Felten-gls-gr,FeltenZenodo} show that Assumptions~\ref{ass:apgls}, \ref{ass:Tx-small} via Strategy~\ref{strat:Tx-sing}, and \ref{ass:agls} are satisfied. Since $d + r = 4 = n$, Lemma~\ref{lem:log-canonical-ci} yields that Assumption~\ref{ass:Gorenstein} is satisfied, and that $\cA_{\widebar\scX/\widebar\scS}^{d,\m{algo}} \cong \cO_{\widebar\scX}$.
	
	The central fiber $V = X_0$ is the same as in Example~\ref{ex:C22-BT5}. For every $\scE \in [\scD^\circ]$, the kink is $\ell(\scE) = 1$.
	
	As in Example~\ref{ex:C22-BT5}, we have $\scZ(x,y)^\circ = \varnothing$, $\scZ(y,z)^\circ = \varnothing$, $\scZ(x,z)^\circ = \varnothing$. Similarly, we have 
	\begin{gather*}
		\scZ(x,w)^\circ = \{w(x^3 + w^3) = 0\}, \ \scZ(y,w)^\circ = \{w(y^3 + w^3) = 0\}, \\
		\scZ(x,u)^\circ = \{u(x^3 + u^3) = 0\}, \ \scZ(y,u)^\circ = \{u(y^3 + u^3) = 0\}
	\end{gather*}
	as in Example~\ref{ex:C22-BT5}. However, for the remaining two components of $\scD^\circ$, we have $\scZ(z,w)^\circ = \{zw = 0,\ w^4 = 0\} = \varnothing$ and $\scZ(z,u)^\circ = \{zu = 0, \ u^4 = 0\} = \varnothing$. We believe but have not proven that the $12$ log singularities in $D^\circ$ are toroidal of type $T(1;1)$.
	
	As in Example~\ref{ex:C22-BT5}, the two points $p_w$ and $p_u$ are not contained in $\m{Sing}(X)$ so that $\phi\colon X^\dagger \to S^\dagger$ is log smooth of the form $xyz = t$ at these points.
	
	In the notation of Section~\ref{sect:non-toroidal-log-sing}, the radical of the annihilator of $Q_1 = \m{coker}(E_1 \to E_1^{\vee\vee})$ is $(xy,yz,xz,w,u)$. Therefore, $\phi_1\colon X_1^\partial \to S_1^\dagger$ is not log toroidal at $p_x$, $p_y$, $p_z$. Compare this to Example~\ref{ex:C22-BT5}, where we do not know if $p_z$ is a log singularity.
	
	The (log) Hodge numbers are
	
	\begin{center}
		\renewcommand{\arraystretch}{1.2}
		\begin{tabular}{l|ccc}
			$t = 0$ & $h^0$ & $h^1$ & $h^2$ \\ \hline
			$\cA^2_{X/S}|_{X_0} \cong \cO_{X_0}(0)$ & 1 & 0 & 1 \\ \hline
			$\cA^1_{X/S}|_{X_0} \cong \Theta^{1,\m{log}}_{X/S}|_{X_0} \otimes \cO_{X_0}(0)$ & 0 & 20 & 0 \\ \hline
			$\cA^0_{X/S}|_{X_0} \cong \cO_{X_0}$ & 1 & 0 & 1 \\
		\end{tabular}
		\hspace{1.6cm}
		\begin{tabular}{l|ccc}
			$t = 1$ & $h^0$ & $h^1$ & $h^2$ \\ \hline
			$\Omega^2_{X_{t = 1}}$& 1 & 0 & 1 \\ \hline
			$\Omega^1_{X_{t = 1}}$ & 0 & 20 & 0 \\ \hline
			$\cO_{X_{t = 1}}$ & 1 & 0 & 1 \\
		\end{tabular}
	\end{center}
	for $t = 0$ and $t = 1$. Thus, $\phi_0\colon X_0^\m{en} \to S_0^\dagger$ is a new example of an enhanced generically log smooth family with unobstructed deformations. \qedex
\end{ex}

\begin{figure}[t]
	\begin{mdframed}
		\includegraphics[scale=.64]{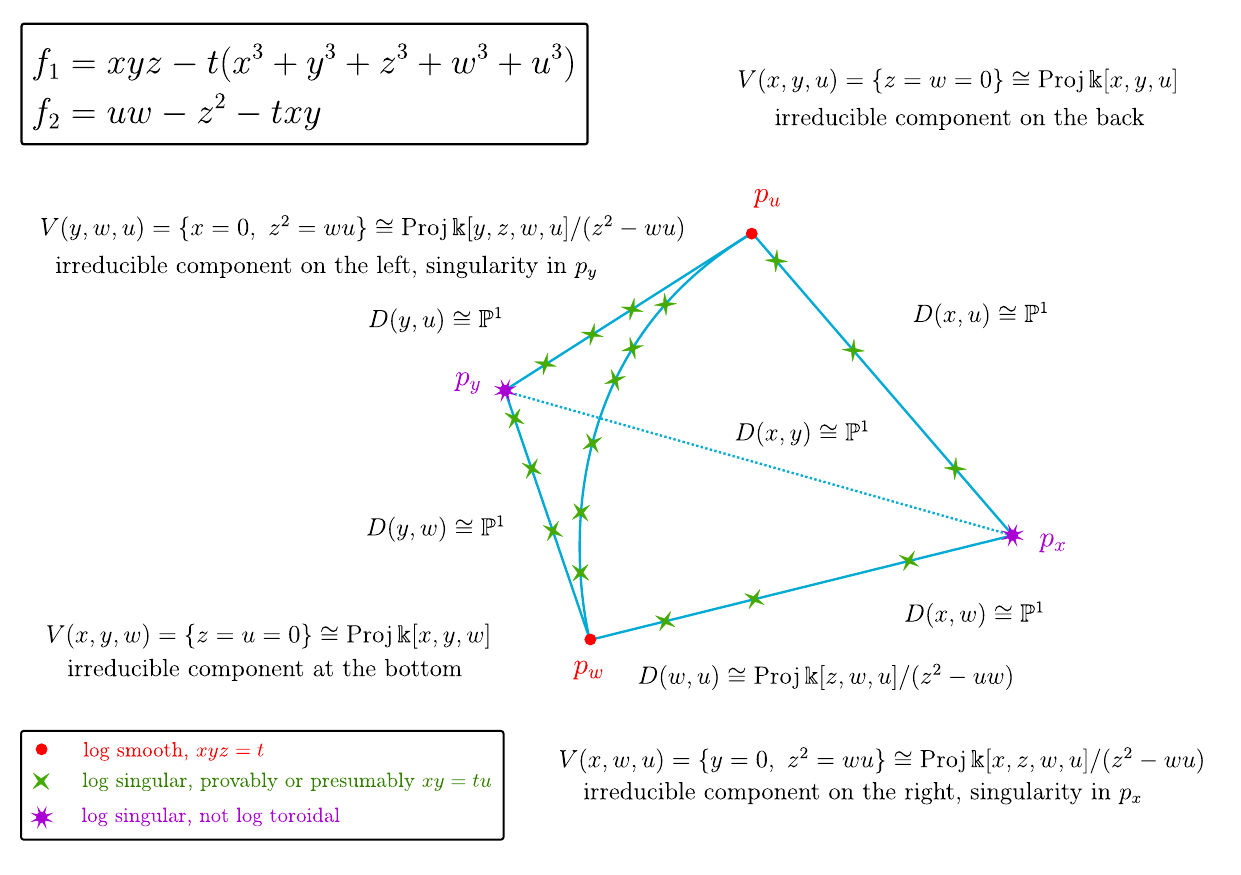}
		\caption{Example~\ref{ex:C22-DCT}}\label{fig:6-2}
	\end{mdframed}
\end{figure}

\begin{ex}[C22-DCT]\label{ex:C22-DCT}\note{ex:C22-DCT}
	We consider $R_\PP^\bullet = \QQ[t,x,y,z,w,u]$,
	\begin{align*}
		f_1 &= xyz - t(x^3 + y^3 + z^3 + w^3 + u^3), \\
		f_2 &= uw - z^2 - txy,
	\end{align*}
	and $d = 2$. The fiber $\scX_{t = 1} = \phi^{-1}(1) \subset \PP^4$ is a smooth sextic K3 surface so that Assumption~\ref{ass:smooth-fiber} is satisfied. The tools of \cite{Felten-gls-gr,FeltenZenodo} show that Assumptions~\ref{ass:apgls}, \ref{ass:Tx-small} via Strategy~\ref{strat:Tx-sing}, and \ref{ass:agls} are satisfied. Since $d + r = 4 = n$, Lemma~\ref{lem:log-canonical-ci} yields that Assumption~\ref{ass:Gorenstein} is satisfied, and that $\cA_{\widebar\scX/\widebar\scS}^{d,\m{algo}} \cong \cO_{\widebar\scX}$.
	
	The central fiber $V = X_0$ has four irreducible components as shown in Figure~\ref{fig:6-2}. Two of them are copies of $\PP^2$ while the other two are projective surfaces of degree $2$. Any two irreducible components intersect in a smooth curve. Five of them are degree-$1$ copies of $\PP^1$ while one of them is a smooth projective curve of degree $2$, isomorphic to $\PP^1$. The reduced singular locus $D = \m{Sing}(V)$ is the union of these six curves. The reduced singular locus $T = \m{Sing}(D)$ coincides with the four points $p_x$, $p_y$, $p_w$, $p_u$ where three irreducible components of $V$ intersect. The component $V(x,w,u)$ has a normal surface singularity at $p_x$, and the component $V(y,w,u)$ has a normal surface singularity at $p_y$.
	
	For every $\scE \in [\scD^\circ]$, the kink is $\ell(\scE) = 1$.
	
	The log singular locus is $\scZ(x,y)^\circ = \{xy = 0\} = \varnothing$,
	\begin{gather*}
		\scZ(x,w)^\circ = \{w(x^3 + w^3) = 0\}, \quad \scZ(y,w)^\circ = \{w(y^3 + w^3) = 0\}, \\
		\scZ(x,u)^\circ = \{u(x^3 + u^3) = 0\}, \quad \scZ(y,u)^\circ = \{u(y^3 + u^3) = 0\},
	\end{gather*}
	and $\scZ(w,u)^\circ = \{w^3 + zwu + u^3 = 0\}$. Thus, after base change along $\Spec \kk \to \Spec \QQ$, the component $D(x,y)^\circ$ has no log singular point, and each of $D(x,w)^\circ$, $D(y,w)^\circ$, $D(x,u)^\circ$, $D(y,u)^\circ$ carries three log singular points. Since $D(w,u) \cong \{z^2 = uw\} \subset \PP^2$ has itself degree $2$, and since $D(w,u) \cap \{w^3 + zwu + u^3 = 0\}$ is smooth (except in the apex of the cone), the equation $w^3 + zwu + u^3 = 0$ cuts out six points; they all lie in $D(w,u)^\circ$. We believe that these $18$ log singularities are toroidal of type $T(1;1)$, but we did not prove this.
	
	A direct computation shows that $V \cap \{w \not= 0\}$ and $V \cap \{u \not= 0\}$ are normal crossing schemes. Moreover, $p_w$ and $p_u$ are not contained in $\m{Sing}(X)$. Therefore, the family $\phi\colon X \to S$ is \'etale locally of the form $xyz = t$ around $p_w$ and $p_u$, and $\phi\colon X^\dagger \to S^\dagger$ is log smooth at $p_w$ and at $p_u$.
	
	In the notation of Section~\ref{sect:non-toroidal-log-sing}, the radical of the annihilator of $Q_1 = \m{coker}(E_1 \to E_1^{\vee\vee})$ is $(xy,z,w,u)$. Therefore, $\phi_1\colon X_1^\partial \to S_1^\dagger$ is not log toroidal at $p_x$ and $p_y$.
	
	At $p_x$ and $p_y$, the scheme $V$ is locally of the form $\{xy = w^2, \: zw = 0\} \subset \IA^4$. When using \cite[Ex.~1.66]{FeltenGLDT}, we see that $V$ has toroidal crossing singularities.
	
	The (log) Hodge numbers are
		\begin{center}
			\renewcommand{\arraystretch}{1.2}
			\begin{tabular}{l|ccc}
				$t = 0$ & $h^0$ & $h^1$ & $h^2$ \\ \hline
				$\cA^2_{X/S}|_{X_0} \cong \cO_{X_0}(0)$ & 1 & 0 & 1 \\ \hline
				$\cA^1_{X/S}|_{X_0} \cong \Theta^{1,\m{log}}_{X/S}|_{X_0} \otimes \cO_{X_0}(0)$ & 0 & 20 & 0 \\ \hline
				$\cA^0_{X/S}|_{X_0} \cong \cO_{X_0}$ & 1 & 0 & 1 \\
			\end{tabular}
			\hspace{1.6cm}
			\begin{tabular}{l|ccc}
				$t = 1$ & $h^0$ & $h^1$ & $h^2$ \\ \hline
				$\Omega^2_{X_{t = 1}}$& 1 & 0 & 1 \\ \hline
				$\Omega^1_{X_{t = 1}}$ & 0 & 20 & 0 \\ \hline
				$\cO_{X_{t = 1}}$ & 1 & 0 & 1 \\
			\end{tabular}
		\end{center}
	for $t = 0$ and $t = 1$, computed over $\kq = \FF_{4021}$ in the case of $t = 0$. We did not succeed in doing the computation over $\kq = \QQ$.  If the numbers are the same over $\kq = \QQ$, then $\cA^\bullet_{X/S}$ is cffc, and $\phi_0\colon X_0^\m{en} \to S_0^\dagger$ is a new example of an enhanced generically log smooth family with unobstructed deformations. \qedex
\end{ex}

\subsection{More log singularities, not Calabi--Yau}

The following example is special because it exhibits non-toroidal log singularities inside the double normal crossing locus of the central fiber.

\begin{ex}[N22-OP5]\label{ex:N22-OP5}
	We consider $R_\PP^\bullet = \QQ[t,x,y,z,w,u]$,
	\begin{align*}
		f_1 &= xy - t(x^2 + y^2 + z^2 + w^2 + u^2) \\
		f_2 &= uw - txy - tz^2
	\end{align*}
	and $d = 2$. The fiber $\scX_{t = -1} = \phi^{-1}(-1) \subset \PP^4$ is smooth so that Assumption~\ref{ass:smooth-fiber} is satisfied. The tools of \cite{Felten-gls-gr,FeltenZenodo} show that Assumptions~\ref{ass:apgls}, \ref{ass:Tx-small} via Strategy~\ref{strat:Tx-sing}, and \ref{ass:agls} are satisfied. Since $d + r = 4 = n$, Lemma~\ref{lem:log-canonical-ci} yields that Assumption~\ref{ass:Gorenstein} is satisfied, and that $\cA_{\widebar\scX/\widebar\scS}^{d,\m{algo}} \cong \cO_{\widebar\scX}(-1)$.
	
	The central fiber $V$ has four irreducible components, each a copy of $\PP^2$ and arranged like the four upper surfaces of a pyramid with a square basis. The apex is given by $\{p_z\} = \{x = y = w = u = 0\}$. Since $p_z$ has a neighborhood $\Spec \kk[x,y]/(xy) \times \Spec \kk[u,w]/(uw)$, the scheme $V$ has toroidal crossing singularities.
	
	In the notation of Section~\ref{sect:non-toroidal-log-sing}, the radical of the annihilator of $Q_1 = \m{coker}(E_1 \to E_1^{\vee\vee})$ is $(z,w,u)$. Therefore, $\phi_1\colon X_1^\partial \to S_1^\dagger$ is not log toroidal at $\{p_x\} = \{y = z = w = u = 0\}$ and $\{p_y\} = \{x = z = w = u = 0\}$. This situation differs from our above examples, which (presumably) do not have non-toroidal log singularities within the double normal crossing locus. We leave it to the interested reader to study the other log singularities in this family.
	
	The (log) Hodge numbers are
	\begin{center}
		\renewcommand{\arraystretch}{1.2}
		\begin{tabular}{l|ccc}
			$t = 0$ & $h^0$ & $h^1$ & $h^2$ \\ \hline
			$\cA^2_{X/S}|_{X_0} \cong \cO_{X_0}(-1)$ & 0 & 0 & 1 \\ \hline
			$\cA^1_{X/S}|_{X_0} \cong \Theta^{1,\m{log}}_{X/S}|_{X_0} \otimes \cO_{X_0}(-1)$ & 0 & 6 & 0 \\ \hline
			$\cA^0_{X/S}|_{X_0} \cong \cO_{X_0}$ & 1 & 0 & 0 \\
		\end{tabular}
		\hspace{1.6cm}
		\begin{tabular}{l|ccc}
			$t = -1$ & $h^0$ & $h^1$ & $h^2$ \\ \hline
			$\Omega^2_{X_{t = -1}}$& 0 & 0 & 1 \\ \hline
			$\Omega^1_{X_{t = -1}}$ & 0 & 6 & 0 \\ \hline
			$\cO_{X_{t = -1}}$ & 1 & 0 & 0 \\
		\end{tabular}
	\end{center}
	for $t = 0$ and $t = -1$. Thus, $\cA^\bullet_{X/S}$ is cffc. \qedex
\end{ex}

\section{Surfaces in $\PP^5$}\label{sect:surfaces-IP5}\note{sect:surfaces-IP5}

In this section, we discuss instances of the log-geometric setup from Section~\ref{sect:log-geo-setup} with $n = 5$ and $r = 3$, i.e., families of complete-intersection surfaces in $\PP^5$. We set $\kq = \QQ$ but leave the algebraically closed extension $\kk$ unspecified.

\begin{figure}[t]
	\begin{mdframed}
		\includegraphics[scale=.64]{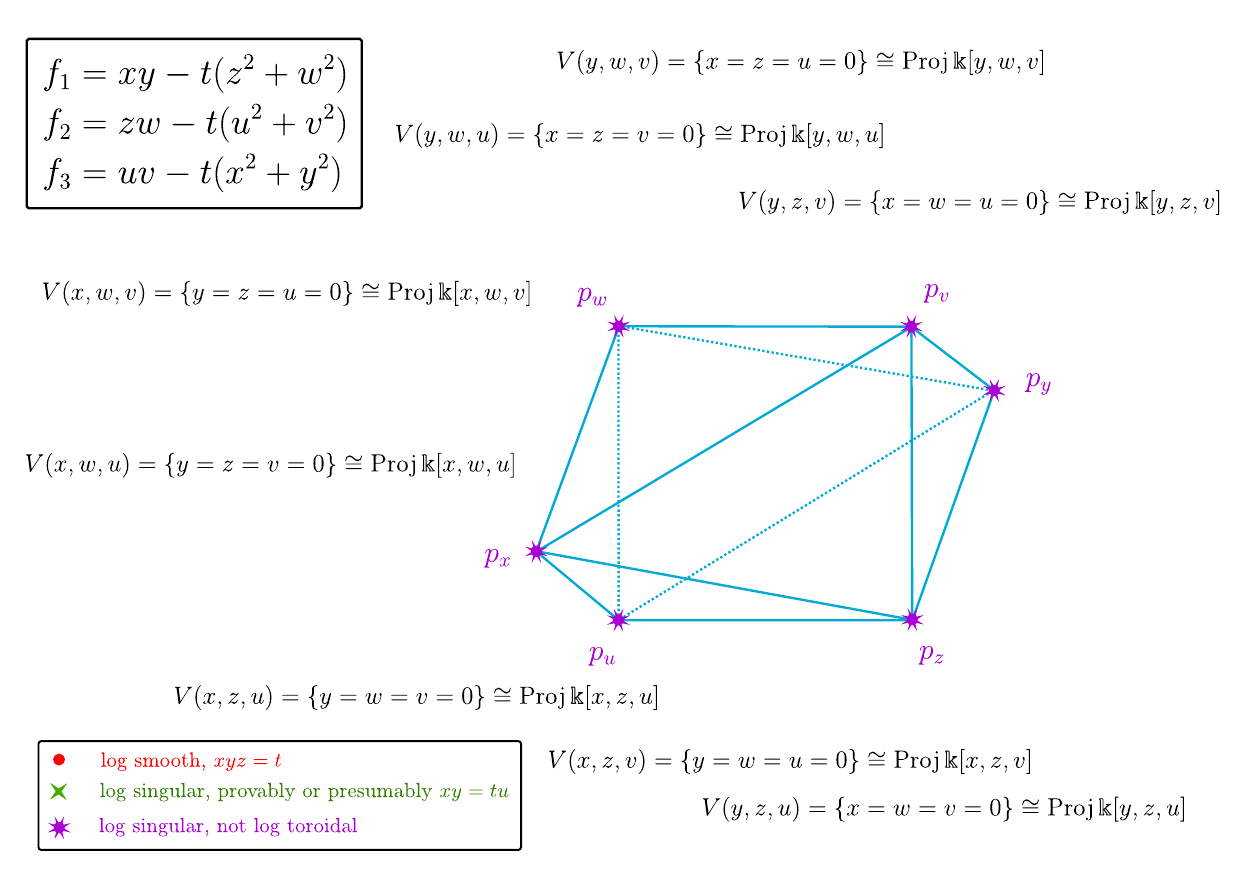}
		\caption{Example~\ref{ex:C23-OCT}}\label{fig:7-1}
	\end{mdframed}
\end{figure}

\begin{ex}[C23-OCT]\label{ex:C23-OCT}\note{ex:C23-OCT}
	We consider $R_\PP^\bullet = \QQ[t,x,y,z,w,u,v]$,
	\begin{align*}
		f_1 &= xy - t(z^2 + w^2), \\
		f_2 &= zw - t(u^2 + v^2), \\
		f_3 &= uv - t(x^2 + y^2),
	\end{align*}
	and $d = 2$. The fiber $\scX_{t = 1} = \phi^{-1}(1) \subset \PP^5$ is a smooth octic K3 surface so that Assumption~\ref{ass:smooth-fiber} is satisfied. The tools of \cite{Felten-gls-gr,FeltenZenodo} show that Assumptions~\ref{ass:apgls}, \ref{ass:Tx-small} via Strategy~\ref{strat:Tx-sing}, and \ref{ass:agls} are satisfied. Since $d + r = 5 = n$, Lemma~\ref{lem:log-canonical-ci} yields that Assumption~\ref{ass:Gorenstein} is satisfied, and that $\cA_{\widebar\scX/\widebar\scS}^{d,\m{algo}} \cong \cO_{\widebar\scX}$.
	
	The central fiber $V = X_0$ has $8$ irreducible components, each a copy of $\PP^2$. For their intersection pattern, see Figure~\ref{fig:7-1}. The reduced singular locus $D = \m{Sing}(V)$ has $12$ irreducible components, each a copy of $\PP^1$. The reduced singular locus $T = \m{Sing}(D)$ consists of six points $p_x$, $p_y$, $p_z$, $p_w$, $p_u$, $p_v$. At each of them, four irreducible components of $V$ intersect. Moreover, $V$ is locally of the form $\Spec \kk[x,y]/(xy) \times \Spec \kk[z,w]/(zw)$ near these points, so $V$ has toroidal crossing singularities.
	
	For every $\scE \in [\scD^\circ]$, the kink is $\ell(\scE) = 1$. We have $\scZ(x,z)^\circ = \{xz = 0\} = \varnothing$, and similarly, we have $\scZ^\scE = \varnothing$ for all $\scE \in [\scD^\circ]$.
	
	In the notation of Section~\ref{sect:non-toroidal-log-sing}, the radical of the annihilator of $Q = \m{coker}(E_1 \to E_1^{\vee\vee})$ cuts out the union of the six lines (affine cones) over the six points of $T = \{p_x,p_y,p_z,p_w,p_u,p_v\}$. Therefore, $\phi_1\colon X_1^\partial \to S_1^\dagger$ is not log toroidal at these six points.
	
	The (log) Hodge numbers are
	\begin{center}
		\renewcommand{\arraystretch}{1.2}
		\begin{tabular}{l|ccc}
			$t = 0$ & $h^0$ & $h^1$ & $h^2$ \\ \hline
			$\cA^2_{X/S}|_{X_0} \cong \cO_{X_0}(0)$ & 1 & 0 & 1 \\ \hline
			$\cA^1_{X/S}|_{X_0} \cong \Theta^{1,\m{log}}_{X/S}|_{X_0} \otimes \cO_{X_0}(0)$ & 0 & 20 & 0 \\ \hline
			$\cA^0_{X/S}|_{X_0} \cong \cO_{X_0}$ & 1 & 0 & 1 \\
		\end{tabular}
		\hspace{1.6cm}
		\begin{tabular}{l|ccc}
			$t = 1$ & $h^0$ & $h^1$ & $h^2$ \\ \hline
			$\Omega^2_{X_{t = 1}}$& 1 & 0 & 1 \\ \hline
			$\Omega^1_{X_{t = 1}}$ & 0 & 20 & 0 \\ \hline
			$\cO_{X_{t = 1}}$ & 1 & 0 & 1 \\
		\end{tabular}
	\end{center}
	for $t = 0$ and $t = 1$. Thus, $\phi_0\colon X_0^\m{en} \to S_0^\dagger$ is a new example of an unobstructed log Calabi--Yau surface. \qedex
\end{ex}

\begin{figure}[t]
	\begin{mdframed}
		\includegraphics[scale=.64]{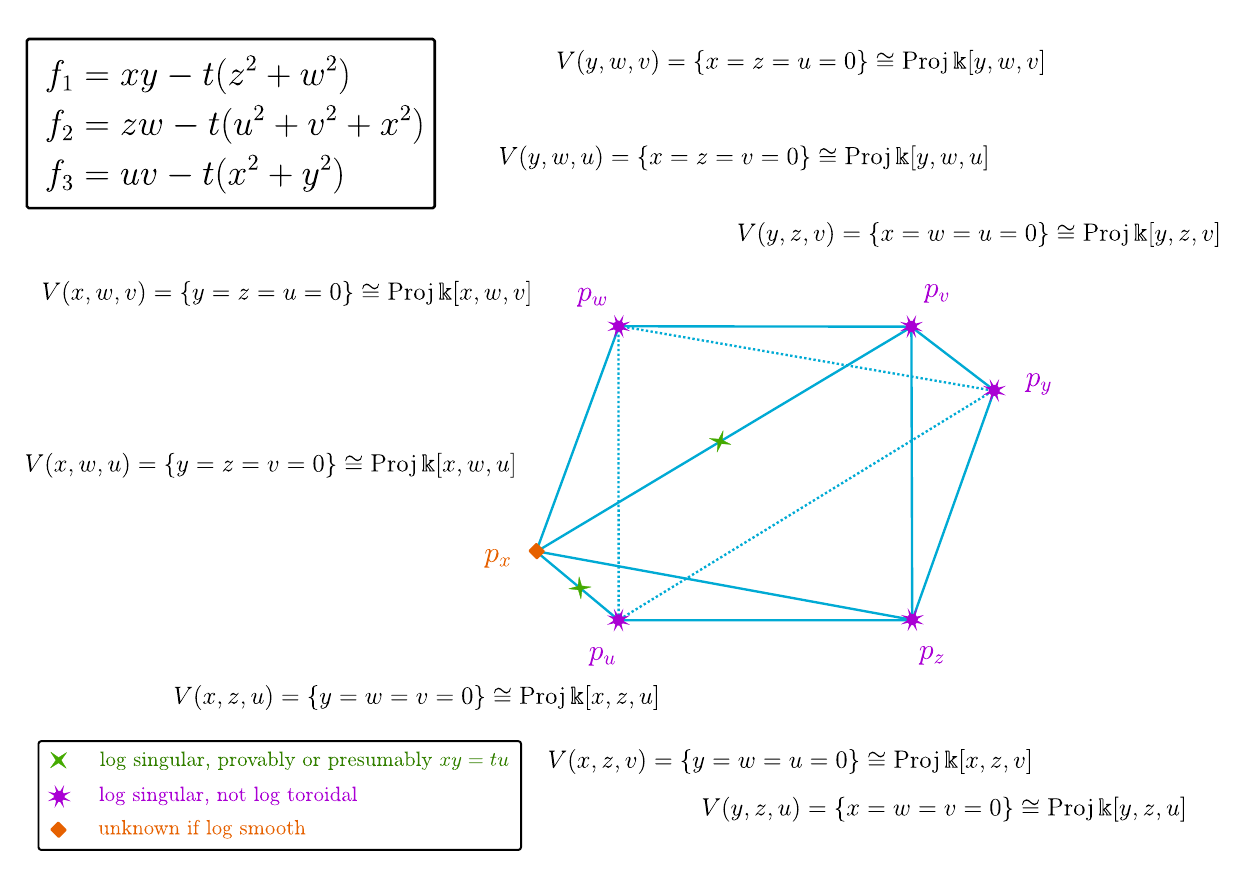}
		\caption{Example~\ref{ex:C23-ZW3}}\label{fig:C23-ZW3}
	\end{mdframed}
\end{figure}

\begin{ex}[C23-ZW3]\label{ex:C23-ZW3}\note{ex:C23-ZW3} 
  Note that this example is currently not reproducible on the 
  official OSCAR development branch; see Remark \ref{rem:NonReproducibleExample}. 
	We consider $R_\PP^\bullet = \QQ[t,x,y,z,w,u,v]$,
	\begin{align*}
		f_1 &= xy - t(z^2 + w^2), \\
		f_2 &= zw - t(u^2 + v^2 + x^2), \\
		f_3 &= uv - t(x^2 + y^2),
	\end{align*}
	and $d = 2$. The fiber $\scX_{t = 1} = \phi^{-1}(1) \subset \PP^5$ is a smooth octic K3 surface so that Assumption~\ref{ass:smooth-fiber} is satisfied. The tools of \cite{Felten-gls-gr,FeltenZenodo} show that Assumptions~\ref{ass:apgls}, \ref{ass:Tx-small} via Strategy~\ref{strat:Tx-sing}, and \ref{ass:agls} are satisfied. Since $d + r = 5 = n$, Lemma~\ref{lem:log-canonical-ci} yields that Assumption~\ref{ass:Gorenstein} is satisfied, and that $\cA_{\widebar\scX/\widebar\scS}^{d,\m{algo}} \cong \cO_{\widebar\scX}$.
	
	The central fiber $V = X_0$ is the same as in Example~\ref{ex:C23-OCT}. For every $\scE \in [\scD^\circ]$, the kink is $\ell(\scE) = 1$. We have $\scZ(r,s)^\circ = \{rs = 0\} = \varnothing$ for $r,s \in \{x,y,z,w,u,v\}$ unless $\{r,s\} = \{x,u\}$ or $\{r,s\} = \{x,v\}$; in this case, we have $\scZ(x,s)^\circ = \{xs(x + s) = 0\} = \{x + s = 0\}$ for $s \in \{u,v\}$.
	
	In the notation of Section~\ref{sect:non-toroidal-log-sing}, the radical of the annihilator of $Q_1 = \m{coker}(E_1 \to E_1^{\vee\vee})$ cuts out the union of the five lines (affine cones) over the five points of $T \setminus \{p_x\} = \{p_y,p_z,p_w,p_u,p_v\}$. Therefore, $\phi_1\colon X_1^\partial \to S_1^\dagger$ is not log toroidal at these five points. We do not know if the family is log smooth at $p_x$.
	
	The (log) Hodge numbers are
	\begin{center}
		\renewcommand{\arraystretch}{1.2}
		\begin{tabular}{l|ccc}
			$t = 0$ & $h^0$ & $h^1$ & $h^2$ \\ \hline
			$\cA^2_{X/S}|_{X_0} \cong \cO_{X_0}(0)$ & 1 & 0 & 1 \\ \hline
			$\cA^1_{X/S}|_{X_0} \cong \Theta^{1,\m{log}}_{X/S}|_{X_0} \otimes \cO_{X_0}(0)$ & 0 & 20 & 0 \\ \hline
			$\cA^0_{X/S}|_{X_0} \cong \cO_{X_0}$ & 1 & 0 & 1 \\
		\end{tabular}
		\hspace{1.6cm}
		\begin{tabular}{l|ccc}
			$t = 1$ & $h^0$ & $h^1$ & $h^2$ \\ \hline
			$\Omega^2_{X_{t = 1}}$& 1 & 0 & 1 \\ \hline
			$\Omega^1_{X_{t = 1}}$ & 0 & 20 & 0 \\ \hline
			$\cO_{X_{t = 1}}$ & 1 & 0 & 1 \\
		\end{tabular}
	\end{center}
	for $t = 0$ and $t = 1$, computed over $\kq = \FF_{4021}$ in the case of $t = 0$.\footnote{The code file is:
	\begin{quote}
		\texttt{Hodge-numbers-C23-ZW3-FF4021-cluster.jl}
	\end{quote}
	We noticed after 67 hours of runtime that the computation was complete, but the actual time until completion likely was much shorter.} If the numbers are the same over $\kq = \QQ$, then $\phi_0\colon X_0^\m{en} \to S_0^\dagger$ is a new example of an enhanced generically log smooth family with unobstructed deformations.
	\qedex
\end{ex}

\begin{figure}[t]
	\begin{mdframed}
		\includegraphics[scale=.64]{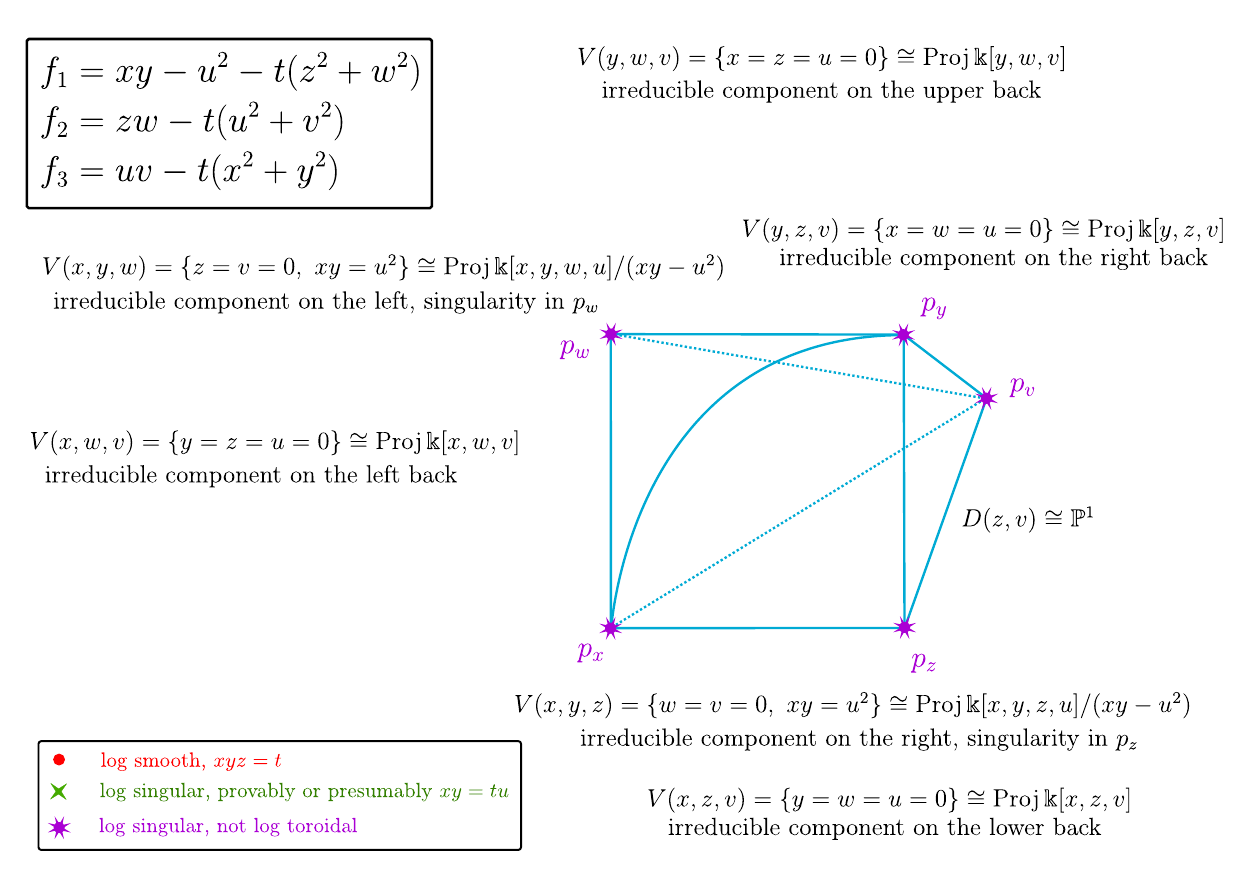}
		\caption{Example~\ref{ex:C23-DCP}}\label{fig:7-2}
	\end{mdframed}
\end{figure}

\begin{ex}[C23-DCP]\label{ex:C23-DCP}\note{ex:C23-DCP}
	We consider $R_\PP^\bullet = \QQ[t,x,y,z,w,u,v]$,
	\begin{align*}
		f_1 &= xy - u^2 - t(z^2 + w^2), \\
		f_2 &= zw - t(u^2 + v^2), \\
		f_3 &= uv - t(x^2 + y^2),
	\end{align*}
	and $d = 2$. The fiber $\scX_{t = 1} = \phi^{-1}(1) \subset \PP^5$ is a smooth octic K3 surface so that Assumption~\ref{ass:smooth-fiber} is satisfied. The tools of \cite{Felten-gls-gr,FeltenZenodo} show that Assumptions~\ref{ass:apgls}, \ref{ass:Tx-small} via Strategy~\ref{strat:Tx-sing}, and \ref{ass:agls} are satisfied. Since $d + r = 5 = n$, Lemma~\ref{lem:log-canonical-ci} yields that Assumption~\ref{ass:Gorenstein} is satisfied, and that $\cA_{\widebar\scX/\widebar\scS}^{d,\m{algo}} \cong \cO_{\widebar\scX}$.
	
	The central fiber $V = X_0$ has six irreducible components. See Figure~\ref{fig:7-2} for their intersection pattern. Four of them are degree-$1$ copies of $\PP^2$. The other two components $V(x,y,z)$ and $V(x,y,w)$ are projective surfaces of degree $2$. The reduced singular locus $D = \m{Sing}(V)$ consists of $9$ smooth curves, $8$ of which are degree-$1$ copies of $\PP^1$. All of them are the intersection of two irreducible components of $V$. The reduced singular locus $T = \m{Sing}(D)$ consists of five points $p_x$, $p_y$, $p_z$, $p_w$, and $p_v$. At $p_x$, $p_y$, and $p_v$, four irreducible components of $V$ intersect; at $p_z$ and $p_w$, only three irreducible components intersect. The component $V(x,y,z)$ has a normal surface singularity at $p_z$, and the component $V(x,y,w)$ has a normal surface singularity at $p_w$.
	
	At $p_x$, $p_y$, and $p_v$, the local form of $V$ is $\Spec \kk[x,y]/(xy) \times \Spec \kk[z,w]/(zw)$; at $p_z$ and $p_w$, the local form of $V$ is $\Spec \kk[x,y,z,w]/(xy - w^2, \, zw)$. Thus, $V$ has toroidal crossing singularities.
	
	For every $\scE \in [\scD^\circ]$, the kink is $\ell(\scE) = 1$. For the special component $\scD(x,y) = \scV(x,y,z) \cap \scV(x,y,w)$ of $\scD$, we have $\scZ(x,y)^\circ = \{u = 0\} = \varnothing$. Next, we have $\scZ(x,z)^\circ = \{xz = 0\} = \varnothing$. Similar formulae hold for the remaining seven components, and there is no log singular point in $D^\circ$.
	
	In the notation of Section~\ref{sect:non-toroidal-log-sing}, the radical of the annihilator of $Q_1 = \m{coker}(E_1 \to E_1^{\vee\vee})$ cuts out the union of the five lines (affine cones) over the five points of $T = \{p_x,p_y,p_z,p_w,p_v\}$. Therefore, $\phi_1\colon X_1^\partial \to S_1^\dagger$ is not log toroidal at these five points.
	
	The (log) Hodge numbers are
	\begin{center}
		\renewcommand{\arraystretch}{1.2}
		\begin{tabular}{l|ccc}
			$t = 0$ & $h^0$ & $h^1$ & $h^2$ \\ \hline
			$\cA^2_{X/S}|_{X_0} \cong \cO_{X_0}(0)$ & 1 & 0 & 1 \\ \hline
			$\cA^1_{X/S}|_{X_0} \cong \Theta^{1,\m{log}}_{X/S}|_{X_0} \otimes \cO_{X_0}(0)$ & 0 & 20 & 0 \\ \hline
			$\cA^0_{X/S}|_{X_0} \cong \cO_{X_0}$ & 1 & 0 & 1 \\
		\end{tabular}
		\hspace{1.6cm}
		\begin{tabular}{l|ccc}
			$t = 1$ & $h^0$ & $h^1$ & $h^2$ \\ \hline
			$\Omega^2_{X_{t = 1}}$& 1 & 0 & 1 \\ \hline
			$\Omega^1_{X_{t = 1}}$ & 0 & 20 & 0 \\ \hline
			$\cO_{X_{t = 1}}$ & 1 & 0 & 1 \\
		\end{tabular}
	\end{center}
	for $t = 0$ and $t = 1$. In particular, $\phi_0\colon X_0^\m{en} \to S_0^\dagger$ is a new example of a log Calabi--Yau surface with unobstructed deformations. \qedex
\end{ex}

\section{Threefolds}\label{sect:threefolds}\note{sect:threefolds}

The following examples illustrate that $\cA^\bullet_{X/S} = \cW^\bullet_{X/S}$ is usually not cffc for degenerations of threefolds, unless $\phi\colon X^\partial \to S^\dagger$ has toroidal log singularities.

\subsection{Calabi--Yau threefolds}

\begin{ex}[C31-FQI]\label{ex:C31-FQI}\note{ex:C31-FQI}
	We consider $R_\PP^\bullet = \QQ[t,x,y,z,w,u]$,
	\begin{align*}
		f_1 &= xyzwu - t(x^5 + y^5 + z^5 + w^5 + u^5),
	\end{align*}
	and $d = 3$. The fiber $\widebar\scX_{t = 1} = \phi^{-1}(1) \subset \PP^4$ is smooth so that Assumption~\ref{ass:smooth-fiber} is satisfied. The tools of \cite{Felten-gls-gr,FeltenZenodo} show that Assumptions~\ref{ass:apgls}, \ref{ass:Tx-small} via Strategy~\ref{strat:Tx-sing}, and \ref{ass:agls} are satisfied. Since $d + r = 4 = n$, Lemma~\ref{lem:log-canonical-ci} yields that Assumption~\ref{ass:Gorenstein} is satisfied, and that $\cA_{\widebar\scX/\widebar\scS}^{d,\m{algo}} \cong \cO_{\widebar\scX}$.
	
	The family $\phi\colon X^\partial \to S^\dagger$ is log toroidal. To see this, first note that $X_0$ is a normal crossing scheme. After computing the singular locus of $\widebar X$, we see that none of the five points in the deepest strata of $X_0$ lies in $\m{Sing}(\widebar X)$. Therefore, the family is semistable at these points. Furthermore, on the regular locus $\m{Reg}(X_0) \subseteq X_0$, the family is smooth, hence log smooth.
	
	Now, let $u \not= 0$ and consider the affine chart $\Spec R$ of $\widebar X$ with
	\[R = \kq[t,x,y,z,w]/(xyzw - t(x^5 + y^5 + z^5 + w^5 + 1)).\]
	Let $A = \kq[\tilde t,\tilde x,\tilde y,\tilde z,\tilde u]/(\tilde x\tilde y\tilde z -\tilde t\tilde u)$, and consider the ring homomorphism
	\begin{align*}
		\psi\colon \enspace A[\tilde r,\tilde s] &\to R, \quad \tilde t \mapsto t, \quad \tilde r \mapsto z, \quad \tilde s \mapsto w, \\
		\tilde x &\mapsto x, \quad \tilde y \mapsto y, \quad \tilde z \mapsto zw, \quad \tilde u \mapsto x^5 + y^5 + z^5 + w^5 + 1.
	\end{align*} 
	It is surjective, and its kernel is generated by $f_1 = \tilde x^5 + \tilde y^5 + \tilde r^5 + \tilde s^5 + 1 - \tilde u$ and $f_2 = \tilde r\tilde s - \tilde z$. Thus, after localizing $R$ in the image $h = 5(z^5 - w^5)$ of the Jacobi determinant $\m{det}(J_\psi) = 5(\tilde r^5 - \tilde s^5)$, the map $\psi\colon A \to R_h$ becomes \'etale, cf.~\cite[Tag~00T8]{stacks}. Thus, $\phi\colon X^\partial \to S^\dagger$ is log toroidal on $\{u \not= 0, z^5 \not= w^5\}$.
	
	By symmetry, if $\phi\colon X^\partial \to S^\dagger$ has a non-toroidal log singularity on $\{u \not= 0\}$, it must lie in $\{u \not= 0, x^5 = y^5 = z^5 = w^5\}$. However, it must also lie in the singular locus $D = \m{Sing}(X_0)$ of the central fiber, so at least two variables are zero, and we have $x = y = z = w = 0$. However, we have already seen that $\phi\colon X \to S$ is semistable at this point. Again by symmetry, $\phi\colon X^\partial \to S^\dagger$ is log toroidal.

	The (log) Hodge numbers are 
	\begin{center}
		\renewcommand{\arraystretch}{1.2}
		\setlength{\tabcolsep}{4pt}
		\begin{tabular}{l|cccc}
			$t = 0$ & $h^0$ & $h^1$ & $h^2$ & $h^3$ \\ \hline
			$\cA^3_{X/S}|_{X_0} \cong \cO_{X_0}(0)$ & 1 & 0 & 0 & 1 \\ \hline
			$\cA^2_{X/S}|_{X_0} \cong \Theta^{1,\m{log}}_{X/S}|_{X_0} \otimes \cO_{X_0}(0)$ & 0 & 101 & 1 & 0 \\ \hline
		 	$\cA^1_{X/S}|_{X_0}$ & 0 & 1 & 101 & 0 \\ \hline
			$\cA^0_{X/S}|_{X_0} \cong \cO_{X_0}$ & 1 & 0 & 0 & 1 \\
		\end{tabular}
		\hspace{1cm}
		\begin{tabular}{l|cccc}
			$t = 1$ & $h^0$ & $h^1$ & $h^2$ & $h^3$ \\ \hline
			$\Omega^3_{X_{t = 1}}$ & 1 & 0 & 0 & 1 \\ \hline
			$\Omega^2_{X_{t = 1}}$ & 0 & 101 & 1 & 0 \\ \hline
			$\Omega^1_{X_{t = 1}}$ & 0 & 1 & 101 & 0 \\ \hline
			$\cO_{X_{t = 1}}$ & 1 & 0 & 0 & 1 \\
		\end{tabular}
	\end{center}
	as guaranteed by Theorem~\ref{thm:log-toroidal-cffc}.\qedex
\end{ex}

\begin{ex}[C31-US1]\label{ex:C31-US1}\note{ex:C31-US1}
	We consider $R_\PP^\bullet = \QQ[t,x,y,z,w,u]$,
	\begin{align*}
		f_1 &= xyzwu - t(xu^4 + x^5 + y^5 + z^5 + w^5),
	\end{align*}
	and $d = 3$. The fiber $\scX_{t = 1} = \phi^{-1}(1) \subset \PP^4$ is smooth so that Assumption~\ref{ass:smooth-fiber} is satisfied. The tools of \cite{Felten-gls-gr,FeltenZenodo} show that Assumptions~\ref{ass:apgls}, \ref{ass:Tx-small} via Strategy~\ref{strat:Tx-sing}, and \ref{ass:agls} are satisfied. Since $d + r = 4 = n$, Lemma~\ref{lem:log-canonical-ci} yields that Assumption~\ref{ass:Gorenstein} is satisfied, and that $\cA_{\widebar\scX/\widebar\scS}^{d,\m{algo}} \cong \cO_{\widebar\scX}$.
	
	In the notation of Section~\ref{sect:non-toroidal-log-sing}, the radical of the annihilator of $Q'_1 = \m{coker}(E'_1 \to E'^{\vee\vee}_1)$ is $(x,y,z,w)$ so that $p_u$ is a non-toroidal log singularity of $\phi_1\colon X_1^\partial \to S_1^\dagger$, for $\{p_u\} = \{t = x = y = z = w = 0\}$.
	
	The (log) Hodge numbers are 
	\begin{center}
		\renewcommand{\arraystretch}{1.2}
		\setlength{\tabcolsep}{4pt}
		\begin{tabular}{l|cccc}
			$t = 0$ & $h^0$ & $h^1$ & $h^2$ & $h^3$ \\ \hline
			$\cA^3_{X/S}|_{X_0} \cong \cO_{X_0}(0)$ & 1 & 0 & 0 & 1 \\ \hline
			$\cA^2_{X/S}|_{X_0} \cong \Theta^{1,\m{log}}_{X/S}|_{X_0} \otimes \cO_{X_0}(0)$ & 0 & 119 & 19 & 0 \\ \hline
			$\cA^1_{X/S}|_{X_0}$ & 0 & 1 & 101 & 0 \\ \hline
			$\cA^0_{X/S}|_{X_0} \cong \cO_{X_0}$ & 1 & 0 & 0 & 1 \\
		\end{tabular}
		\hspace{1cm}
		\begin{tabular}{l|cccc}
			$t = 1$ & $h^0$ & $h^1$ & $h^2$ & $h^3$ \\ \hline
			$\Omega^3_{X_{t = 1}}$ & 1 & 0 & 0 & 1 \\ \hline
			$\Omega^2_{X_{t = 1}}$ & 0 & 101 & 1 & 0 \\ \hline
			$\Omega^1_{X_{t = 1}}$ & 0 & 1 & 101 & 0 \\ \hline
			$\cO_{X_{t = 1}}$ & 1 & 0 & 0 & 1 \\
		\end{tabular}
	\end{center}
	so that $\cA^\bullet_{X/S}$ is not cffc. \qedex
\end{ex}

\begin{ex}[C32-BCU]\label{ex:C32-BCU}\note{ex:C32-BCU}
	We consider $R_\PP^\bullet = \QQ[t,x,y,z,w,u,v]$,
	\begin{align*}
		f_1 &= xyz - t(u^3 + v^3 + w^3), \\
		f_2 &= uvw - t(x^3 + y^3 + z^3), 
	\end{align*}
	and $d = 3$. The fiber $\scX_{t = 1} \subset \PP^5$ is smooth so that Assumption~\ref{ass:smooth-fiber} is satisfied. The tools of \cite{Felten-gls-gr,FeltenZenodo} show that Assumptions~\ref{ass:apgls}, \ref{ass:Tx-small} via Strategy~\ref{strat:Tx-sing}, and \ref{ass:agls} are satisfied. Since $d + r = 5 = n$, Lemma~\ref{lem:log-canonical-ci} yields that Assumption~\ref{ass:Gorenstein} is satisfied, and that $\cA_{\widebar\scX/\widebar\scS}^{d,\m{algo}} \cong \cO_{\widebar\scX}$.
	
	The central fiber $V$ is locally of the form $\Spec \kk[x,y,z]/(xyz) \times \Spec \kk[u,w]/(uw)$, so it has toroidal crossing singularities.
	
	In the notation of Section~\ref{sect:non-toroidal-log-sing}, the radical of the annihilator of $Q'_1 = \m{coker}(E'_1 \to E'^{\vee\vee}_1)$ cuts out the lines over the six points $\{p_x,p_y,p_z,p_w,p_u,p_v\}$. Thus, the six points are non-toroidal log singularities of $\phi_1\colon X_1^\partial \to S_1^\dagger$. 
	
	The (log) Hodge numbers are
	\begin{center}
		\renewcommand{\arraystretch}{1.2}
		\setlength{\tabcolsep}{4pt}
		\begin{tabular}{l|cccc}
			$t = 0$ & $h^0$ & $h^1$ & $h^2$ & $h^3$ \\ \hline
			$\cA^3_{X/S}|_{X_0} \cong \cO_{X_0}(0)$ & 1 & 0 & 0 & 1 \\ \hline
			$\cA^2_{X/S}|_{X_0} \cong \Theta^{1,\m{log}}_{X/S}|_{X_0} \otimes \cO_{X_0}(0)$ & 0 & 115 & 43 & 0 \\ \hline
			$\cA^1_{X/S}|_{X_0}$ & $*$ & $*$ & $*$ & $*$  \\ \hline
			$\cA^0_{X/S}|_{X_0} \cong \cO_{X_0}$ & 1 & 0 & 0 & 1 \\
		\end{tabular}
		\hspace{1cm}
		\begin{tabular}{l|cccc}
			$t = 1$ & $h^0$ & $h^1$ & $h^2$ & $h^3$ \\ \hline
			$\Omega^3_{X_{t = 1}}$ & 1 & 0 & 0 & 1 \\ \hline
			$\Omega^2_{X_{t = 1}}$ & 0 & 73 & 1 & 0 \\ \hline
			$\Omega^1_{X_{t = 1}}$ & 0 & 1 & 73 & 0 \\ \hline
			$\cO_{X_{t = 1}}$ & 1 & 0 & 0 & 1 \\
		\end{tabular}
	\end{center}
	for $t = 0$ and $t = 1$, where $*$ denotes a number that we did not determine. In particular, $\cA^\bullet_{X/S}$ is not cffc. \qedex
\end{ex}

\subsection{The class $\scC^\m{adm}$}

Recall from Section~\ref{sect:Fano-smoothing} the conjectural class of log singularities $\scC^\m{adm}$, which is supposed to play a crucial role in smoothing reducible toric Fano schemes $\PP_0(\widetilde\Delta)$. Also recall that the family in Example~\ref{ex:adm-model} is conjecturally contained in $\scC^\m{adm}$. Furthermore, we write $\scC^\m{tor}$ for the set of all toric local models $A_{\itPi/\NN}^\partial \to A_\NN^\dagger$ as in Definition~\ref{defn:toroidal-log-sing}. Here is a projective morphism of partial log schemes $\phi\colon X^\partial \to S^\dagger$ which is (presumably) of class $\scC^\m{adm} \cup \scC^\m{tor}$:

\begin{ex}[N32-ADM]\label{ex:N32-ADM}\note{ex:N32-ADM}
	We consider $R_\PP^\bullet = \QQ[t,x,y,z,w,u,v]$,
	\begin{align*}
		f_1 &= xy - z^2 + t^2(u^2 + v^2 + w^2), \\
		f_2 &= zw - t(uv - x^2 - y^2), 
	\end{align*}
	and $d = 3$. The fiber $\scX_{t = 1} = \phi^{-1}(1) \subset \PP^5$ is smooth so that Assumption~\ref{ass:smooth-fiber} is satisfied. The tools of \cite{Felten-gls-gr,FeltenZenodo} show that Assumptions~\ref{ass:apgls}, \ref{ass:Tx-small} via Strategy~\ref{strat:Tx-sing}, and \ref{ass:agls} are satisfied. Since $d + r = 5 = n$, Lemma~\ref{lem:log-canonical-ci} yields that Assumption~\ref{ass:Gorenstein} is satisfied, and that $\cA_{\widebar\scX/\widebar\scS}^{d,\m{algo}} \cong \cO_{\widebar\scX}(-2)$. 
	
	The geometry of this family was discussed by the first author in \cite[Ex.~2.24]{Felten-gls-gr}, where the reader can also find a graphic representation. The equations were, up to a minor modification, originally proposed by Alessio Corti as a projective family with log singularities modeled on Example~\ref{ex:adm-model}.
	
	The central fiber $V = X_0$ has three irreducible components $V(x,y,z)$, $V(x,w)$, and $V(y,w)$. Any two of them intersect in a copy of $\PP^2$, denoted by $D(x)$, $D(y)$, and $D(w)$. All three irreducible components intersect in $T \cong \PP^1$, which is also the complement of the double normal crossing locus. Locally near points in $T$, the scheme $V$ is of the form $\Spec \kk[x,y,z,w,u]/(xy - z^2, \, zw)$, so $V$ has toroidal crossing singularities.
	
	In \cite[Ex.~2.24]{Felten-gls-gr}, the first author computed a log indeterminate locus
	\[Z_X = Z(x) \cup Z(y) \cup Z(w) \cup T\]
	with $Z(x) = \{x^2 = uv\} \subset D(x)$, $Z(y) = \{y^2 = uv\} \subset D(y)$, and
	\[Z(w) = \{u^2v^2 - u^2w^2 - v^2w^2 - w^4 = 0\} \subset D(w).\]
	From the theory in \cite{Felten-gls-gr}, it follows that $Z(x) \cup Z(y) \cup Z(w) \subset \itSigma_X \subset Z_X$, where $\itSigma_X \subset X$ is the log singular locus. We will now study the local structure of $\phi\colon X^\partial \to S^\dagger$, which was not completed in \cite[Ex.~2.24]{FeltenGLDT}.
	
	Let
	\[q_0 = [0:0:0:0:0:1]\quad \text{and} \quad q_\infty = [0:0:0:0:1:0]\]
	be the points where all components of $Z_X$ intersect. We claim the following: The points $q_0$ and $q_\infty$ are non-toroidal log singularities modeled on Example~\ref{ex:adm-model}.
	
	We prove it for $q_0$, the case of $q_\infty$ being analogous. Since $q_0 \in \{v \not= 0\}$, let us consider 
	\[R = \QQ[t,x,y,z,w,u]/(xy - z^2 + t^2(u^2 + 1 + w^2), \: zw - t(u - x^2 - y^2)).\]
	Then, we set $S = R[a^\pm]/(a^2 + u^2 + 1 + w^2)$. The map $\chi\colon R \to S$ is \'etale, cf.~\cite[Tag~00T8]{stacks}. After setting
	\[\{q_0'\} = \{t = x = y = z = w = u = 0,\: a^2 + 1 = 0\} \subset \Spec S,\]
	we have $\Spec(\chi)(q_0') = q_0$. Note that $q_0'$ is a $\QQ[i]$-valued point.
	
	Now, let 
	\[A = \QQ[\tilde t,\tilde x,\tilde y,\tilde z,\tilde w,\tilde u]/(\tilde x\tilde y - \tilde z^2 - \tilde t^2, \: \tilde z \tilde w - \tilde t\tilde u),\]
	and consider the ring homomorphism
	\begin{align*}
		\psi\colon \enspace A[\tilde a, \tilde b] &\to S, \qquad \tilde t \mapsto t, \quad  \tilde a \mapsto a, \quad \tilde b \mapsto \frac{1}{a} \\
		& \tilde x \mapsto \frac{x}{a}, \quad \tilde y \mapsto \frac{y}{a}, \quad \tilde z \mapsto \frac{z}{a}, \quad \tilde w \mapsto wa, \quad \tilde u \mapsto u - x^2 - y^2.
	\end{align*}
	It is surjective, and its kernel is generated by $f_1 = \tilde a\tilde b - 1$ and 
	\[f_2 = \tilde x^4\tilde a^2 + 2\tilde x^2\tilde y^2\tilde a^2 + \tilde y^4\tilde a^2 + \tilde w^2\tilde b^4 + \tilde u^2\tilde b^2 + 2\tilde x^2\tilde u + 2\tilde y^2\tilde u + \tilde b^2 + 1.\]
	We check this with Macaulay2 in the file \texttt{geometry-N32-ADM.m2} in \cite{FZ26}. The Jacobi determinant is 
	\[\m{det}(J_\psi) = 2\tilde b^2(2\tilde b^2\tilde w^2 + \tilde u^2 + 1) - 2\tilde a^2(\tilde x^2 + \tilde y^2)^2.\]
	Its image under $\psi$ is 
	\[h = \frac{2}{a^2}(2a^2 + w^2 - 2u(x^2 + y^2)).\]
	Since $h(q_0') = 4 \in \QQ[i]$, the map $\Spec S \to \Spec A$ induced by $\psi$ is \'etale in a neighborhood of $q_0'$.
	
	We leave it to the interested reader to determine the local structure at the other points of $Z_X$. We believe that the log singular locus is $\itSigma_X = Z(x) \cup Z(y) \cup Z(w)$, and that every log singularity other than $q_0$, $q_\infty$ is toroidal of type $T(1;1)$ or $T(2;1)$. In particular, we believe that $\phi\colon X^\partial \to S^\dagger$ is of the conjectural class $\scC^\m{adm} \cup \scC^\m{tor}$.

	The (log) Hodge numbers are
		\begin{center}
		\renewcommand{\arraystretch}{1.2}
		\setlength{\tabcolsep}{4pt}
		\begin{tabular}{l|cccc}
			$t = 0$ & $h^0$ & $h^1$ & $h^2$ & $h^3$ \\ \hline
			$\cA^3_{X/S}|_{X_0} \cong \cO_{X_0}(-2)$ & 0 & 0 & 0 & 1 \\ \hline
			$\cA^2_{X/S}|_{X_0} \cong \Theta^{1,\m{log}}_{X/S}|_{X_0} \otimes \cO_{X_0}(-2)$ & 0 & 2 & 1 & 0 \\ \hline
			$\cA^1_{X/S}|_{X_0}$ & $*$ & $*$ & $*$ & $*$  \\ \hline
			$\cA^0_{X/S}|_{X_0} \cong \cO_{X_0}$ & 1 & 0 & 0 & 0 \\
		\end{tabular}
		\hspace{1cm}
		\begin{tabular}{l|cccc}
			$t = 1$ & $h^0$ & $h^1$ & $h^2$ & $h^3$ \\ \hline
			$\Omega^3_{X_{t = 1}}$ & 0 & 0 & 0 & 1 \\ \hline
			$\Omega^2_{X_{t = 1}}$ & 0 & 2 & 1 & 0 \\ \hline
			$\Omega^1_{X_{t = 1}}$ & 0 & 1 & 2 & 0 \\ \hline
			$\cO_{X_{t = 1}}$ & 1 & 0 & 0 & 0 \\
		\end{tabular}
	\end{center}
	for $t = 0$ and $t = 1$, as far as we computed them (using $\kq = \QQ$ for $t = 0$).\footnote{The code file of a failed attempt to compute $\cA^{1,\m{algo}}_{\widebar\scX/\widebar\scS}$ over $\kq = \FF_{4021}$ can be found in \cite{FZ26}:
	\begin{quote}
		\texttt{Hodge-numbers-N32-ADM-FF4021-cluster.jl}
	\end{quote}
	We aborted the attempt after $> 3$ weeks of runtime.}
	
	While we have not been able to compute $\cA^1_{X_0/S_0} = \cA^1_{X/S}|_{X_0}$, we have a short exact sequence 
	\[0 \to \cA^1_{X/S}|_{X_0} \to \cW^1_{X_0/S_0} \to \cQ_{X_0} \to 0,\]
	where $\cW^1_{X_0/S_0} = j_*\Omega^{1,\m{log}}_{U_0/S_0} \cong \cH om(\Theta^{1,\m{log}}_{X/S}|_{X_0},\cO_{X_0})$. Injectivity on the left follows from the fact that $\m{depth}(\cA_{X/S,p}) \geq 2$ at every point $p \in X$ with $\m{depth}(\cO_{X,p}) \geq 2$. The cohomology of $\cW^1_{X_0/S_0}$ is $[0\ 1\ 4\ 0]$.
	
	In the local model $Y^\partial \to A_\NN^\dagger$ from Example~\ref{ex:adm-model}, we can compute $\cQ_{Y_0}$ and obtain that it is isomorphic to $\cO_{Y_0}/(x,y,z,w,u)$ as an $\cO_{Y_0}$-module. If indeed all log singularities except $q_0$, $q_\infty$ are toroidal, then $\cQ_{X_0}|_{X_0 \setminus \{q_0,q_\infty\}} = 0$, and $\cQ_{X_0}$ is isomorphic to $\cO_{X_0}/\cI_{\{q_0,q_\infty\}}$ as an $\cO_{X_0}$-module. In this case, the cohomology of $\cQ_{X_0}$ is $[2\ 0\ 0\ 0]$, and the long exact sequence yields that the cohomology of $\cA^1_{X_0/S_0}$ is $[0\ 3\ 4\ 0]$. In particular, if our assumption is correct, $\cA^\bullet_{X/S}$ is not cffc. \qedex
\end{ex}

Now that $\scC^\m{adm} \cup \scC^\m{tor}$ does not seem to work, we may speculate that an alternative class of log singularities $\scC^\m{alt}$ is better. For inspiration, we look at \cite[\S1.12.3]{FeltenGLDT}, where the first author considers $\PP_0(\widetilde\Delta)$ for a specific choice of $\Delta \subset \RR^3$ and constructs a partial log smooth structure $\phi_0\colon \PP_0(\widetilde\Delta)^\partial \to S_0^\dagger$ with non-toroidal log singularities of the following form.

\begin{ex}[A32-CSP]\label{ex:almost-adm-model}\note{ex:almost-adm-model}
	We consider the affine family $\mu\colon Y = \Spec R \to \IA^1_t$ for 
	\[R = \kk[t,x,y,z,w,u]/(xy - z^2 - t^2w, \: zw - tu).\]
	For a graphic, see \cite[Ex.~1.159]{FeltenGLDT}. The central fiber $Y_0$ is the only singular fiber, and it is the same as in Example~\ref{ex:adm-model}. The methods of \cite{Felten-gls-gr} yield a log indeterminate locus 
	\[Z_Y = Z(x) \cup Z(y) \cup Z(w) \cup T\]
	which makes $\mu\colon Y^\partial \to A_\NN^\dagger$ generically log smooth. Here, $T \subset Y_0$ is the triple locus, and we have $Z(x) = \{u = 0\} \subset D(x)$, $Z(y) = \{u = 0\} \subset D(y)$, and $Z(w) = \{w^3 + u^2 = 0\} \subset D(w)$.
	
	It follows from the methods of \cite{Felten-gls-gr} that $\mu\colon Y^\dagger \to A_\NN^\dagger$ is indeed log singular at $p \in Z_Y \setminus T$. We claim that $\mu\colon Y^\dagger \to A_\NN^\dagger$ is log smooth at $p \in T \setminus \{0\}$. To see this, we consider 
	\begin{align*}
		\psi\colon \enspace &R_{u(u^2 + w^3)} \to \kk[\tilde t,\tilde x,\tilde y,\tilde z,\tilde w,\tilde u]_{\tilde u(1 + \tilde w^3\tilde u))}/(\tilde x\tilde y - \tilde z^2, \: \tilde z\tilde w - \tilde t), \\
		t &\mapsto \tilde t, \quad x \mapsto \tilde x, \quad y \mapsto \tilde y(1 + \tilde w^3\tilde u), \quad z \mapsto \tilde z, \quad w \mapsto \tilde w\tilde u, \quad u \mapsto \tilde u.
	\end{align*}
	One may easily check that this is a well-defined surjective ring homomorphism. Since both are integral rings of the same Krull dimension, $\psi$ is an isomorphism. We have $T \setminus \{0\} \subset \{u(w^3 + u^2) \not= 0\}$, and the target of $\psi$ corresponds to a log smooth and saturated morphism, see \cite[Ex.~1.66]{FeltenGLDT}. Therefore, the log singular locus is $\itSigma_Y = Z(x) \cup Z(y) \cup Z(w)$. From now on, we use $\itSigma_Y$ as the log indeterminate locus of $\mu\colon Y^\partial \to A_\NN^\dagger$.
	
	It is immediately clear that $\mu\colon Y^\partial \to A_\NN^\dagger$ has a toroidal log singularity of type $T(1;1)$ at every $p \in Z(x) \setminus \{0\}$ and $p \in Z(y) \setminus \{0\}$. Furthermore, it is not difficult to show that every $p \in Z(w) \setminus \{0\}$ is a toroidal log singularity of type $T(2;1)$, by writing down a suitable local model explicitly, in analogy with our discussion in Example~\ref{ex:C31-FQI}. The point $0 \in Y$ is a non-toroidal log singularity since $\cA_{Y_0/S_0}^1$ is not reflexive at $0$. \qedex
\end{ex}

In this local model, the branch $Z(w)$ of the log singular locus is a cuspidal cubic with the cusp lying in the triple locus $T$. Therefore, the central fiber $\mu_0\colon Y_0^\partial \to S_0^\dagger$ is not in $\scC_0^\m{adm}$. Namely, \cite[Defn.~12]{CGR2025} requires that all irreducible components of the log indeterminate locus are smooth.

\begin{ex}[N32-CSP]\label{ex:N32-CSP}\note{ex:N32-CSP}
	We consider $R_\PP^\bullet = \QQ[t,x,y,z,w,u,v]$,
	\begin{align*}
		f_1 &= xy - z^2 + t^2(uv + uw + vw), \\
		f_2 &= zw - t(uv - x^2 - y^2), 
	\end{align*}
	and $d = 3$, a variant of Example~\ref{ex:N32-ADM}. The fiber $\scX_{t = 1} \subset \PP^5$ is smooth so that Assumption~\ref{ass:smooth-fiber} is satisfied. The tools of \cite{Felten-gls-gr,FeltenZenodo} show that Assumptions~\ref{ass:apgls}, \ref{ass:Tx-small} via Strategy~\ref{strat:Tx-sing}, and \ref{ass:agls} are satisfied. Since $d + r = 5 = n$, Lemma~\ref{lem:log-canonical-ci} yields that Assumption~\ref{ass:Gorenstein} is satisfied, and that $\cA_{\widebar\scX/\widebar\scS}^{d,\m{algo}} \cong \cO_{\widebar\scX}(-2)$. 
	
	The central fiber is the same as in Example~\ref{ex:N32-ADM}. The log indeterminate locus is 
	\[Z_X = Z(x) \cup Z(y) \cup Z(w) \cup T\]
	with $Z(x) = \{x^2 = uv\} \subset D(x)$, $Z(y) = \{y^2 = uv\} \subset D(y)$, and
	\[Z(w) = \{u^2v^2 - w^3u - w^3v - w^2uv = 0\} \subset D(w).\]
	The log indeterminate locus has two special points $q_0$ and $q_\infty$, as in Example~\ref{ex:N32-ADM}.
	
	We claim that the log singularities in $q_0$ and $q_\infty$ are locally modeled on Example~\ref{ex:almost-adm-model}. By symmetry, it suffices to show this for $q_0$. Let 
	\[R = \Spec \QQ[t,x,y,z,w,u]/(xy - z^2 + t^2(u + w + uw), \: zw - t(u - x^2 - y^2))\]
	so that $q_0 = 0 \in \Spec R$. The following comparison maps were suggested by AI.\footnote{We used ChatGPT but did not record the version.} After setting $c = 1 + w$ and $a = t^2c$, we consider 
	\[\chi\colon \enspace R \to S \coloneq R[r]_{(1 + r^2)}/(r - a(1 + r^2)).\]
	This map is \'etale at $q'_0 = 0 \in \Spec S$, and we have $\Spec(\chi)(q_0') = q_0$.
	
	Next, let 
	\[A = \QQ[\tilde t,\tilde x,\tilde y,\tilde z,\tilde w,\tilde u]/(\tilde x\tilde y - \tilde z^2 - \tilde t^2\tilde w,\: \tilde z\tilde w - \tilde t\tilde u)\]
	be the local model. After setting $b = u - x^2 - y^2$, $d = 1 + \frac{wc^2}{4}$, and $s = \frac{1}{1 + r^2}$, we obtain a ring homomorphism 
	\begin{align*}
		\psi\colon \enspace &A \to S, \quad \tilde t \mapsto t, \quad \tilde x \mapsto x + ry, \quad \tilde y \mapsto s(y + rx), \\
		\tilde z &\mapsto z - {\frac{twc}{2}}, \quad \tilde w \mapsto -wd, \quad \tilde u \mapsto d\left({\frac{w^2c}{2}} - b\right).
	\end{align*}
	One may check by a direct computation that it is well-defined. Furthermore, we have $\Spec(\psi)(q_0') = 0 \in \Spec A$.
	
	We show that $\Spec(\psi)$ is \'etale at $q_0'$. First, a computation in Macaulay2 shows that both $A$ and $S$ are integral and normal of dimension $4$. Then, we compute $\Omega^1_{S/A}$; the tensor product along $S \to \QQ$ corresponding to the inclusion $\{q_0'\} \to \Spec S$ is zero, so we find that $\Spec(\psi)$ is unramified at $q_0'$. In particular, we can find an open neighborhood $q_0' \in U \subseteq \Spec S$ such that $U \to \Spec A$ is unramified, hence quasi-finite. But then, $U \to \Spec A$ is dominant, and \cite[Tag~0GS9]{stacks} yields that $U \to \Spec A$ is \'etale.
	
	We leave it to the interested reader to determine the local structure at the other points of $Z_X$. We believe that the log singular locus is $\itSigma_X = Z(x) \cup Z(y) \cup Z(w)$, and that all log singularities other than $q_0$, $q_\infty$ are toroidal of type $T(1;1)$ or $T(2;1)$.
	
	The (log) Hodge numbers are
	\begin{center}
		\renewcommand{\arraystretch}{1.2}
		\setlength{\tabcolsep}{4pt}
		\begin{tabular}{l|cccc}
			$t = 0$ & $h^0$ & $h^1$ & $h^2$ & $h^3$ \\ \hline
			$\cA^3_{X/S}|_{X_0} \cong \cO_{X_0}(-2)$ & 0 & 0 & 0 & 1 \\ \hline
			$\cA^2_{X/S}|_{X_0} \cong \Theta^{1,\m{log}}_{X/S}|_{X_0} \otimes \cO_{X_0}(-2)$ & 0 & 3 & 2 & 0 \\ \hline
			$\cA^1_{X/S}|_{X_0}$ & $*$ & $*$ & $*$ & $*$  \\ \hline
			$\cA^0_{X/S}|_{X_0} \cong \cO_{X_0}$ & 1 & 0 & 0 & 0 \\
		\end{tabular}
		\hspace{1cm}
		\begin{tabular}{l|cccc}
			$t = 1$ & $h^0$ & $h^1$ & $h^2$ & $h^3$ \\ \hline
			$\Omega^3_{X_{t = 1}}$ & 0 & 0 & 0 & 1 \\ \hline
			$\Omega^2_{X_{t = 1}}$ & 0 & 2 & 1 & 0 \\ \hline
			$\Omega^1_{X_{t = 1}}$ & 0 & 1 & 2 & 0 \\ \hline
			$\cO_{X_{t = 1}}$ & 1 & 0 & 0 & 0 \\
		\end{tabular}
	\end{center}
	for $t = 0$ and $t = 1$, as far as we computed them (using $\kq = \QQ$ for $t = 0$). In particular, $\cA^\bullet_{X/S}$ is not cffc. \qedex
\end{ex}

\bibliography{log-sing-surf-smooth-moduli.bib}\bibliographystyle{plain}

\end{document}